\RequirePackage{fix-cm}

\let\oldvec\vec
\documentclass[smallextended]{./styles/svjour3}       
\let\vec\oldvec
\smartqed  
\makeatletter
\def\cl@chapter{\@elt {theorem}}
\makeatother

\usepackage{graphicx}
\graphicspath{{./figs/pdf/}}
\DeclareGraphicsExtensions{.pdf}

\usepackage{xcolor}

\usepackage[draft]{todonotes}   

\usepackage[caption=false,font=footnotesize]{subfig}
\usepackage{lineno}
\usepackage{picins}  
\usepackage{subfig}

\usepackage{multirow}

\usepackage{algorithm}
\usepackage{algpseudocode}
\usepackage{commath} 

\usepackage[capitalize,sort&compress]{cleveref}

\usepackage{amssymb}

\usepackage{algorithm}
\usepackage{algpseudocode}
\algtext*{EndWhile}
\algtext*{EndIf}
\algtext*{EndFor}
\algtext*{EndProcedure}

\usepackage{xspace}

\newcommand{\Tij}{T_\mathrm{ij}}
\newcommand{\Dij}{D_\mathrm{ij}}
\newcommand{\Fij}{F_\mathrm{ij}}
\newcommand{\Ti}{T_\mathrm{i}}
\newcommand{\Tj}{T_\mathrm{j}}
\newcommand{\Ttilde}{\widetilde{T}_\mathrm{i}}
\newcommand{\Tjtilde}{\widetilde{T}_\mathrm{j}}

\newcommand{\Fi}{F_\mathrm{i}}
\newcommand{\Fx}{F}
\newcommand{\Tx}{T_\mathrm{x}}
\newcommand{\Ty}{T_\mathrm{y}}

\newcommand{\Td}{T_\mathrm{d}}
\newcommand{\Deltax}{\Delta_\mathrm{x}}
\newcommand{\Deltay}{\Delta_\mathrm{y}}
\newcommand{\Deltaz}{\Delta_\mathrm{z}}

\newcommand{\x}{\mathbf{x}}
\newcommand{\xii}{\x_\mathrm{i}}
\newcommand{\xjj}{\x_\mathrm{j}}
\newcommand{\xiij}{\x_\mathrm{ij}}
\newcommand{\xmin}{\x_{\mathrm{min}}\xspace}
\newcommand{\Xstart}{\mathcal X_{\mathrm{s}}\xspace}
\newcommand{\Xgoal}{\mathcal X_{\mathrm{g}}\xspace}
\newcommand{\grid}{\mathcal X\xspace}
\newcommand{\F}{\mathcal F\xspace}
\newcommand{\T}{\mathcal T\xspace}
\newcommand{\Narrow}{\texttt{Narrow}\xspace}
\newcommand{\Unknown}{\texttt{Unknown}\xspace}
\newcommand{\Frozen}{\texttt{Frozen}\xspace}
\newcommand{\Octe}{$\mathcal{O}(1)$\xspace}
\newcommand{\Ologn}{$\mathcal{O}(\log{n})$\xspace}
\newcommand{\On}{$\mathcal{O}(n)$\xspace}
\newcommand{\Onlogn}{$\mathcal{O}(n\log{n})$\xspace}

\journalname{Journal}

\begin{document}

\title{Fast Methods for Eikonal Equations: an Experimental Survey
\thanks{This work was supported by the Spanish Ministry of Science and Innovation under project DPI2010-17772 and by the Comunidad Autonoma de Madrid and Structural Funds of the EU under project S2013/MIT-2748.}
}

\author{Javier V. G\'omez\and
        David \'Alvarez \and
        Santiago Garrido \and
        Luis Moreno
}

\institute{ \at RoboticsLab. 
              Carlos III University of Madrid Avda. de la Universidad 30, 28911, Leganes, Madrid, Spain. \\
              Tel.: +34-916248812\\
              \email{{jvgomez,dasanche,sgarrido,moreno}@ing.uc3m.es}           
}

\date{Received: date / Accepted: date}
\maketitle

\begin{abstract}
The Fast Marching Method is a very popular algorithm to compute times-of-arrival maps (distances map measured in time units). Since their proposal in 1995, it has been applied to many different applications such as robotics, medical computer vision, fluid simulation, etc. Many alternatives have been proposed with two main objectives: to reduce its computational time and to improve its accuracy. In this paper, we collect the main approaches which improve the computational time of the standard Fast Marching Method, focusing on single-threaded methods and isotropic environments. 9 different methods are studied under a common mathematical framework and experimentally in representative environments: Fast Marching Method with binary heap, Fast Marching Method with Fibonacci Heap, Simplified Fast Marching Method, Untidy Fast Marching Method, Fast Iterative Method, Group Marching Method, Fast Sweeping Method, Lock Sweeping Method and Double Dynamic Queue Method.

\end{abstract}

\section{Introduction}
The Fast Marching Method and its derivatives (Fast Methods in short) have been extensively applied since they were firstly proposed in 1995 \cite{Tsitsiklis95} as a solution to isotropic control problems using first-order semi-Langragian discretizations on Cartesians grids. Their main field of application are robotics \cite{Gomez12Thesis,Do14,Liu15} and computer vision \cite{Forcadel08}, mainly medical image segmentation \cite{Basu14,Alzaben15}. However, it has proven to be useful in many other applications such as tomography \cite{Xinxin14} or seismology \cite{Zhang11}.

The first approach was proposed by Tsitsiklis \cite{Tsitsiklis95} but the most popular solution was given few months later by Sethian \cite{Sethian96} using first-order upwind-finite differences in the context of isotropic front propagation. Differences and similarities between both works can be found in \cite{Sethian03}.

Fast Methods were originally proposed to simulate a wavefront propagation through a regular discretization of the space. However, many different approaches have been proposed, extending these methods to other discretizations and formulations. For a more detailed history of Fast Marching methods, we refer the interested readers to \cite{Chacon14}. Despite the vast amount of Fast Methods literature, there is a lack of in-depth comparison and benchmarking. 

\begin{figure}[ht]
	\centering
	\includegraphics[width=0.5\columnwidth]{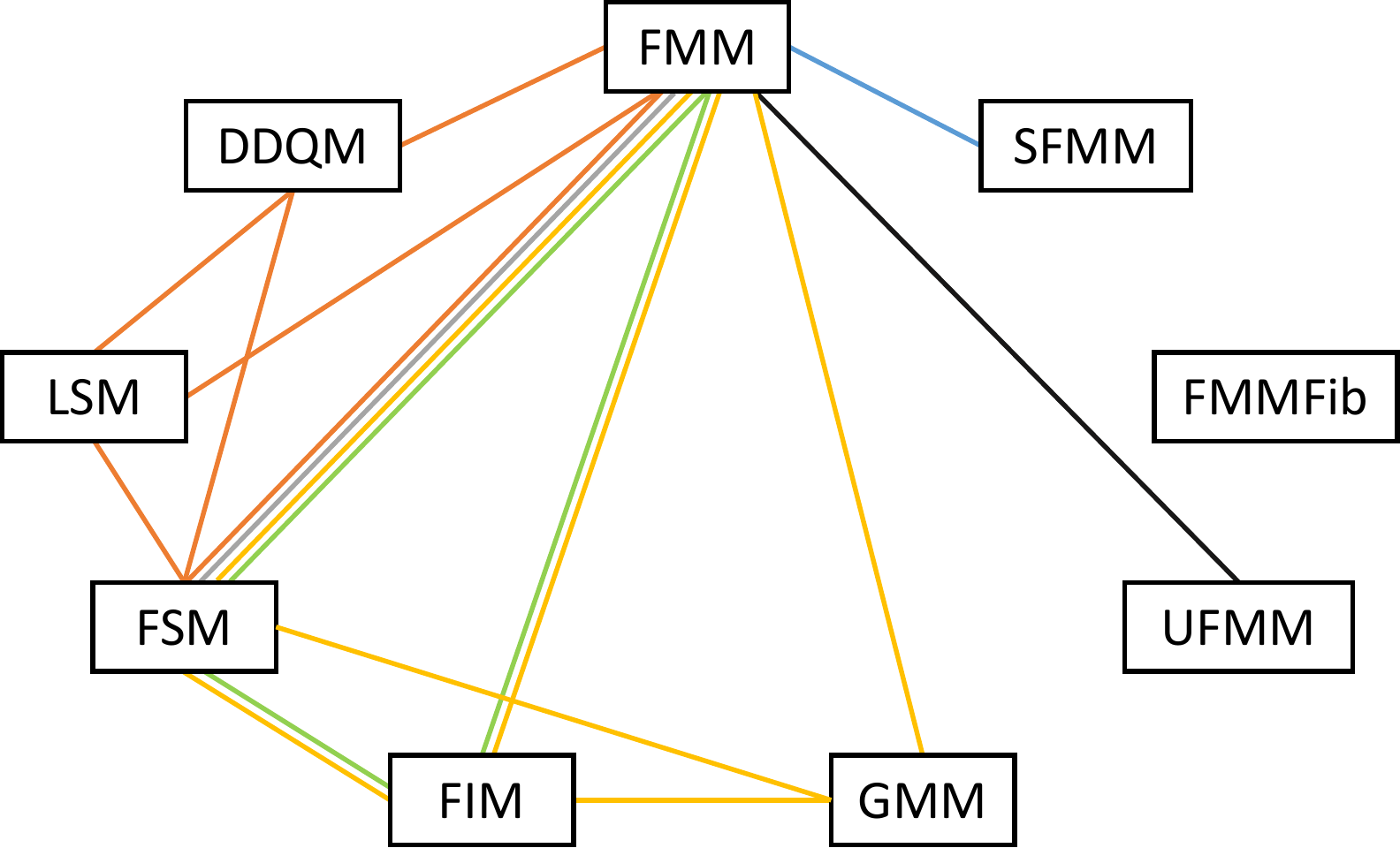}
	\caption{Comparisons among algorithms. Colors refer to different works: orange \cite{Bak10}, gray \cite{Gremaud06}, yellow \cite{Jeong08}, green \cite{Capozzoli13}, black \cite{Yatziv05}, and blue \cite{Jones06}.}
\label{fig:comp}
\end{figure}

In this work 9 sequential (mono-thread), isotropic, grid-based Fast Methods are detailed in the following sections: Fast Marching Method (FMM), Fibonacci-Heap FMM (FMMFib), Simplified FMM (SFMM), Untidy FMM (UFMM), Group Marching Method (GMM), Fast Iterative Method (FIM), Fast Sweeping Method (FSM), Lock Sweeping Method (LSM) and Double Dynamic Queue Method (DDQM). All these algorithms provide exactly the same solution except for UFMM and FIM, which have bounded errors. However, the question of which one is the best for which applications is still open. For example, \cite{Gremaud06} compares only FMM and FSM in spite of the fact that GMM and UFMM were already published. This survey \cite{Jones06} mentions most of the algorithms but only compares FMM and SFMM. A more recent work compares FMM, FSM and FIM in 2D \cite{Capozzoli13}. However, FIM is paralellized and implemented in CUDA providing a biased comparison. \cref{fig:comp} schematically shows the comparisons among algorithms carried out in the literature. Methods such as UFMM have been barely compared to their counterparts while others as FMM and FSM are compared in many works despite the fact that it is well known when each perform better: FSM is faster in simple environments with constant velocity. Also, results from one work cannot be directly extrapolated to other works since the performance of these methods highly depend on their implementation.

\textbf{Statement of contributions:} Three main contributions are included in this paper: 1) Based on previous works, a common formulation and notation is given for all the algorithms. This way, it is possible to easily understand their working principles and mathematical formulation. More detailed formulations are available in the literature but we decided to focus from a practical perspective. 2) A recent survey of the work on designing sequential Fast Methods during the last years. We are not taking into account parallel or high-accuracy approaches as these fields are dense enough to fill another survey. 3) Extensive and systemic comparison among the mentioned methods with experiments designed taking into account their applications and results previously reported.

It is important to note that we have not taken into account the three methods proposed in \cite{Chacon14} (Fast Marching-Sweeping Method, Heap Cell Method and Fast Heap Cell Method) for 3 different reasons: 1) their error is not mathematically bounded, 2) the error highly depends on the environment, and 3) they assume that the velocities are approximately the same in regions of arbitrary size, which is a strong assumption.

This paper is organized as follows: next section introduces the common notation and formulation to be used throughout the document. Following, Fast Methods are detailed, classified attending to their algorithm family. Thus, \cref{sec:fmm} includes Fast Marching-like algorithms, \cref{sec:fsm} Fast Sweeping-based, and \cref{sec:other} contains other Fast Methods. The benchmark and its results are included in section \cref{sec:expcomp}, followed by a discussion in section \cref{sec:disc}. Finally, \cref{sec:conc} outlines the conclusions and proposed future works.

\section{Problem Formulation}
\label{sec:problem}
Fast Methods are built to solve the nonlinear boundary value problem\footnote{This problem formulation closely follows \cite{Sethian99b}}. That is, given a domain $\Omega$ and a function $F:\Omega \rightarrow \mathbb{R}_+$ which represents the local speed of the motion, drive a system from a starting set $\Xstart \subset \Omega$ to a goal set $\Xgoal \subset \delta \Omega$ through the fastest possible path. The Eikonal equation computes the minimum time-of-arrival function $T(\x)$ as follows:
\begin{equation}
\begin{array}{c}
\abs{\nabla T(\x)}F(\x) = 1, \text{on}~\Omega \subset \mathbb{R}^N\\
T(\x) = 0,  \x~\text{in}~\Xstart
\end{array}
\label{eq:eikonal}
\end{equation}

Once solved, $T(\x)$ represents a distances (times-of-arrival) field containing the time it takes to go from any point $\x$ to the closest point in $\Xstart$ following the velocities on $F(\x)$.

We assume, without loss of generality, that the domain is a unit hypercube of $N$ dimensions: $\Omega = [0,1]^N$. The domain is represented with a rectangular Cartesian grid $\grid \subset \mathbb{R}^N$, containing the discretizations of the functions $F(\x)$ and $T(\x)$, $\F$ and $\T$ respectively. We refer to grid points $\xiij = (x_i, y_i), \xiij \in \grid$ as the point $\x = (x,y)$ in the space corresponding to a cell $(i,j)$ of the grid (for the 2D case). For notation simplicity, we will denote $\Tij = T(\xiij) \approx T(\x), \Tij \in \T$, that is, $\Tij$ represents an approximation to the real value of the function $T(\x)$. Analogously, $\Fij = F(\xiij) \approx F(\x), \Fij \in \F$. For a general grid of $N$-dimensions, we will refer to cells by their index (or key) $i$ as $\xii$, since a flat representation is more efficient for such datastructure. We will denote the set of Von-Neumann (4-connectivity in 2D) neighbors of grid point $\xiij$ as $\mathcal{N}(\xiij)$.

\subsection{n-Dimensional Discrete Eikonal Equation}
\label{sec:problem:eikonal}
In this section the most common first-order discretization of the Eikonal equation is detailed. There exist many other first-order and higher-order approaches on grids, meshes and manifolds \cite{Sethian98,Sethian00,Ahmed11,Luo14}. 

We are deriving the discrete Eikonal equation in 2D for better understanding. The most common first-order dicretization is given in \cite{Osher88}, which uses an upwind-difference scheme to approximate partial derivatives of $T(\x)$ ($\Dij^{\pm x}$ represents the one-sided partial difference operator in direction $\pm x$):

\begin{equation}
\begin{array}{l}
	T_x(\x) \approx \Dij^{\pm x}T=\frac{T_{i\pm1,j},\Tij}{\pm\Deltax}\\\\
	T_y(\x) \approx \Dij^{\pm y}T=\frac{T_{i,j\pm1},\Tij}{\pm\Deltay}\\
	\label{eq:disc}
\end{array}
\end{equation}

\begin{equation}
\left\{ \begin{array}{c}
	\max(\Dij^{-x}T,0)^{2}+\min(\Dij^{+x}T,0)^{2}+\\
	\max(\Dij^{-y}T,0)^{2}+\min(\Dij^{+y}T,0)^{2}
\end{array}\right\} = \frac{1}{\Fij^2}
\label{eq:eikonal_disc}
\end{equation}

\noindent Simpler but less accurate solution to \cref{eq:eikonal_disc} is proposed in \cite{Sethian99}:
\begin{equation}
\left\{ \begin{array}{c}
	\max(\Dij^{-x}T,-\Dij^{+x}T,0)^{2}+\\
	\max(\Dij^{-y}T,-\Dij^{+y}T,0)^{2}
\end{array}\right\} = \frac{1}{\Fij^2}
\label{eq:eikonal_disc_simpler}
\end{equation}

\noindent and $\Delta x$ and $\Delta y$ are the grid spacing in the
$x$ and $y$ directions. Substituting \cref{eq:disc} in \cref{eq:eikonal_disc_simpler}
and letting
\begin{equation}
\begin{array}{l}
	T=T_{i,j}\\
	\Tx=\min(T_{i-1,j},T_{i+1,j})\\
	\Ty=\min(T_{i,j-1},T_{i,j+1})
\end{array}
\end{equation}

\noindent we can rewrite the Eikonal Equation, for a discrete 2D space as:

\begin{equation}
\max\left(\frac{T-\Tx}{\Deltax},0\right)^2+\max\left(\frac{T-\Ty}{\Deltay},0\right)^2=\frac{1}{\Fij^2}
\label{eq:eikonal_final2d}
\end{equation}

Since we are assuming that the speed of the front is positive ($F>0$), $T$ must be greater than $\Tx$ and $\Ty$ whenever the front wave has not already passed over the coordinates $(i,j)$. Therefore, it is safe to simplify \cref{eq:eikonal_final2d} to:

\begin{equation}
\left(\frac{T-\Tx}{\Delta x}\right)^{2}+\left(\frac{T-\Ty}{\Delta y}\right)^{2}=\frac{1}{\Fij^2}
\label{eq:eikonal_quad2d}
\end{equation}

\cref{eq:eikonal_quad2d} is a regular quadratic equation of the form $aT^2+bT+c=0$, where:
\begin{equation}
\begin{array}{l}
	a = \Deltax^2+\Deltay^2\\
	b = -2(\Deltay^2\Tx+\Deltax^2\Ty)\\
	c = \Deltay^2\Tx^2+\Deltax^2\Ty^2-\frac{\Deltax^2\Deltay^2}{\Fij^2}
\end{array}
\end{equation}

In order to simplify the notation for the n-dimensional case, let us assume that the grid is composed by cubic cells, that is, $\Deltax = \Deltay = \Deltaz = \dots = h$. Let us denote $\Td$ as the generalization of $\Tx$ or $\Ty$ for dimension $d$, up to $N$ dimensions. We also denote as $\Fx$ the propagation velocity for point with coordinates $(i,j,k,\dots)$. Therefore, the discretization of the Eikonal is a quadratic equation with parameters:

\begin{equation}
\begin{array}{l}
	a = N\\
	b = -2\sum\limits_{d=1}^N\Td\\
	c = (\sum\limits_{d=1}^N\Td^2) - \frac{h^2}{\Fx^2}
\end{array}
\label{eq:ndeikparams}
\end{equation}

\subsection{Solving the nD discrete Eikonal equation}
\label{sec:problem:solvingEikonal}
Wavefront propagation follows causality. That is, in order to reach a point with higher time of arrival, it should firstly travel through neighbors of such point with smaller values. The opposite would imply a jump in time continuity and therefore the solutions would be erroneous.

The proposed Eikonal solution (quadratic equation with parameters of \cref{eq:ndeikparams}) does not guarantee the causality of the resulting distance map, as $\Fx$ and $h$ can have arbitrary values. Therefore, before accepting a solution as valid its causality has to be checked. For instance, in 2D the Eikonal is solved as:

\begin{equation}
T = \frac{\Tx+Ty}{2}+\frac{1}{2}\sqrt{\frac{2h^2}{\Fx^2}-\left(\Tx-\Ty\right)^2}
\label{eq:2deikonalsolution}
\end{equation}

\noindent called the \emph{two-sided} update, as both parents $\Tx$ and $\Ty$ are taken into account. The solution is only accepted if $T \geq \max\left(\Tx,\Ty\right)$. The \emph{upwind condition} \cite{Chacon14} shows that:

\begin{equation}
T \geq \max\left(\Tx,\Ty\right) \Longleftrightarrow \abs{\Tx-\Ty} \leq \frac{h}{\Fx}
\label{eq:upwindcondition}
\end{equation}

If this condition fails, the \emph{one-sided update} is applied instead:

\begin{equation}
T = \min\left(\Tx,\Ty\right)+\frac{h}{\Fx}
\label{eq:onesidedeikonal}
\end{equation}

This is a top-down approach: the parents are iteratively discarded until a causal solution is found. To generalize \cref{eq:upwindcondition} is complex. Therefore, we choose to use a bottom-up approach: \cref{eq:onesidedeikonal} is solved and parents are iteratively included until the time of the next parent is higher than the current solution: $T_k > T$. The procedure is detailed in \cref{alg:eikonal,alg:solvendims}. The $\Call{MinTDim}$ function returns the minimum time of the neighbors in a given dimension (left and right for $dim = 1$, bottom and top for $dim=2$, etc.). Our experiments found this approach more robust for 3 or more dimensions with negligible impact on the computational performance.

\begin{algorithm}
\begin{algorithmic}[1]
\Procedure{SolveEikonal}{$\xii, \T, \F$}
	\State $a \leftarrow N$
	\For{$dim = 1:N$}
		\State $min_T \leftarrow \Call{MinTDim}{dim}$
		\If{$min_T \neq \infty~\textbf{and}~min_T < \Ti$}
			\State $T_{\texttt{values}}{\texttt{.push}}(min_T)$
		\Else
			\State $a \leftarrow a-1$
		\EndIf
	\EndFor
	\If{$a = 0$} \Comment{FSM can cause this situation.}
		\State \textbf{return} $\infty$
	\EndIf
	\State $\T_{\texttt{values}} \leftarrow \Call{Sort}{T_{\texttt{values}}}$
	\For{$dim = 1:a$}
		\State $\Ttilde \leftarrow \Call{SolveNDims}{\xii, dim, T_{\texttt{values}}, \F}$
		\If{$dim = a~\textbf{or}~\Ttilde < \T_{\texttt{values},dim+1}$}
			\State \textbf{break}
		\EndIf
	\EndFor
	\State \textbf{return} $\Ttilde$
\EndProcedure
\end{algorithmic}
\caption{Solve Eikonal Equation}
\label{alg:eikonal}
\end{algorithm}

\begin{algorithm}
\begin{algorithmic}[1]
\Procedure{SolveNDims}{$\xii, dim, T_{\texttt{values}}, \F$}
	\If{$dim = 1$}
		\State \textbf{return} $T_{\texttt{values},1}+\frac{h}{\Fi}$
	\EndIf
	\State $sumT \leftarrow \sum\limits_{i=1}^{dim} T_{\texttt{values},i}$
	\State $sumT^2 \leftarrow \sum\limits_{i=1}^{dim} T^2_{\texttt{values},i}$
	\State $a \leftarrow dim$
	\State $b \leftarrow -2sumT$
	\State $c \leftarrow sumT^2 - \frac{h^2}{\Fi}$
	\State $q \leftarrow b^2-4ac$
	\If{$q < 0$} \Comment{Non-causal solution}
		\State \textbf{return} $\infty$
	\Else
		\State \textbf{return} $\frac{-b+sqrt(q)}{2a}$
	\EndIf
\EndProcedure
\end{algorithmic}
\caption{Solve Eikonal for $n$ dimensions}
\label{alg:solvendims}
\end{algorithm}

\section{Fast Marching Methods}
\label{sec:fmm}
The Fast Marching Method (FMM) \cite{Sethian96} is the most common Eikonal solver. It can be classified as a label-setting, Dijkstra-like algorithm \cite{Dijkstra}. It uses a first-order upwind-finite difference scheme to simulate an isotropic front propagation. The main difference with Dijsktra's algorithm is the operation carried out on every node. Dijkstra's algorithm is designed to work on graphs. Therefore, the value for every node $\xii$ only depends on one parent $\xjj$, following the Bellman's optimality principle \cite{Bellman}:

\begin{equation}
\Ti = \min_{\xii\in\mathcal{N}(\xii)}(c_{ij}+\Tj)
\label{eq:bellman}
\end{equation}

In other words, a node $\xii$ is connected to the parent $\xjj$ in its neighborhood $\mathcal{N}(\xii)$ which minimizes (or maximizes) the function value (in this case $\Ti$) composed by the value of $\Tj$ plus the addition of the cost of \emph{traveling} from $\xjj$ to $\xii$, represented as $c_{ij}$.

The FMM follows Bellman's optimality principle but the value for every node is computed following first-order upwind discretization of the Eikonal equation detailed in \cref{sec:problem}. This discretization takes into account the spatial representation (i.e. a rectangular grid) and the value of all the causal upwind neighbors. Thus, the times-of-arrival field computed by FMM is more accurate than Dijkstra's.

The algorithm divides the cells in three different sets: 1) $\Frozen$: those cells which value is computed and cannot change, 2) $\Unknown$: cells with no value assigned, to be evaluated, and 3) $\Narrow$ band (or just $\Narrow$): frontier between $\Frozen$ and $\Unknown$ containing those cells with a value assigned that can be improved. These sets are mutually exclusive, that is, a cell cannot belong to more than one of them at the same time. The implementation of the $\Narrow$ set is a critical aspect of FMM. A more detailed discussion will be carried out in \cref{sec:fmm:heap}.

The procedure is detailed in \cref{alg:fmm}. Initially, all points\footnote{From now on, we will indistinctly use point, cell or node to refer to each element of the grid.} in the grid belong to the Unknown set with infinite arrival time. The initial points (wave sources) are assigned a value 0 and introduced in $\Frozen$ (lines \ref{fmm:init:init}-\ref{fmm:init:end}). Then, the main FMM loop starts by choosing the element with minimum arrival time from $\Narrow$ (line \ref{fmm:min}). All its non-$\Frozen$ neighbors are evaluated: for each of them the Eikonal is solved and the new arrival time value is kept if it is improved. In case the cell is in $\Unknown$, it is transferred to $\Narrow$ (lines \ref{fmm:prop:init}-\ref{fmm:prop:end}). Finally, the previously chosen point from $\Narrow$ is transferred to $\Frozen$ (lines \ref{fmm:end1} and \ref{fmm:end2}) and a new iteration starts until the $\Narrow$ set is empty. The arrival times map $\T$ is returned as the result of the procedure.

\begin{algorithm}
\begin{algorithmic}[1]
\Procedure{FMM}{$\grid, \T, \F, \Xstart$}
	\Statex Initialization:
	\State $\Unknown \leftarrow \grid, \Narrow \leftarrow 
 \emptyset, \Frozen \leftarrow \emptyset$ \label{fmm:init:init}
	\State $\Ti \leftarrow \infty~\forall \xii \in \grid$
	\For{$\xii \in \Xstart$}
		\State $\Ti \leftarrow 0$
		\State $\Unknown \leftarrow \Unknown \backslash \{\xii\}$
		\State $\Narrow \leftarrow \Narrow \cup \{\xii\}$ \label{fmm:init:end}
	\EndFor
	\Statex
	\Statex Propagation:
	\While{$\Narrow \neq \emptyset$}
		\State $\xmin \leftarrow \arg\min_{\xii\in\Narrow}{\{\Ti\}}$ \Comment{$\Narrow$ top operation.}\label{fmm:min}
		\For{$\xii \in (\mathcal{N}(\xmin) \cap \grid\backslash\Frozen) $} \Comment{For all neighbors not in Frozen.} \label{fmm:prop:init}
			\State $\Ttilde \leftarrow \Call{SolveEikonal}{\xii, \T, \F}$
			\If{$\Ttilde < \Ti$} 
				\State $\Ti \leftarrow \Ttilde$ \Comment{$\Narrow$ increase operation if $\xii\in\Narrow$.}
			\EndIf
			\If{$\xii \in \Unknown$} \Comment{$\Narrow$ push operation.}
				\State $\Narrow \leftarrow \Narrow \cup \{\xii\}$
				\State $\Unknown \leftarrow \Unknown \backslash \{\xii\}$ \label{fmm:prop:end}
			\EndIf
		\EndFor	
		\State $\Narrow \leftarrow \Narrow \backslash \{\xmin\}$ \label{fmm:end1}\Comment{$\Narrow$ pop operation: add to Frozen.}
		\State $\Frozen \leftarrow \Frozen \cup \{\xmin\}$ \label{fmm:end2}
	\EndWhile
	\State \textbf{return} $\T$
\EndProcedure
\end{algorithmic}
\caption{Fast Marching Method}
\label{alg:fmm}
\end{algorithm}

\subsection{Binary and Fibonacci Heaps}
\label{sec:fmm:heap}
FMM requires the implementation of the $\Narrow$ set to have four different operations: 1) Push: to insert a new element to the set, 2) Increase: to reorder an element already existing in the set which value has been improved, 3) Top: retrieve the element with minimum value, and 4) Pop: remove the element with minimum value. As stated before, this is the most critical aspect of the FMM implementation. The most efficient way to implement $\Narrow$ is by using a min-heap data structure. A heap is an ordered tree in which every parent is ordered with respect to its children. In a min-heap, the minimum value is at the root of the tree and the children have higher values. This is satisfied for any parent node of the tree.

Among all the existing heaps, FMM is usually implemented with a binary heap \cite{Porter74}. However, the Fibonacci Heap \cite{Fredman87} has a better amortized time for Increase and Push operations. However, it has additional computational overhead with respect other heaps. For relatively small grids, where the narrow band is composed by few elements and the performance is still far from its asymptotic behavior, the binary heap performs better. \cref{tab:heaps} summarizes the time complexities for these heaps\footnote{http://bigocheatsheet.com/} (the priority queue will be detailed in \cref{sec:fmm:sfmm}). Note that $n$ for the methods complexity is the number of cells in the map, as the worst case is to have all the cells in the heap.

\begin{table}[ht]
\caption{Summary of amortized time complexities for common heaps used in FMM ($n$ is the number of elements in the heap). }
\begin{tabular}{r|cccc}
\multicolumn{1}{l|}{} & Push      & Increase  & Top  	& Pop       \\ \hline
Fibonacci             & \Octe     & \Octe			& \Octe	& \Ologn \\
Binary                & \Ologn 		& \Ologn	  & \Octe	& \Ologn \\
Priority Queue        & \Ologn 		& --    		& \Octe & \Ologn \\ \hline
\end{tabular}
\label{tab:heaps}
\end{table}

Each cell is pushed and popped at most once in the heap. For each loop, the top of $\Narrow$ is accessed (\Octe), the Eikonal is solved for at most $2^N$ neighbors (\Octe~for a given $N$), these cells are pushed or increased (\Ologn~in the worst case), and finally the top cell is popped (\Ologn). Therefore each loop is at most \Ologn. Since this loop is executed at most $n$ times, the total FMM complexity is \Onlogn, where $n$ represents the total number of cells of the grid in the worst case scenario.

\subsection{Simplified Fast Marching Method}
\label{sec:fmm:sfmm}
The Simplified Fast Marching Method (SFMM) \cite{Jones06} is a relatively unknown variation of the standard FMM but with an impressive performance. SFMM is a reduced version of FMM, where $\Narrow$, implemented as a simple priority queue, can contain different instances of the same cell with different values. Additionally, it can happen that the same cell belongs to $\Narrow$ and $\Frozen$ at the same time. The simplification occurs since no Increase operation is required. Every time a cell has an updated value, it is pushed to the priority queue. Once it is popped and inserted in $\Frozen$, the remaining instances in the queue are simply ignored.

The advantage is that all the increase operations are substituted by push operations. Although both have the same computational complexity, the constant for push is much lower (increase requires removal and Push operations). Note that the computational complexity is maintained, \Onlogn.

\begin{algorithm}
\begin{algorithmic}[1]
\Procedure{SFMM}{$\grid, \T, \F, \Xstart$}
	\Statex Initialization as FMM (in \cref{alg:fmm})
	\Statex
	\Statex Propagation:
	\While{$\Narrow \neq \emptyset$}
		\State $\xmin \leftarrow \arg\min_{\xii\in\Narrow}{\{\Ti\}}$\label{alg:sfmm:xmin} \Comment{$\Narrow$ top operation.}
		\If{$\xmin \in \Frozen$}
			\State $\Narrow \leftarrow \Narrow \backslash \{\xmin\}$
		\Else
			\For{$\xii \in (\mathcal{N}(\xmin) \cap \grid\backslash\Frozen) $} \Comment{For all neighbors not in Frozen.}
				\State $\Ttilde \leftarrow \Call{SolveEikonal}{\xii, \T, \F}$
				\If{$\Ttilde < \Ti$} \Comment{Update arrival time.}
					\State $\Ti \leftarrow \Ttilde$
					\State $\Narrow \leftarrow \Narrow \cup \{\xii\}$ \Comment{$\Narrow$ push operation.}
				\EndIf
				\If{$\xii \in \Unknown$}
					\State $\Unknown \leftarrow \Unknown \backslash \{\xii\}$
				\EndIf
			\EndFor	
			\State $\Narrow \leftarrow \Narrow \backslash \{\xmin\}$ \Comment{$\Narrow$ pop operation.}
			\State $\Frozen \leftarrow \Frozen \cup \{\xmin\}$
		\EndIf
	\EndWhile
	\State \textbf{return} $\T$
\EndProcedure
\end{algorithmic}
\caption{Simplified Fast Marching Method}
\label{alg:sfmm}
\end{algorithm}

\subsection{Untidy Fast Marching Method}
\label{sec:fmm:ufmm}
The Untidy Fast Marching Method (UFMM) \cite{Yatziv05,Rasch08} follows exactly the same procedure as FMM. However, a special heap structure is used which reduces the overall computational complexity to \On: the \emph{untidy} priority queue.

This untidy priority queue is closer to a look-up table than to a tree. It assumes that the $\F$ values are bounded, hence the $\T$ values are also bounded. The untidy queue is a circular array which divides the maximum range of $\T$ into a set of $k$ consecutive buckets. Each bucket contains an unordered list of cells with similar $\Ti$ value. The threshold values for each bucket evolve with the algorithm, trying to maintain an uniform distribution of the elements in $\Narrow$ among the buckets.

Since the index of the corresponding bucket can be analytically computed, Push is \Octe~as well as Top and Pop. The Increase operation is, in average, \Octe~as long as $\#_{buckets} <$\On. Therefore, the total UFMM complexity is \On. However, since elements within a bucket are not sorted (FIFO strategy is applied in each bucket), errors are being introduced in the final result. In fact, it is shown that the accumulated additional error is bounded by $\mathcal{O}(h)$, which is the same order of magnitude as in the original FMM.

\begin{figure}[ht]
	\centering
	\includegraphics[width=\columnwidth]{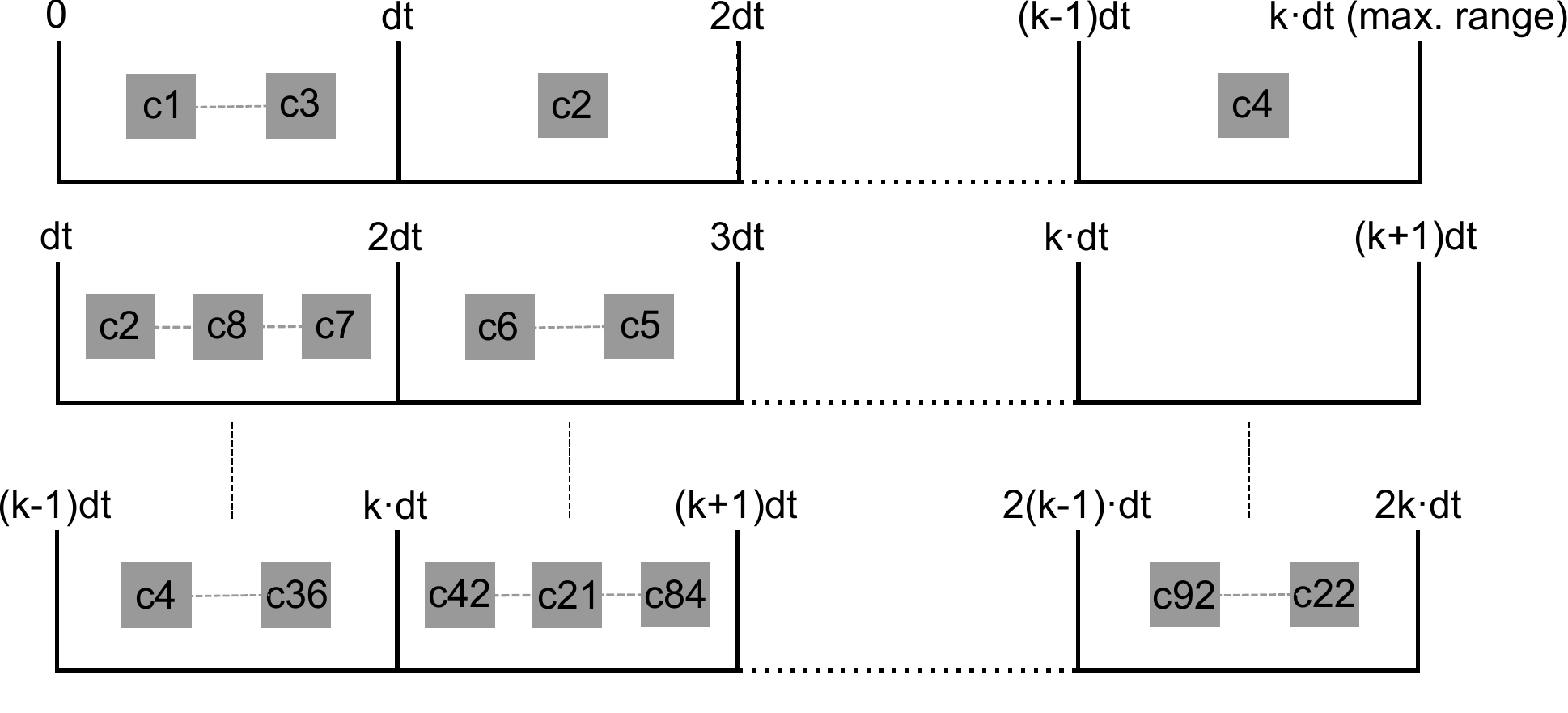}
	\caption{Untidy priority queue representation. Top: first iteration, the four neighbors of the initial point are pushed. Middle: the first bucket became empty, so the circular array advances one position. Cell c2 is first evaluated because it was the first pushed in the bucket. Bottom: after a few iterations, a complete loop on the queue is about to be completed.}
\label{fig:ufmm}
\end{figure}

\section{Fast Sweeping Methods}
\label{sec:fsm}
\begin{figure}[ht]
	\centering
	\includegraphics[width=\columnwidth]{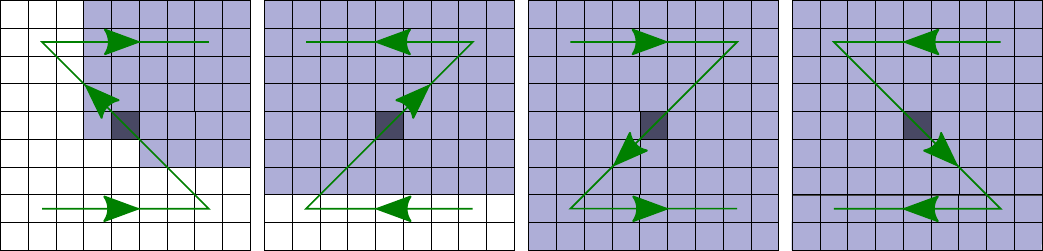}
	\caption{FSM sweep directions in 2D represented with arrows. The darkest cell is the initial point and the shaded cells are those analyzed by the current sweep (time improved or maintained).}
\label{fig:fsm:sweeps}
\end{figure}

The Fast Sweeping Method (FSM) \cite{Tsai03,Zhao05} is an iterative algorithm which computes the times-of-arrival map by successively \emph{sweeping} (traversal) the whole grid following a specific order. FSM performs Gauss-Seidel iterations in alternating directions. These directions are chosen so that all the possible characteristic curves of the solution to the Eikonal are divided into the possible quadrants (or octants in 3D) of the environment. For instance, a bi-dimensional grid has 4 possible Gauss-Seidel iterations (the combinations of traversing $x$ and $y$ dimensions forwards and backwards): are North-East, North-West, South-East and South-West, as shown in \cref{fig:fsm:sweeps}.

The FSM is a simple algorithm: it performs sweeps until no value is improved. In each sweep, the Eikonal equation is solved for every cell. However, to generalize this algorithm to $N$-dimensions is complex. Up to our knowledge, there are only 2D and 3D implementations. However, in \cref{alg:fsm} we introduce an $N$-dimensional version. We will denote the sweeping directions as a binary array $\texttt{SweepDirs}$ with elements $1$ or $-1$, with $1$ ($-1$) meaning forwards (backwards) traversal in that dimension. This array is initialized to $1$ (North-East in the 2D case or North-East-Top in 3D) and the grid is initialized as in FMM (lines \ref{fsm:init:init}-\ref{fsm:init:end}). The main loop updates $\texttt{SweepDirs}$ and a sweep is performed in the new direction (lines \ref{fsm:getsweep}-\ref{fsm:sweep}).

The $\Call{GetSweepDirs}$ procedure (see \cref{alg:fsm:gsi}) is in charge of generating the appropriate Gauss-Seidel iteration directions. If a 3D $\texttt{SweepDirs} = [1,1,1]$ vector is given, the following sequence will be generated:

\begin{equation}
\begin{split}
&1 - [-1,-1,-1]\\
&2 - [1,-1,-1]\\
&3 - [-1,1,-1]\\
\end{split}
\quad
\begin{split}
&4 - [1,1,-1]\\
&5 - [-1,-1,1]\\
&6 - [1,-1,1]\\
\end{split}
\quad
\begin{split}
&7 - [-1,1,1]\\
&8 - [1,1,1]\\
\end{split}
\end{equation}

Note that this sequence creates a sweep pattern which is not exactly the same as detailed in the literature, but it is equally valid as the same directions are visited and the same number of sweeps are required to cover the whole grid.

Finally, the $\Call{Sweep}$ procedure (see \cref{alg:fsm:sweep}) recursively generates the Gauss-Seidel iterations following the traversal directions specified by the corresponding value of $\texttt{SweepDirs}$ (\cref{fsm:sweep:for}). Once the most inner loop is reached, the corresponding cell is evaluated and its value updated if necessary (lines \ref{fsm:sweep2}-\ref{fsm:sweep3}).

\begin{algorithm}
\begin{algorithmic}[1]
\Procedure{FSM}{$\grid, \T, \F, \Xstart$}
	\Statex Initialization.
	\State $\texttt{SweepDirs} \leftarrow [1,\dots,1]$ \Comment{Initialize sweeping directions.} \label{fsm:init:init}
	\State $\Ti \leftarrow \infty~\forall \xii \in \grid$
	\For{$\xii \in \Xstart$}
		\State $\Ti \leftarrow 0$ \label{fsm:init:end}
	\EndFor
	\Statex
	\Statex Propagation:
	\State $\texttt{stop} \leftarrow \texttt{False}$
	\While{$\texttt{stop} \neq \texttt{True}$}
		\State $\texttt{SweepDirs} \leftarrow \Call{getSweepDirs}{\grid, \texttt{SweepDirs}}$ \label{fsm:getsweep}
		\State $\texttt{stop} \leftarrow \Call{Sweep}{\grid, \T, \F, \texttt{SweepDirs}, N}$ \label{fsm:sweep}
	\EndWhile
	\State \textbf{return} $\T$
\EndProcedure
\end{algorithmic}
\caption{Fast Sweeping Method}
\label{alg:fsm}
\end{algorithm}

\begin{algorithm}
\begin{algorithmic}[1]
\Procedure{getSweepDirs}{$\grid, \texttt{SweepDirs}$}
	\For{$i = 1:N$}
		\State $\texttt{SweepDirs}_\mathrm{i} \leftarrow \texttt{SweepDirs}_\mathrm{i} + 2$
		\If{$\texttt{SweepDirs}_\mathrm{i} \leq 1$}
			\State \textbf{break} \Comment{Finish For loop.}
		\Else
			\State $\texttt{SweepDirs}_\mathrm{i} \leftarrow -1$
		\EndIf
	\EndFor
	\State \textbf{return} \texttt{SweepDirs}
\EndProcedure
\end{algorithmic}
\caption{Sweep directions algorithm}
\label{alg:fsm:gsi}
\end{algorithm}

\begin{algorithm}
\begin{algorithmic}[1]
\Procedure{Sweep}{$\grid, \T, \F, \texttt{SweepDirs}, n$}
	\State $\texttt{stop} \leftarrow \texttt{True}$
	\If{$n > 1$}
		\For{$i \in \grid_\mathrm{n}~\text{following}~\texttt{SweepDirs}_\mathrm{n}$} \label{fsm:sweep:for}
			\State $\texttt{stop} \leftarrow \Call{Sweep}{\grid, \T, \F, \texttt{SweepDirs}, n-1}$
		\EndFor
	\Else
		\For{$i \in \grid_\mathrm{1}~\text{following}~\texttt{SweepDirs}_\mathrm{1}$} \label{fsm:sweep2}
			\State $\Ttilde \leftarrow \Call{SolveEikonal}{\xii, \T, \F}$ \Comment{$\xii$ is the corresponding cell.}
			\If{$\Ttilde < \Ti$}
				\State $\Ti \leftarrow \Ttilde$
				\State $\texttt{stop} \leftarrow \texttt{False}$ \label{fsm:sweep3}
			\EndIf
		\EndFor
	\EndIf
	\State \textbf{return} \texttt{stop}
\EndProcedure
\end{algorithmic}
\caption{Recursive sweeping algorithm}
\label{alg:fsm:sweep}
\end{algorithm}

The FSM carries out as many grid traversals as necessary until the value $\Ti$ for every cell has converged. Since no ordering is implied, the evaluation of each cell is \Octe. As there are $n$ cells, the total computational complexity of FSM is \On. However, note that the constants highly depend on the velocity function $F(\x)$. In the case of an empty environment with constant $F(\x)$, only $2^N$ sweeps will be required as the characteristic directions are straight lines. However, for environment with obstacles or complex velocities functions, where the characteristics directions change frequently, the number of sweeps required can be much higher and therefore FSM will take longer to return a solution. Note as well that the $\T$ returned by FSM is exactly the same as all the FMM-like algorithms (except UFMM).

\subsection{Lock Sweeping Methods}
\label{sec:fsm:lsm}
The Lock Sweeping Method (LSM) \cite{Bak10} is a natural improvement over FSM. The FSM might spend computation time recomputing $\Ti$ even if none of the neighbors of $\xii$ has improved their value since the last sweep. LSM labels a cell as \emph{unlocked} if any of its neighbors has changed and thus its value can be improved. Otherwise, the cell is labeled as \emph{locked} and it will be skipped.

The LSM procedure is detailed in \cref{alg:lsm}. It is basically the same as FSM but with the addition of tracking if a point is locked (\Frozen) or unlocked (\Narrow). Analogously, the $\Call{LockSweep}$ (see \cref{alg:lsm:sweep}) procedure is similar to $\Call{Sweep}$ with two differences: 1) if a point is not unlocked it is skipped (see \cref{alg:lsm:skip}), and 2) neighbors of $\xii$ are unlocked if the new value $\Ti$ is better that their current value.

\begin{algorithm}
\begin{algorithmic}[1]
\Procedure{LSM}{$\grid, \T, \F, \Xstart$}
	\Statex Initialization.
	\State $\Frozen \leftarrow \grid, \Narrow \leftarrow \emptyset$
	\State $\Ti \leftarrow \infty~\forall \xii \in \grid$
	\State $\texttt{SweepDirs} \leftarrow [1,\dots,1]$ \Comment{Initialize sweeping directions.}
	\For{$\xii \in \Xstart$}
			\State $\Ti \leftarrow 0$
			\For{$\xjj \in \mathcal{N}(\xii)$}\Comment{Unlocking neighbors of starting cells.}
				\State $\Frozen \leftarrow \Frozen \backslash \{\xjj\}$
				\State $\Narrow \leftarrow \Narrow \cup \{\xjj\}$
			\EndFor			
	\EndFor
	\Statex
	\Statex Propagation:
	\State $\texttt{stop} \leftarrow \texttt{False}$
	\While{$\texttt{stop} \neq \texttt{True}$}
		\State $\texttt{SweepDirs} \leftarrow \Call{getSweepDirs}{\grid, \texttt{SweepDirs}}$
		\State $\texttt{stop} \leftarrow \Call{LockSweep}{\grid, \T, \F, \texttt{SweepDirs}, N}$
	\EndWhile
	\State \textbf{return} $\T$
\EndProcedure
\end{algorithmic}
\caption{Lock Sweeping Method}
\label{alg:lsm}
\end{algorithm}

\begin{algorithm}
\begin{algorithmic}[1]
\Procedure{LockSweep}{$\grid, \T, \F, \texttt{SweepDirs}, n$}
	\State $\texttt{stop} \leftarrow \texttt{True}$
	\If{$n > 1$}
		\For{$i \in \grid_\mathrm{n}~\text{following}~\texttt{SweepDirs}_\mathrm{n}$}
			\State $\texttt{stop} \leftarrow \Call{LockSweep}{\grid, \T, \F, \texttt{SweepDirs}, n-1}$
		\EndFor
	\Else
		\For{$i \in \grid_\mathrm{1}~\text{following}~\texttt{SweepDirs}_\mathrm{1}$}
			\If{$\xii \in \Narrow$} \label{alg:lsm:skip}
				\State $\Ttilde \leftarrow \Call{SolveEikonal}{\xii, \T, \F}$ \Comment{$\xii$ is the corresponding cell.}
				\If{$\Ttilde < \Ti$}
					\State $\Ti \leftarrow \Ttilde$
					\State $\texttt{stop} \leftarrow \texttt{False}$
					\For{$\xjj \in \mathcal{N}(\xii)$}
						\If{$\Ti < \Tj$} \Comment{Add improvable neighbors to Narrow.}
							\State $\Frozen \leftarrow \Frozen \backslash \{\xjj\}$
							\State $\Narrow \leftarrow \Narrow \cup \{\xjj\}$
						\EndIf
					\EndFor
				\EndIf
				\State $\Narrow \leftarrow \Narrow \backslash \{\xii\}$ \Comment{Add $\xii$ to Frozen.}
				\State $\Frozen \leftarrow \Frozen \cup \{\xii\}$
			\EndIf
		\EndFor
	\EndIf
	\State \textbf{return} \texttt{stop}
\EndProcedure
\end{algorithmic}
\caption{Recursive sweeping algorithm}
\label{alg:lsm:sweep}
\end{algorithm}

Note that the asymptotic computational complexity of FSM is kept, \On. The number of required sweeps is also maintained. However, in practice it turns out that most of the cells are locked during a sweep. Therefore, the computation time saved is important.

\section{Other Fast Methods}
\label{sec:other}
\subsection{Group Marching Method}
\label{sec:other:gmm}
The Group Marching Method (GMM) \cite{Kim01} is an FMM-based Eikonal solver which solves a group of grid points in $\Narrow$ at once, instead of sorting them in a heap structure.

Consider a front propagating. At a given time, $\Narrow$ will be composed by the set of cells belonging to the wavefront. GMM selects a group $G$ out of $\Narrow$ composed by the global minimum and the local minima in $\Narrow$. Then, every neighboring cell to $G$ is evaluated and added to $\Narrow$. These points in $G$ have to be chosen carefully so that causality is not violated since GMM does not sort the $\Narrow$ set. For that, GMM selects those points following \cref{eq:gmm:select}:

\begin{equation}
G = \{ \xii \in \Narrow : \Ti \leq \min(T_{\Narrow}) + \delta_\tau \}
\label{eq:gmm:select}
\end{equation}

\noindent where

\begin{equation}
\delta_\tau = \frac{1}{\max(\F)}
\label{eq:gmm:delta}
\end{equation}

In the original GMM work \cite{Kim01}, $\delta_\tau = \frac{h}{\max(\F)\sqrt{N}}$. However, we have chosen \cref{eq:gmm:delta} as referred in \cite{Jeong08}. Although this second formula is not mathematically proven, the results for the original $\delta_\tau$ are much worse than FMM in most of the cases, reaching an order of magnitude of difference.

If the time difference between two adjacent cells is larger than $\delta_\tau$, their values will barely affect each other since the wavefront propagation direction is more perpendicular than parallel to the line segment formed by both cells. However, the downwind points (those to be evaluated in future iterations) can be affected by both adjacent cells. Therefore, points in $G$ are evaluated twice to avoid instabilities.

GMM is detailed in \cref{alg:gmm}. Its initialization is FMM-like. Note that $\delta_\tau$ depends on the maximum velocity value in the grid. The main loop updates the threshold $T_\mathrm{m}$ every iteration. Firstly, it carries out a reverse traversal through the selected points, computing and updating their value (lines \ref{gmm:rev_trav:init}-\ref{gmm:rev_trav:end}). Then,  lines \ref{gmm:for_trav:init}-\ref{gmm:for_trav:end} perform a forward traversal with the same operations as the reverse traversal but updating the $\Narrow$ and $\Frozen$ sets in the same way as FMM.

\begin{algorithm}
\begin{algorithmic}[1]
\Procedure{GMM}{$\grid, \T, \F, \Xstart$}
	\Statex Initialization:
	\State $\Unknown \leftarrow \grid, \Narrow \leftarrow
 \emptyset, \Frozen \leftarrow \emptyset$
	\State $\Ti \leftarrow \infty~\forall \xii \in \grid$
	\State $\delta_\tau \leftarrow \frac{1}{\max(\F)}$
	\For{$\xii \in \Xstart$}
		\State $\Ti \leftarrow 0$
		\State $\Unknown \leftarrow \Unknown \backslash \{\xii\}$
		\State $\Frozen \leftarrow \Frozen \cup \{\xii\}$
		\For{$\xjj \in \mathcal{N}(\xii)$} \Comment{Adding neighbors of starting points to Narrow.}
			\State $\Ti \leftarrow \Call{SolveEikonal}{\xjj, \T, \F}$
			\If{$\Ti < T_\mathrm{m}$}
				\State $T_\mathrm{m} \leftarrow \Ti$
			\EndIf
			\State $\Unknown \leftarrow \Unknown \backslash \{\xii\}$
			\State $\Narrow \leftarrow \Narrow \cup \{\xii\}$
		\EndFor
	\EndFor
	\Statex
	\Statex Propagation:
	\While{$\Narrow \neq \emptyset$}
		\State $T_\mathrm{m} \leftarrow T_\mathrm{m} + \delta_\tau$
		\For{$\xii \in (\Narrow \leq T_\mathrm{m})~\texttt{REVERSE}$} \Comment{Reverse traversal.} \label{gmm:rev_trav:init}
			\For{$\xjj \in (\mathcal{N}(\xii)\cap\grid\backslash\Frozen)$}
				\State $\Ttilde \leftarrow \Call{SolveEikonal}{\xjj, \T, \F}$
				\If{$\Ttilde < \Ti$}
					\State $\Ti \leftarrow \Ttilde$  \label{gmm:rev_trav:end}
				\EndIf
			\EndFor
		\EndFor
		\For{$\xii \in (\Narrow \leq T_\mathrm{m})~\texttt{FORWARD}$} \Comment{Forward tranversal.}  \label{gmm:for_trav:init}
			\For{$\xjj \in (\mathcal{N}(\xii)\cap\grid\backslash\Frozen)$}
				\State $\Ttilde \leftarrow \Call{SolveEikonal}{\xjj, \T, \F}$
				\If{$\Ttilde < \Ti$}
					\State $\Ti \leftarrow \Ttilde$
				\EndIf
				\If{$\xii \in \Unknown$}
					\State $\Unknown \leftarrow \Unknown \backslash \{\xii\}$
					\State $\Narrow \leftarrow \Narrow \cup \{\xii\}$
				\EndIf
			\EndFor
			\State $\Narrow \leftarrow \Narrow \backslash \{\xii\}$
			\State $\Frozen \leftarrow \Frozen \cup \{\xii\}$  \label{gmm:for_trav:end}
		\EndFor
	\EndWhile
	\State \textbf{return} $\T$
\EndProcedure
\end{algorithmic}
\caption{Group Marching Method}
\label{alg:gmm}
\end{algorithm}

Note that GMM returns the same solution as FMM. GMM evaluates twice every node before inserting it in $\Frozen$ while FMM only evaluates it once. However, GMM does not require any sorting. Therefore, GMM is an \On~ iterative algorithm that converges in only 2 iterations (traversals). The value of $\delta_\tau$ can be modified: higher $\delta_\tau$ would require more iterations to converge. However, smaller $\delta_\tau$ will require also 2 traversals but the group $G$ will be composed by less cells. As GMM authors point out, GMM can be interpreted as an intermediary point between FMM ($\delta_\tau = 0$) and a purely iterative method \cite{Rouy92} ($\delta_\tau = \infty$).

Additionally, a generalized $N$-dimensional implementation is straightforward.

\subsection{Dynamic Double Queue Method}
\label{sec:other:ddqm}
The Dynamic Double Queue Method (DDQM) \cite{Bak10} is inspired in the LSM but resembles to GMM. 
DDQM is conceptually simple. $\Narrow$ is divided into two non-sorted FIFO queues: one with cells to be evaluated sooner and the other one with cells to be evaluated later. Every iteration takes an element from the first queue and evaluates it. If the time is improved, the neighboring cells with higher time are unlocked and added to the first or second queue depending on the value of the cell updated. Once the first queue is empty, queues are swapped and the algorithm continues. The purpose is to achieve a pseudo-ordering of the cells, so that cells with lower value are evaluated first.

Since queues are not sorted, it could require to solve many times the same cell until its value converges. DDQM dynamically computes the threshold value depending on the number of points inserted in each queue, trying to reach an equilibrium. The original paper includes an important analysis about this threshold update. Initially, the step value of the threshold is increased every iteration and is computed as:

\begin{equation}
step = \frac{1.5hn}{\sum\limits_{i}\Fi}
\label{eq:ddqm:step}
\end{equation}

\noindent where $n$ is the total number of cells in the grid. Originally, this step was proposed as $step = \frac{1.5n}{h\sum\limits_{i}\frac{1}{\Fi}}$. However, step should have time units and this expression $[t^{-1}]$ (probably an error due to the ambiguity of using speed $F$ or slowness $f=\frac{1}{F}$). Therefore, we propose as alternative \cref{eq:ddqm:step}.

Every time the first queue is empty, $\Call{UpdateStep}$ (see \cref{alg:ddqm:step}) is called, with the value of the current $step$, $c_1$, and $c_{\mathrm{total}}$ which are the number of cells inserted in the first queue and the total number of cells, correspondingly. $step$ is modified so that the number of cells inserted in the first queue are between 65$\%$ and 75$\%$. This is a conservative approach, since the closer this percentage is to 50$\%$ the faster DDQM is. However, the penalization for percentages lower than 50$\%$ is much important than for higher percentages. 

Note that the step is increased by a factor 1.5 but decreased by a factor of 2. This makes $step$ to converge to a value instead of overshooting around the optimal value. Dividing by a larger number causes the first queue to become empty earlier. Thus, next iteration will finish faster and a better $step$ value can be computed.

\begin{algorithm}
\begin{algorithmic}[1]
\Procedure{UpdateStep}{$step, c_1, c_{\mathrm{total}}$}
	\State $\texttt{m} \leftarrow 0.65$
	\State $\texttt{M} \leftarrow 0.75$
	\State $Perc \leftarrow 1$
	\If{$c_1 > 0$}
		\State $Perc \leftarrow \frac{c_1}{c_{\mathrm{total}}}$
	\EndIf
	\If{$Perc \leq m$}
		\State $step \leftarrow step*1.5 $
	\ElsIf{$Perc \geq M$}
		\State $step \leftarrow \frac{step}{2}$
	\EndIf
	\State \textbf{return} $step$
\EndProcedure
\end{algorithmic}
\caption{DDQM Threshold Increase}
\label{alg:ddqm:step}
\end{algorithm}

DDQM is detailed in \cref{alg:ddqm}. As in LSM, points are locked (\Frozen) or unlocked (\Narrow). Initialization sets all points as frozen except the neighbors of the start points, which are added to the first queue (lines \ref{ddqm:init:init}-\ref{ddqm:init:end}). While first queue is not empty, its front element is extracted and evaluated (lines \ref{ddqm:prop:1}-\ref{ddqm:prop:2}). If its value is improved, all its locked neighbors with higher value are unlocked and added to its corresponding queue.

\begin{algorithm}
\begin{algorithmic}[1]
\Procedure{DDQM}{$\grid, \T, \F, \Xstart$}
	\Statex Initialization:
	\State $\Frozen \leftarrow \grid, \Narrow \leftarrow \emptyset$ \label{ddqm:init:init}
	\State $\mathcal{Q}_1 \leftarrow \emptyset, \mathcal{Q}_2 \leftarrow \emptyset$
	\State $c_1 \leftarrow 0$ \Comment{Counters}
	\State $c_{\mathrm{total}} \leftarrow 0$
	\State $step = \frac{1.5hn}{\sum\limits_{i}\Fi}$ \Comment{$n$ is the total number of cells.}
	\State $th \leftarrow step$
	\State $\Ti \leftarrow \infty~\forall \xii \in \grid$
	\For{$\xii \in \Xstart$}
		\State $\Ti \leftarrow 0$
		\For{$\xjj \in \mathcal{N}(\xii)$}
			\State $\mathcal{Q}_1 \leftarrow \mathcal{Q}_1 \cup \{\xjj\}$
			\State $\Unknown \leftarrow \Unknown \backslash \{\xii\}$
			\State $\Narrow \leftarrow \Narrow \cup \{\xii\}$ \label{ddqm:init:end}
		\EndFor
	\EndFor
	\Statex
	\Statex Propagation:
	\While{$\mathcal{Q}_1 \neq \emptyset~\textbf{or}~\mathcal{Q}_2 \neq \emptyset$}
		\While{$\mathcal{Q}_1 \neq \emptyset$} \label{ddqm:prop:1}
			\State $\xii \leftarrow \mathcal{Q}_1\Call{.front}$ \Comment{Extracts the front element.}
			\State $\Ttilde \leftarrow \Call{SolveEikonal}{\xii, \T, \F}$ \label{ddqm:prop:2}
			\If{$\Ttilde < \Ti$}  \label{ddqm:prop:3}
				\State $\Ti \leftarrow \Ttilde$
				\For{$\xjj \in (\mathcal{N}(\xii)\cap\Frozen$)}
					\If{$\Ti < \Tj$} \Comment{Add improvable neighbors to corresponding queue.}
						\State $\Frozen \leftarrow \Frozen \backslash \{\xjj\}$
						\State $\Narrow \leftarrow \Narrow \cup \{\xjj\}$
						\State $c_{\mathrm{total}} \leftarrow c_{\mathrm{total}}+1$
						\If{$\Ti \leq th$}
							\State $\mathcal{Q}_1 \leftarrow \mathcal{Q}_1 \cup \{\xjj\}$
							\State $c_1 \leftarrow c_1+1$
						\Else
							\State $\mathcal{Q}_2 \leftarrow (\mathcal{Q}_2 \cup \{\xjj\} )$  \label{ddqm:prop:4}
						\EndIf
					\EndIf
				\EndFor
			\EndIf
			\State $\Narrow \leftarrow \Narrow \backslash \{\xii\}$ \label{ddqm:prop:5}
			\State $\Frozen \leftarrow \Frozen \cup \{\xii\}$
		\EndWhile
		\State $step \leftarrow \Call{UpdateStep}{step, c_1, c_{\mathrm{total}}}$
		\State $\Call{swap}{\mathcal{Q}_1, \mathcal{Q}_2}$
		\State $c_1 \leftarrow 0$
		\State $c_{\mathrm{total}} \leftarrow 0$
		\State $th \leftarrow th+step$ \label{ddqm:prop:6}
	\EndWhile	
	\State \textbf{return} $\T$
\EndProcedure
\end{algorithmic}
\caption{Double Dynamic Queue Method}
\label{alg:ddqm}
\end{algorithm}

In the original work, three methods were proposed: 1) single-queue (SQ), and therefore simpler algorithm, 2) two-queue static (TQS), where the $step$ is not updated, and 3) two-queue dynamic (which we call DDQM). SQ and TQS slightly improve DDQM in some experiments, but when DDQM improves SQ and TQS (for instance environments with noticeable speed changes) the difference can reach one order of magnitude. Therefore, we decided to include DDQM instead of SQ and TQS since it has shown a more adaptive behaviour. In any case, any of these methods return the same solution as FMM.

In the worst case, the whole grid is contained in both queues and traversed many times during the propagation. However, since queue insertion and deletion are \Octe~ operations, the overall complexity is \On. Note that $\Call{SWAP}$ can be efficiently implemented in \Octe~ as a circular binary index, or updating references (or pointers). There is not need for a real swap operation.

\subsection{Fast Iterative Method}
\label{sec:other:fim}
The Fast Iterative Method (FIM) \cite{Jeong08} is based on the iterative method proposed by \cite{Rouy92} but inspired in FMM. It also resembles to DDQM (concretely to its single queue variant). It iteratively evaluates every point in $\Narrow$ until it converges. Once a node has converged its neighbors are inserted into $\Narrow$ and the process continues. $Narrow$ is implemented as a non-sorted list. The algorithm requires a convergence parameter $\epsilon$: if $\Ti$ is improved less than $\epsilon$, it is considered converged. FIM has been also proposed for triangulated surfaces \cite{Fu11}.

FIM is designed to be efficient for parallel computing, since all the elements in $\Narrow$ can be evaluated simultaneously. However, we are focusing on its sequential implementation in order to have a fair comparison with other methods.

\cref{alg:fim} details FIM. Its initialization is the same as FMM. Then, for each element in $\Narrow$, its value is updated (lines \ref{fim:prop:1}-\ref{fim:prop:2}). If the value difference is less than $\epsilon$, the neighbors are evaluated and added to $\Narrow$ in case their value is improved (lines \ref{fim:prop:3}-\ref{fim:prop:4}). Since $\Narrow$ is a list, the new elements should be inserted just before the point being currently evaluated, $\xii$. Finally, this point is removed from $\Narrow$ and labeled as $\Frozen$ (lines \ref{fim:prop:5} and \ref{fim:prop:6}).

\begin{algorithm}
\begin{algorithmic}[1]
\Procedure{FIM}{$\grid, \T, \F, \Xstart, \epsilon$}
	\Statex Initialization:
	\State $\Frozen \leftarrow \grid, \Narrow \leftarrow \emptyset$
	\State $\Ti \leftarrow \infty~\forall \xii \in \grid$
	\For{$\xii \in \Xstart$}
		\State $\Ti \leftarrow 0$
		\For{$\xjj \in (\mathcal{N}(\xii)\cap\Unknown)$}
			\State $\Frozen \leftarrow \Frozen \backslash \{\xii\}$
			\State $\Narrow \leftarrow \Narrow \cup \{\xii\}$
		\EndFor
	\EndFor
	\Statex
	\Statex Propagation:
	\While{$\Narrow \neq \emptyset$}
		\For{$\xii \in \Narrow$}
		\State $\Ttilde \leftarrow \Ti$ \label{fim:prop:1}
		\State $\Ti \leftarrow \Call{SolveEikonal}{\xii, \T, \F}$ \label{fim:prop:2}
			\If{$\abs{\Ti - \Ttilde} < \epsilon$} \label{fim:prop:3}
				\For{$\xjj \in (\mathcal{N}(\xii)\cap\Frozen)$}
					\State $\Tjtilde \leftarrow \Call{SolveEikonal}{\xjj, \T, \F}$
					\If{$\Tjtilde < \Tj$}
						\State $\Tj \leftarrow \Tjtilde$
						\State $\Frozen \leftarrow \Frozen \backslash \{\xjj\}$
						\State $\Narrow \leftarrow \Narrow \cup \{\xjj\}$ \Comment{Insert in the list just before $\xii$} \label{fim:prop:4}
					\EndIf
				\EndFor
				\State $\Narrow \leftarrow \Narrow \backslash \{\xii\}$ \label{fim:prop:5}
				\State $\Frozen \leftarrow \Frozen \cup \{\xii\}$ \label{fim:prop:6}
			\EndIf
		\EndFor
	\EndWhile	
	\State \textbf{return} $\T$
\EndProcedure
\end{algorithmic}
\caption{Fast Iterative Method}
\label{alg:fim}
\end{algorithm}

A node can be added several times to $\Narrow$ during FIM execution, since every time an upwind (parent) neighbor is updated, the node can improve its value. In the worst case, $\Narrow$ contains the whole grid and the loop would go through all the points several times. Operations on the list are \Octe. Therefore, the overall computational complexity of FIM is \On.

For a small enough $\epsilon$ (depending on the environment), FIM will return the same solution as FMM. However, it can be speed up allowing small errors bounded by $\epsilon$.

\section{Experimental Comparison}
\label{sec:expcomp}
\subsection{Experimental setup}
\label{sex:expcomp:setup}
In order to give an impartial and meaningful comparison, all the algorithms have been implemented from scratch, in C++11 using the Boost.Heap library\footnote{The source code is available at https://github.com/jvgomez/fastmarching}. An automatic benchmarking application has been also created so that the experiments are carried out and evaluated in the most systematic possible way.

This implementation is focused on time performance and compiled using G++ 4.9.2 with optimizations flag -Ofast. However, no special optimizations have been included. All algorithms use the same primitive functions for grid and cell computations. The times reported correspond to an Ubuntu 14.04 64 bits computer running in a Dual Core 3.3 GHz with 4Gb of RAM. However, all experiments were carried out in one core. Only propagation times are taking into account. The computation time used in the initialization has been omitted since it can be done offline, besides it is similar for all algorithms and represents a little percentage of the total computation time.

Although the algorithms are deterministic (their result only depends on the grid, which is not modified), the results shown are the mean of 10 runs for every algorithm, so that the deviation of the results is practically 0.

For UFMM, the default is a maximum range of $\T$ of 2 units and 1000 buckets (the checkerboard experiment required different parameters, see \cref{res:expcom:setup:checker}). The $\epsilon$ parameter for FIM is set to 0 (actually $10^{-47}$ to provide robust 64bit double comparison).

Although error analysis is not in the scope of this paper, these can be compared among the existing papers since they are implementation-independent. UFMM errors are reported in those experiments with non-constant velocity. Usually, $L_1$ and $L_{\infty}$ norms of the error are reported. Most of the works compute norm $L_1$ as:

\begin{equation}
|\T|_1 = \sum\limits_{\xii\in\grid} |\Ti|
\label{eq:l1_Bad}
\end{equation}

\noindent where $\grid$~is treated as a regular vector. However, following \cite{Chacon14}, we treat the numerical solutions as elements of $L_p$ spaces (generalization of the $p$-norm to vector spaces), where $L_1$ norm is defined as an integral over the function. The result is a norm closely related to its physical meaning and independent of the cell size. $L_1$ is numerically integrated over the domain and therefore computed as (assuming cubic cells and grids):

\begin{equation}
|\T|_1 = \sum\limits_{\xii\in\grid} |\Ti h^N| = h^N\sum\limits_{\xii\in\grid} |\Ti|
\label{eq:l1}
\end{equation}

Four different experiments have been carried out, which represent characteristic cases. By combining these problems it is possible to get close to any situation. These experiments were chosen attending also to the most common situations tested in the literature.

\subsubsection{Empty map}
This experiment is designed to show the performance of the methods in the most basic situation, were most of the algorithms perform best. An empty map with constant velocity represents the simplest possible case for the Fast Methods. In fact, analytical methods could be implemented by computing the euclidean distance from every point to the initial point. However, it is interesting because it shows the performance of the algorithms on open spaces which, in a real application, can be part of large environments.

The same environment is divided into a different number of cells to study how the algorithms behave as the number of cell increases. Composed by an empty 2D, 3D and 4D hyper-cubical environment of size $[0,1]^N$, with $N=2,3,4$. Constant velocity $\Fi = 1~\text{on}~\Omega$. The wavefront starts at the center of the grid.

The number of cells  was chosen so that an experiment has the same (or as close as possible) number of cells in all dimensions. For instance, a 50x50 2D grid has 2500 cells. Therefore, the equivalent 3D grid is 14x14x14 (2744) and in 4D is 7x7x7x7 (2401). This way, it is possible to also analyze the performance of the algorithms for a different number of dimensions. Thus, we have chosen the following number of cells for each dimension for 2D grid:

\begin{equation*}
2D: \{50, 100, 200, 400, 800, 1000, 1500, 2000, 2500, 3000, 4000\}
\end{equation*}

Consequently, the 3D and 4D cells are:

\begin{equation*}
\begin{split}
3D:& \{14, 22, 34, 54, 86, 100, 131, 159, 184, 208, 252\}\\
4D:& \{7, 10, 14, 20, 28, 32, 39, 45, 450, 55, 63\}
\end{split}
\end{equation*}

\subsubsection{Alternating barriers}
In this case, we want to analyze how the algorithms behave with obstacles ($\Fi = 0$) in a constant velocity environment ($\Fi = 1$). The obstacles cause the characteristics to change.

\begin{figure}[ht]
	\centering
	\subfloat[1 barrier.]{\includegraphics[width=0.32\columnwidth]{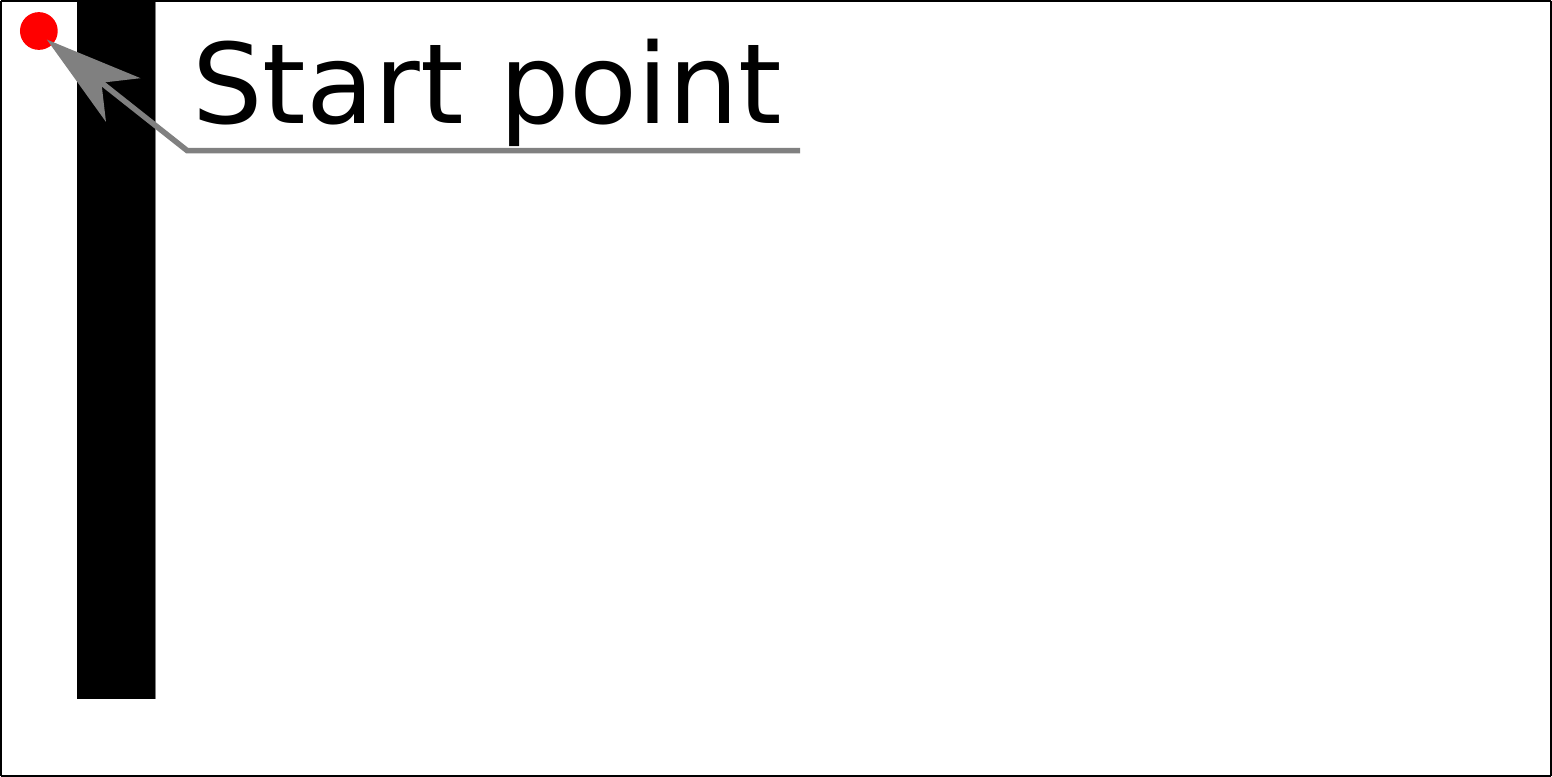}}
	\subfloat[5 barriers.]{\includegraphics[width=0.32\columnwidth]{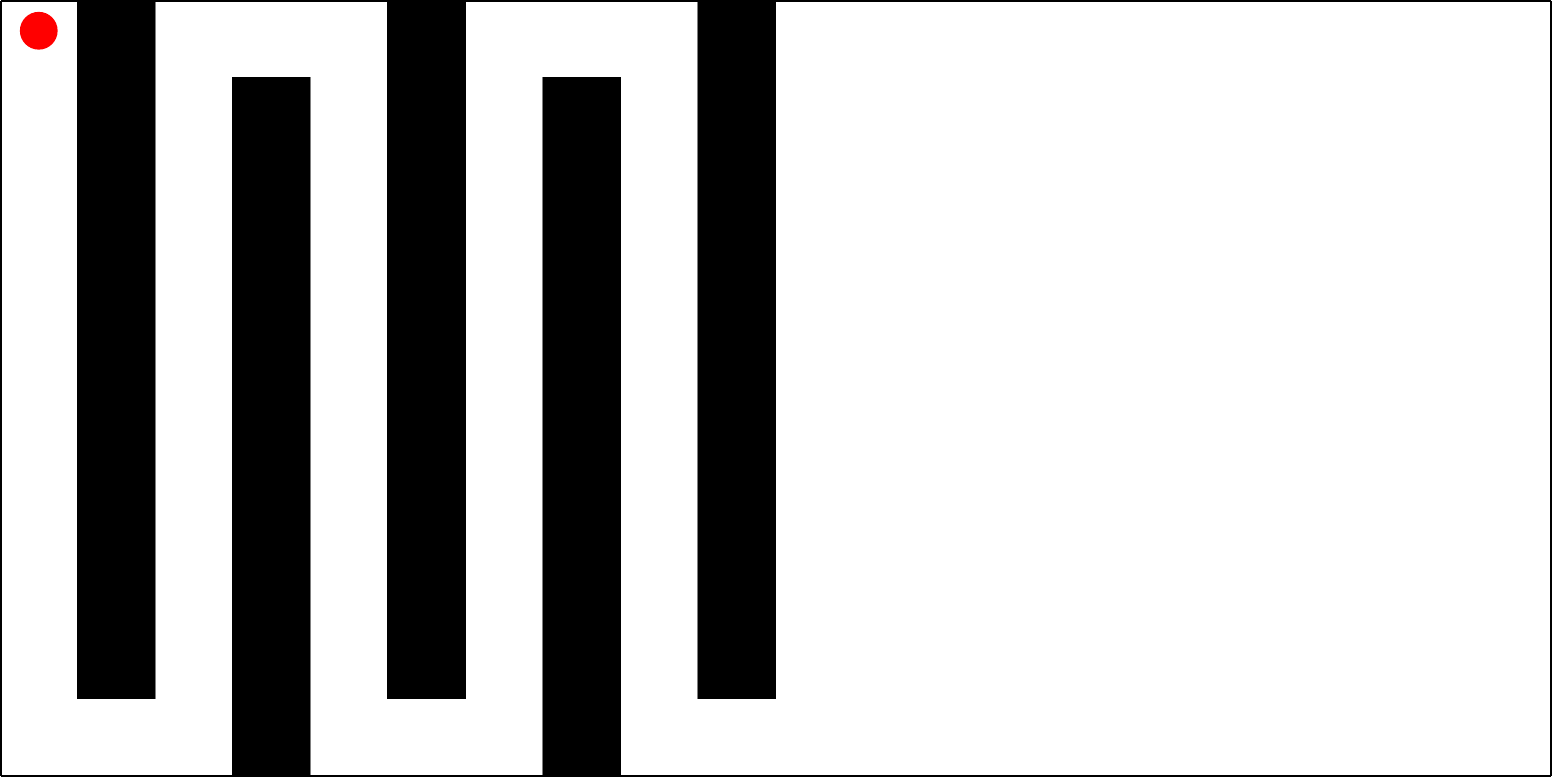}}
	\subfloat[9 barriers.]{\includegraphics[width=0.32\columnwidth]{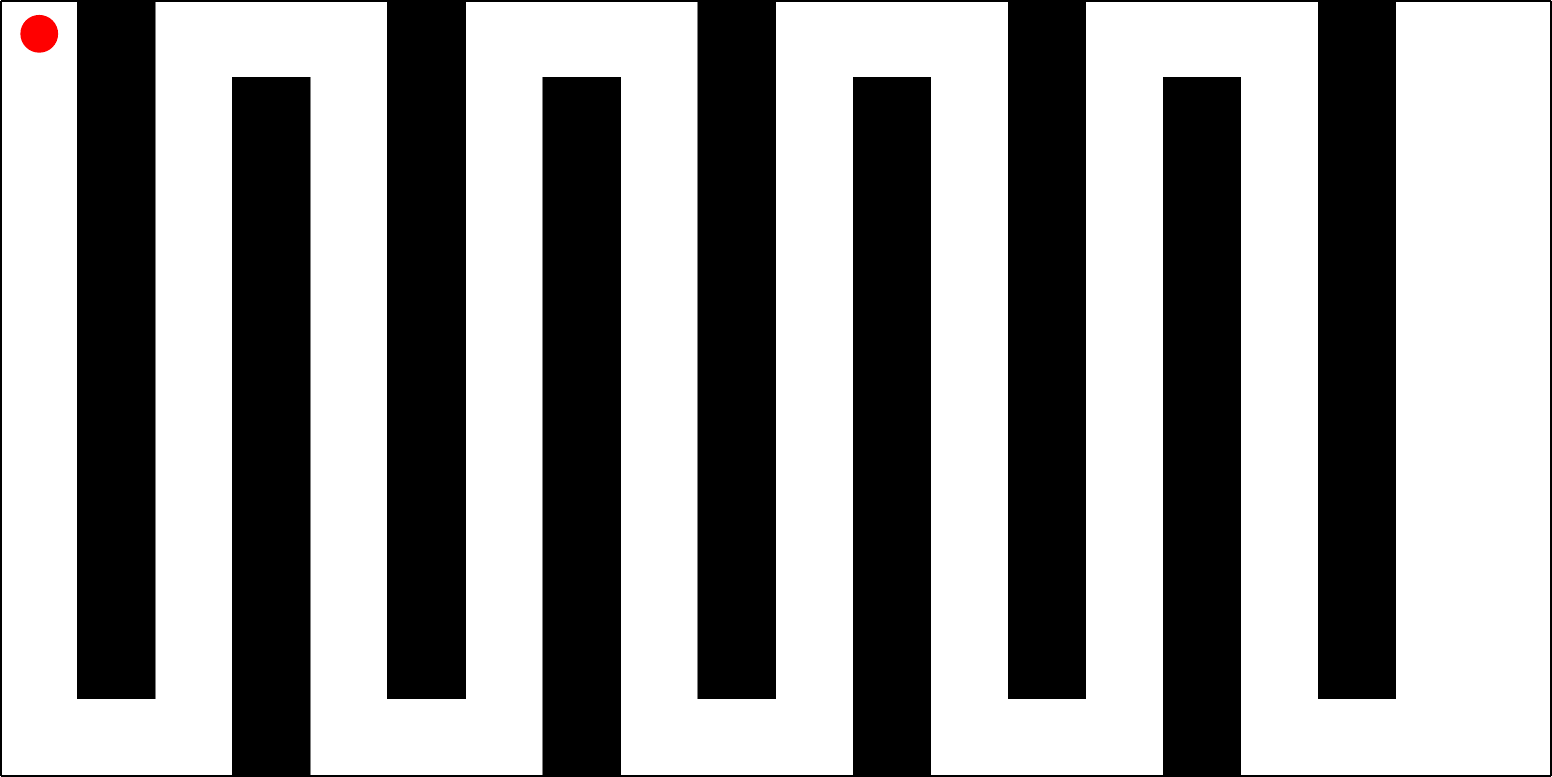}}
	\caption{2D alternating barriers environments.}
\label{fig:barriers}
\end{figure}

The experiment contains a 2D environment of constant size $[0,1]x[0,2]$ discretized in a 1000x2000 grid. A variable number of alternating barriers are equally distributed along the longest dimension. The number of barriers goes from 0 to 9. Examples are shown in \cref{fig:barriers}.  Analogously, in 3D a $[0,1]x[0,1]x[0,2]$ environment represented by a 100x100x200 is chosen, with equally-distributed alternating barriers (from 0 to 9) along the $z$ axis. The wavefront starts in all cases close to a corner of the map.

\subsubsection{Random velocities}
In order to test the performance of the algorithms with random velocities, or noisy images as in the case of medical computer vision. This experiment creates a 2D, 3D and 4D environment of size $[0,1]^N$ with $N=2,3,4$ discretized in a 2000x2000 grid in 2D, 159x159x159 in 3D and 45x45x45x45 in 4D. These discretizations are chosen so that it is possible to compare directly with the empty map problem for the corresponding grid sizes. The wavefront starts in the center of the grid.

\begin{figure}[ht]
	\centering
	\subfloat[Max. velocity = 30]{\includegraphics[width=0.32\columnwidth]{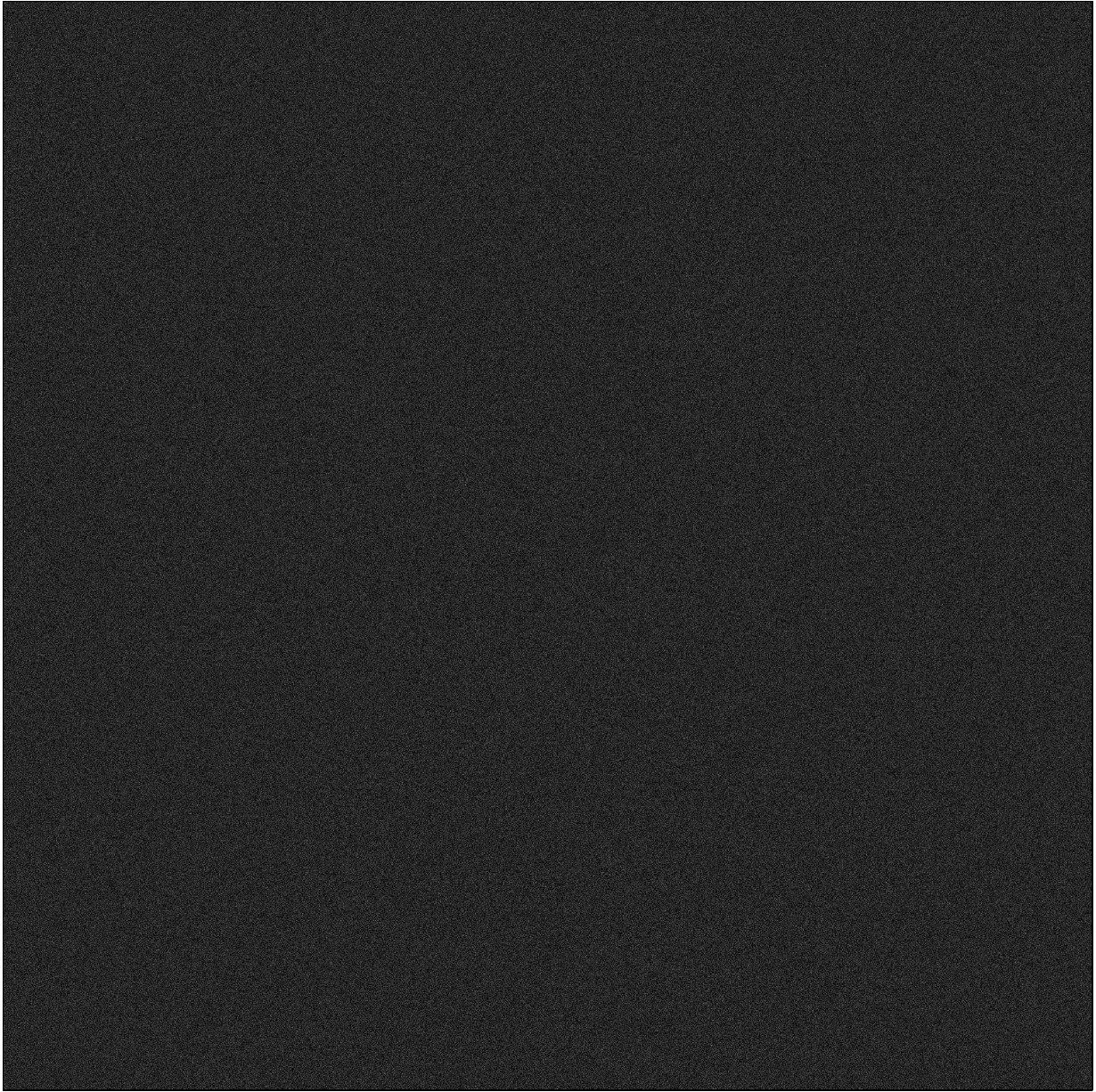}}
	\subfloat[Max. velocity = 60.]{\includegraphics[width=0.32\columnwidth]{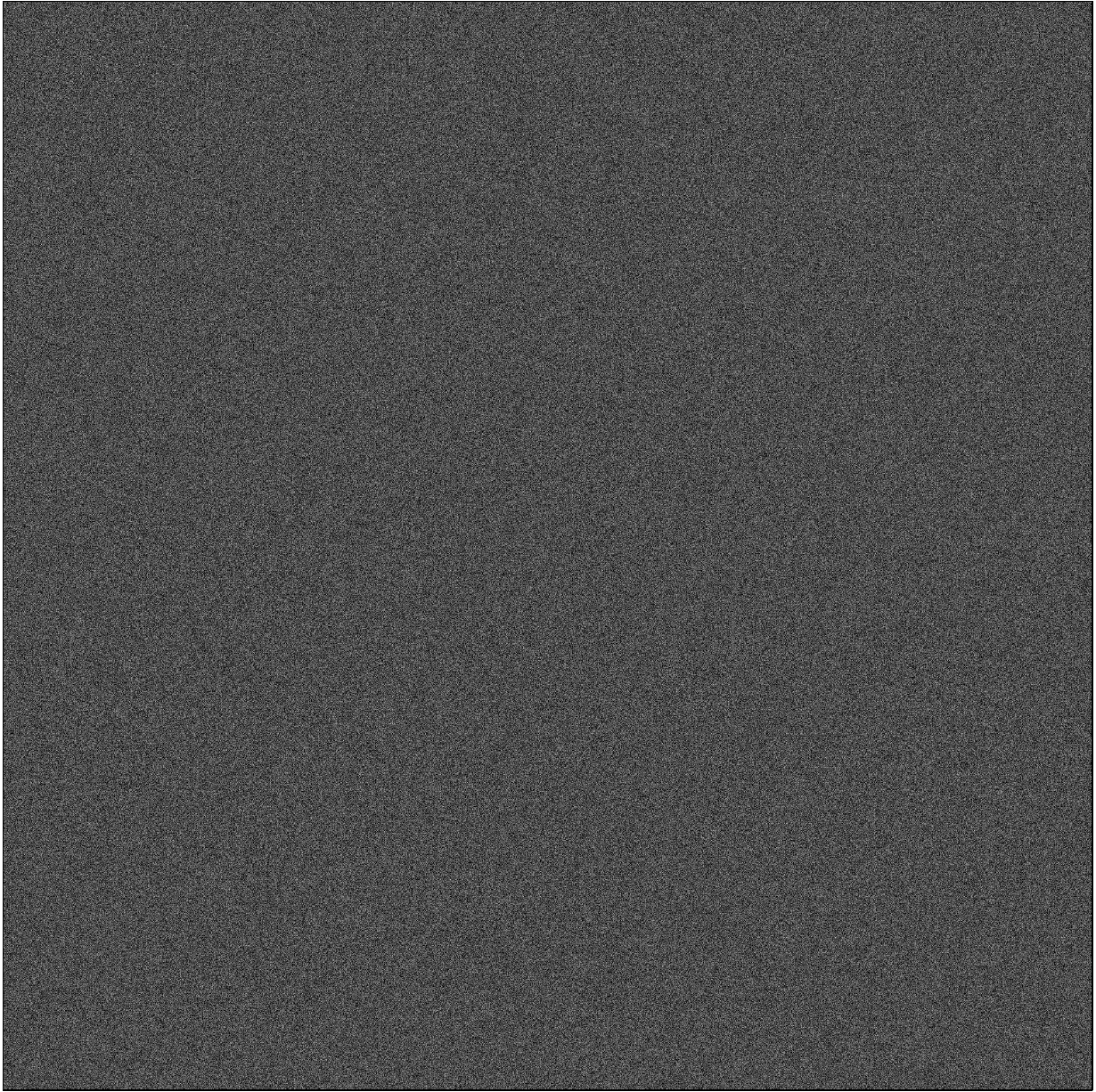}}
	\subfloat[Max. velocity = 100.]{\includegraphics[width=0.32\columnwidth]{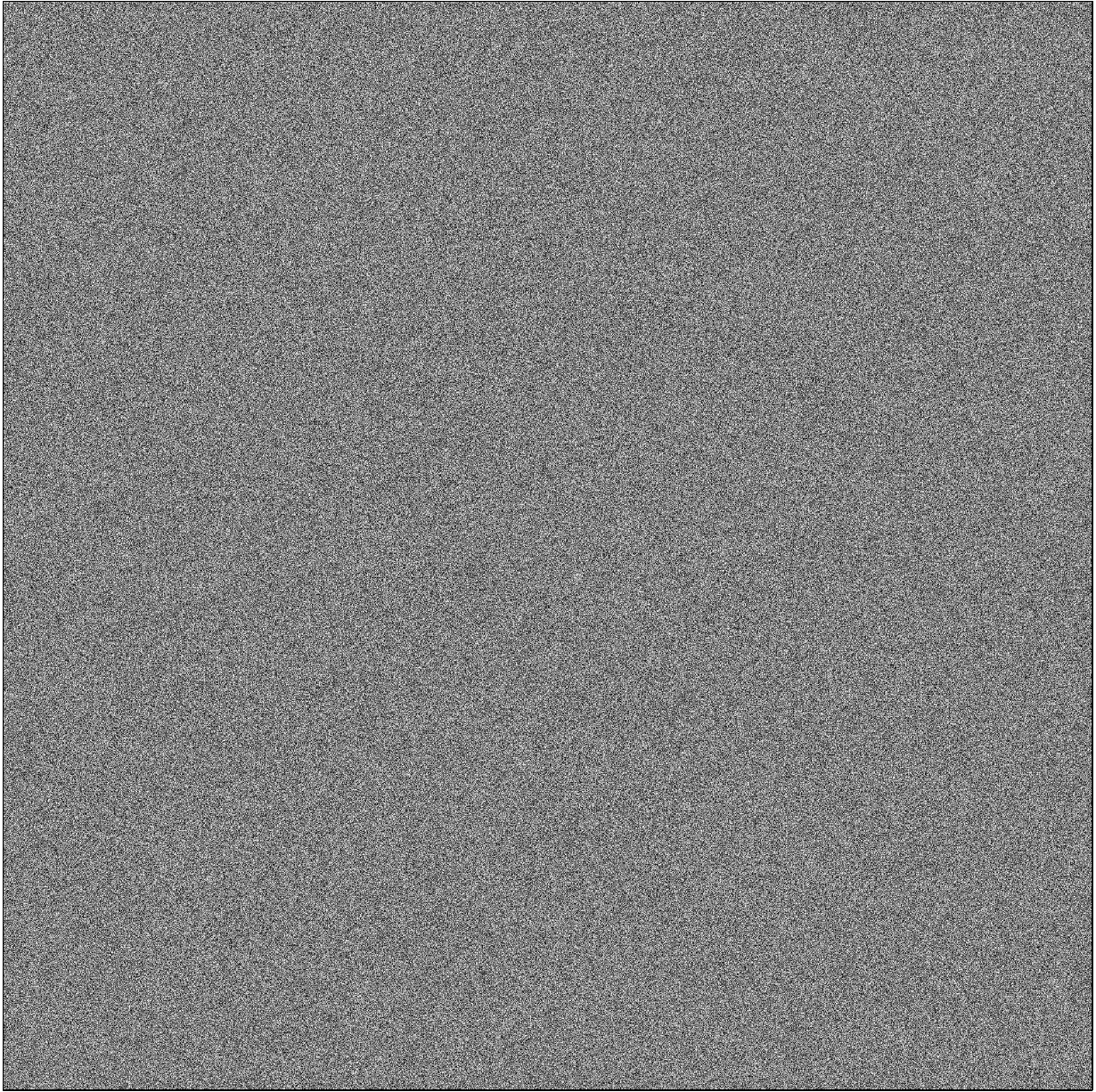}}
	\caption{2D random velocities environments. Lighter color means faster wave propagation.}
\label{fig:random}
\end{figure}

Additionally, the maximum velocity is increased from 10 to 100 (in steps of 10 units) to analyze how the algorithms behave with increasing velocity changes. 2D examples are shown in \cref{fig:random}.

\subsubsection{Checkerboard}
\label{res:expcom:setup:checker}
The random velocities experiment already tested changes in velocities. However, those are high-frequency changes because it is unlikely to have 2 adjacent cells with the same velocities. In this experiment low-frequency changes are studied. The same environment and discretization as in random velocities is now divided like a checkerboard, alternating minimum and maximum velocities. Analogously, the maximum velocity is increased from 10 to 100, while the minimum velocity is always 1. There are 10 checkerboard divisions on each dimension. The wavefront starts in the center of the grid. 2D examples are shown in \cref{fig:random}.

\begin{figure}[ht]
	\centering
	\subfloat[Max. velocity = 30]{\includegraphics[width=0.32\columnwidth]{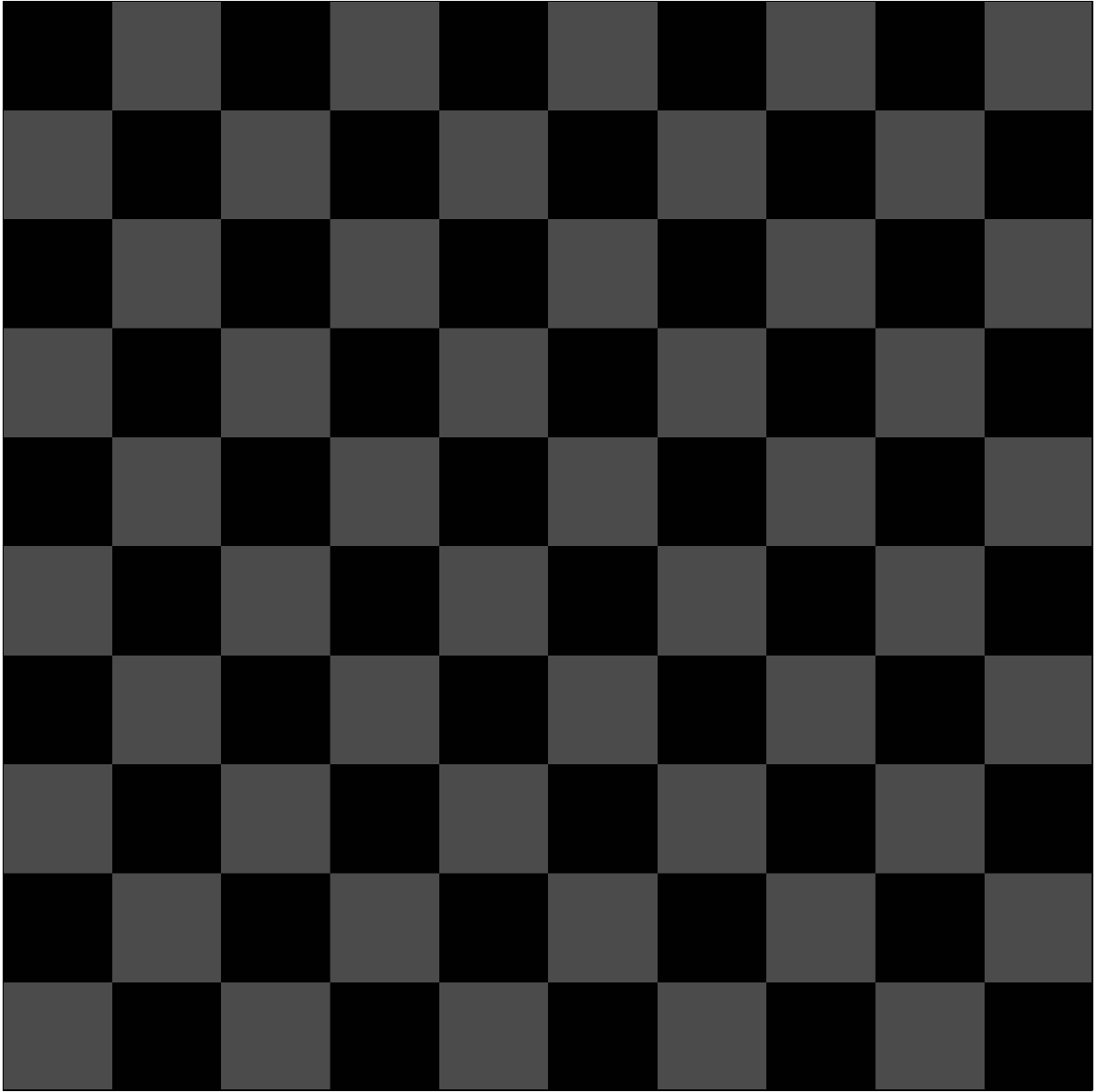}}
	\subfloat[Max. velocity = 60.]{\includegraphics[width=0.32\columnwidth]{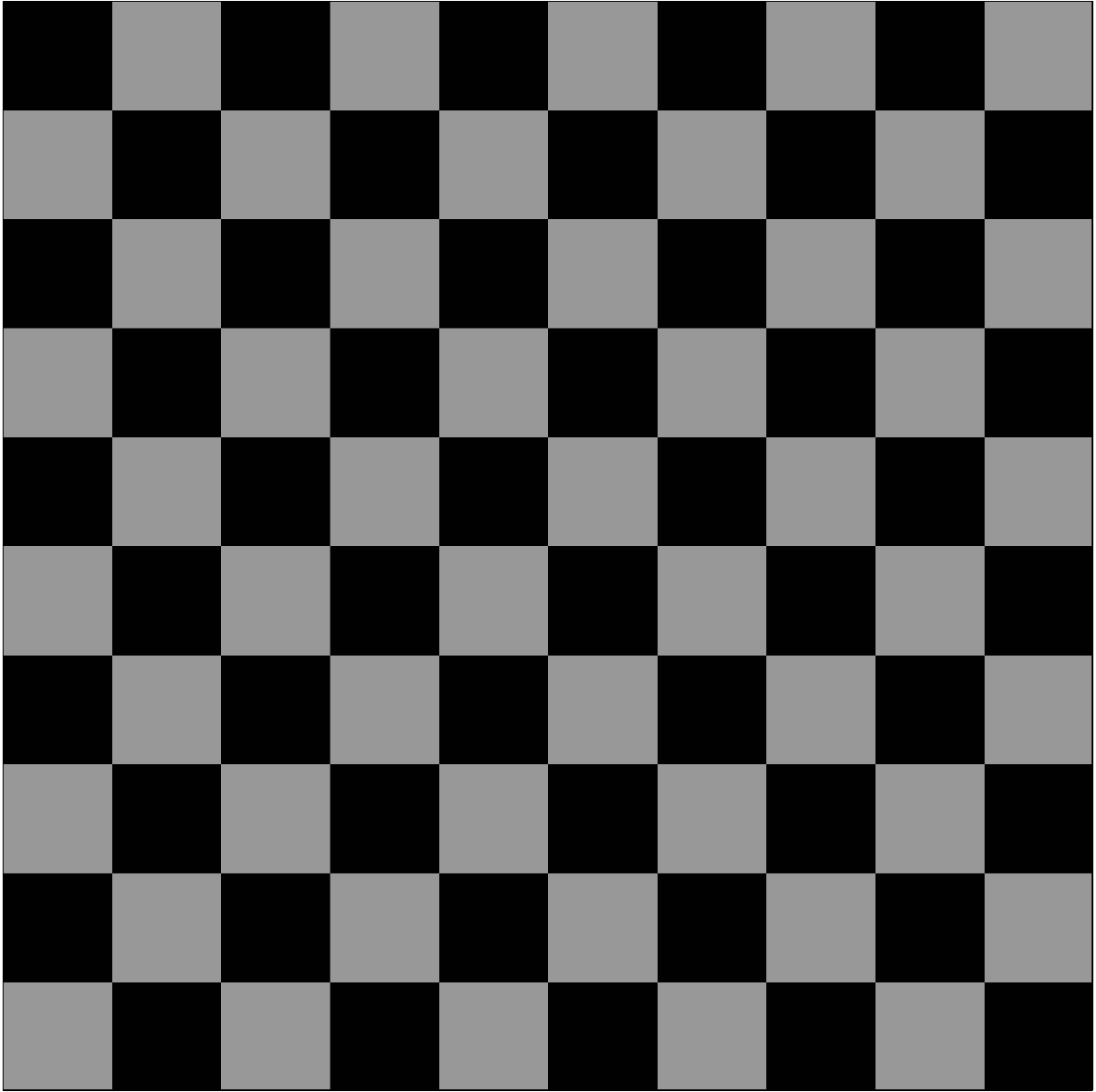}}
	\subfloat[Max. velocity = 100.]{\includegraphics[width=0.32\columnwidth]{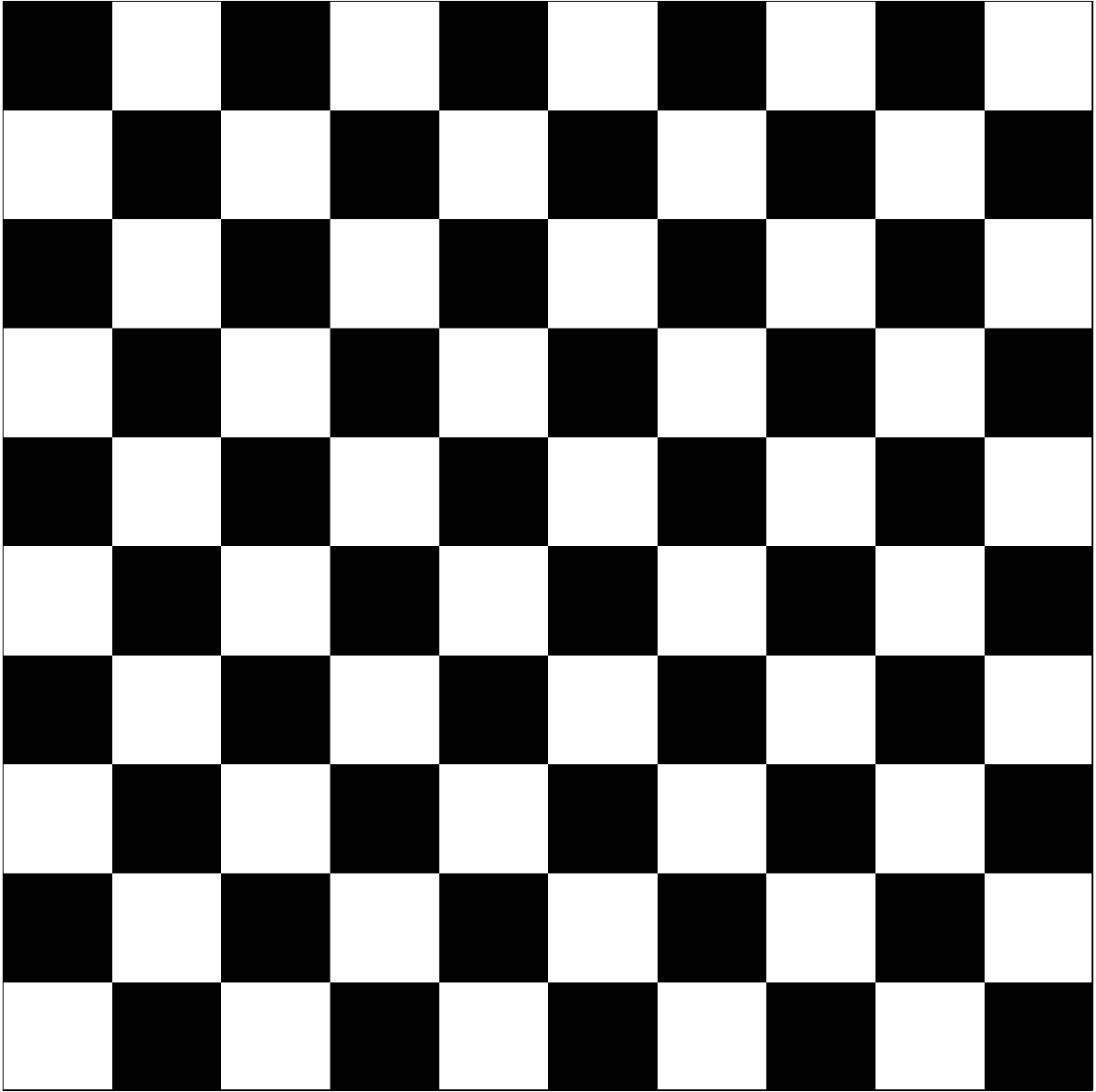}}
	\caption{2D checkerboard environments. Lighter color means faster wave propagation.}
\label{fig:checker}
\end{figure}

In this case, UFMM in 3D and 4D performed very poorly with the default parameters. Our previous experiments show that the best parameters for UFMM are approximately 1000 buckets with a maximum range of 0.01 in 3D. In the 4D case, 20000 buckets and maximum range of 0.025 are used.

\subsection{Results}
\label{sec:expcomp:results}
\subsubsection{Empty map}
\begin{figure}[ht]
	\centering
	\subfloat[50x50]{\includegraphics[width=0.3\columnwidth]{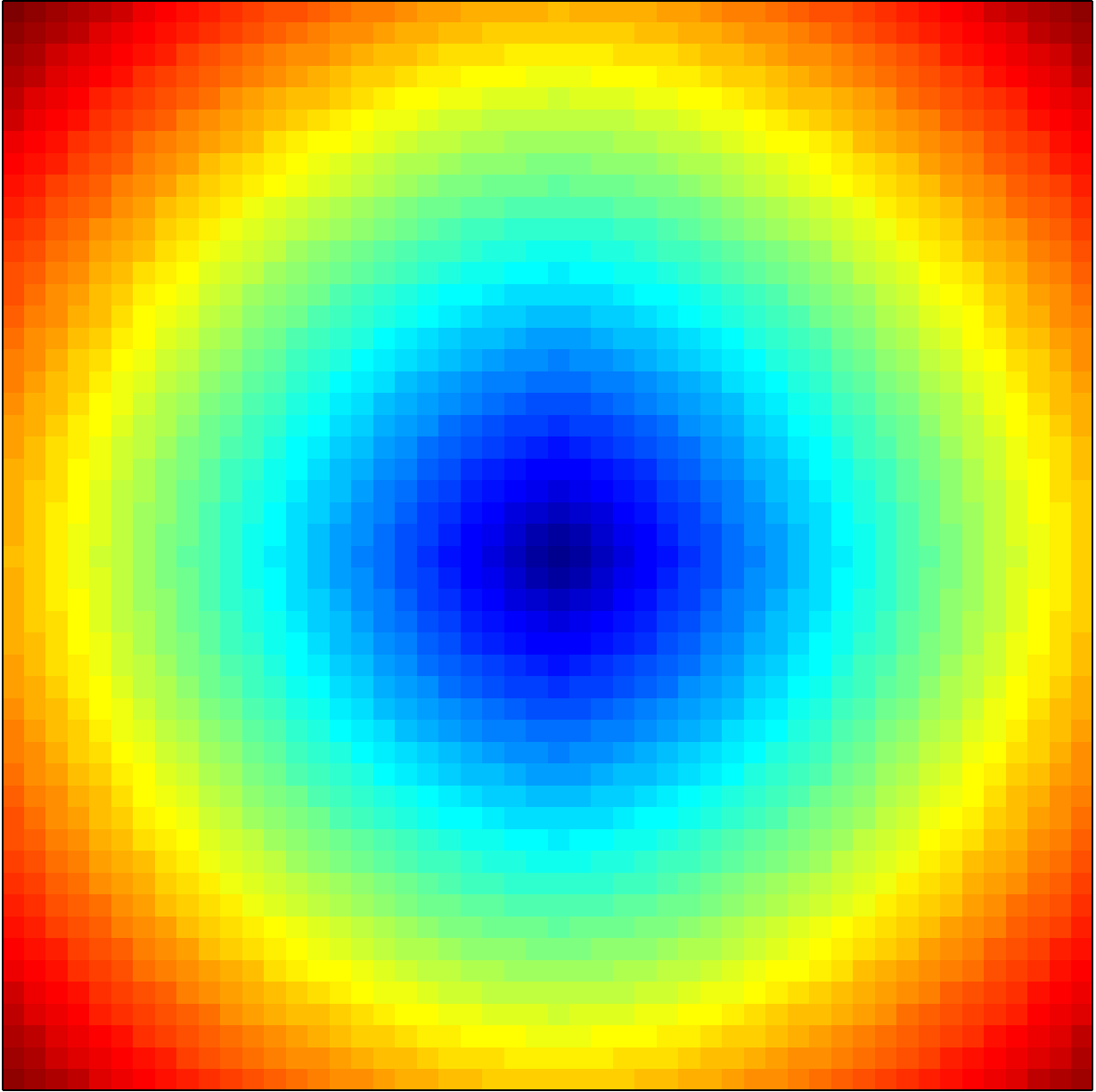}}
	\subfloat[800x800]{\includegraphics[width=0.3\columnwidth]{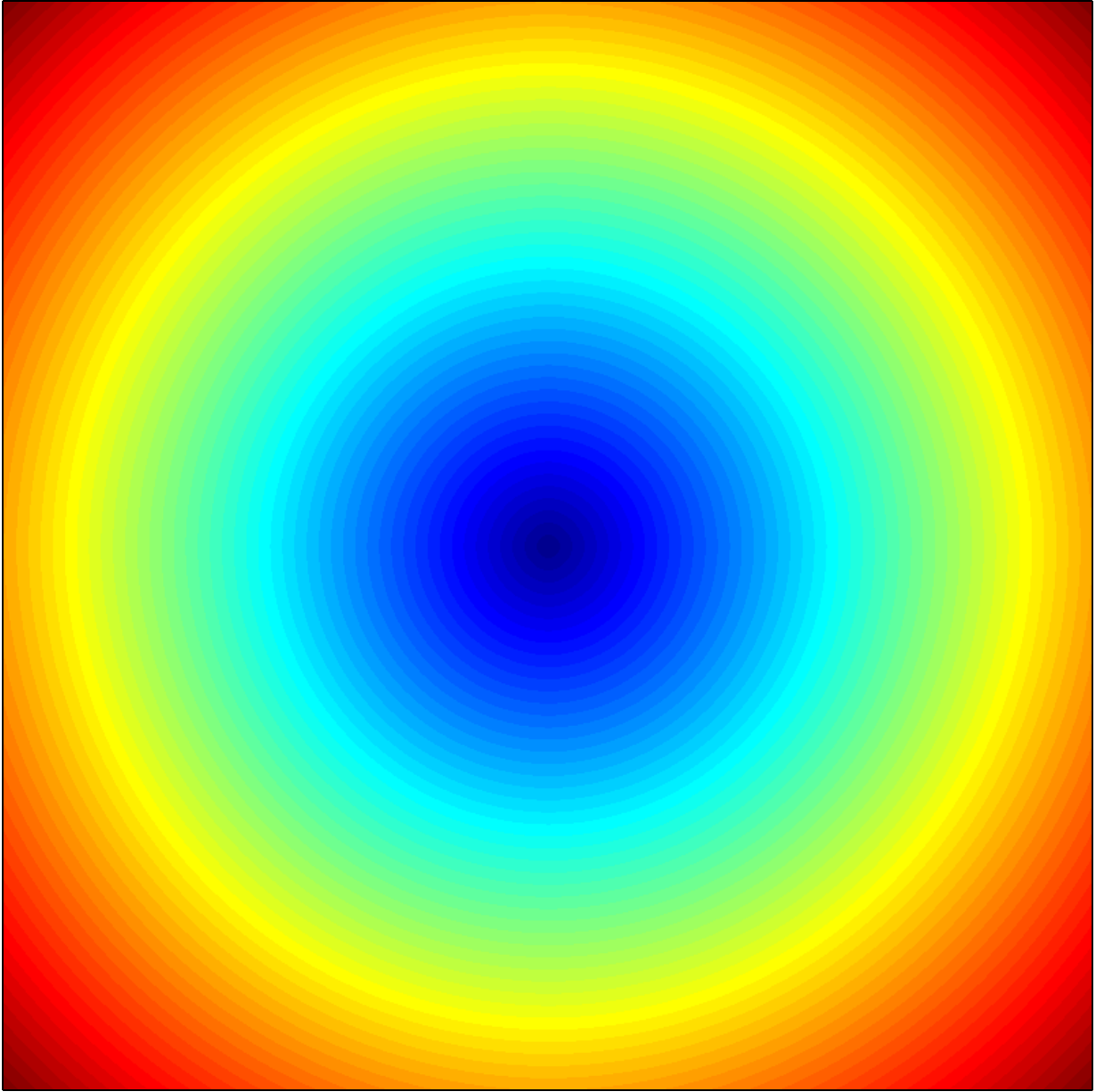}}
	\subfloat[4000x4000]{\includegraphics[width=0.3\columnwidth]{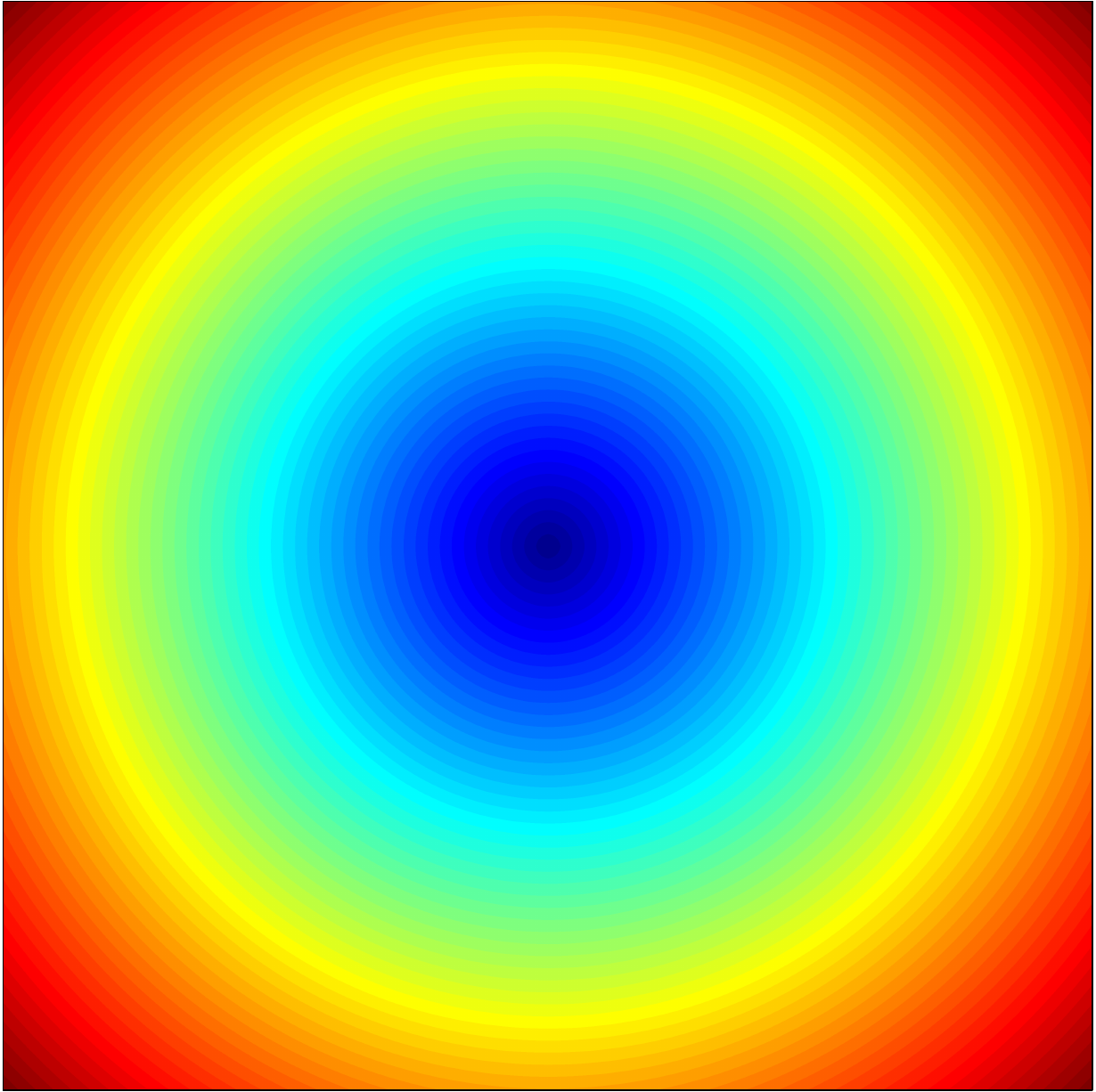}}
	\includegraphics[width=0.0425\columnwidth]{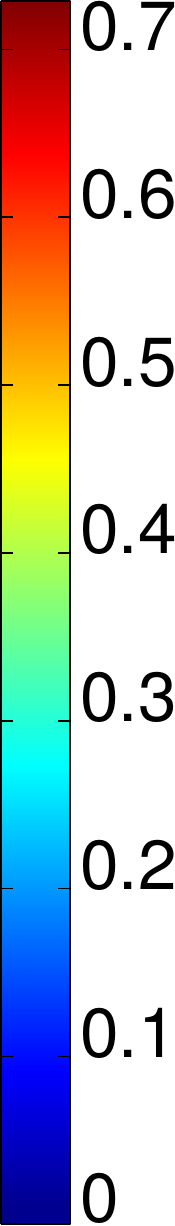}
	\caption{Example of the resulting times-of-arrival maps applying FMM to the empty environment in 2D.}
\label{fig:res:empty}
\end{figure}

An example of the times-of-arrival field computed by FMM is shown in \cref{fig:res:empty}. Note that all algorithms provide the same exact solution in this case. The higher the resolution the better the accuracy.

\begin{figure}[ht]
	\centering
	\subfloat[Computation times.]{\includegraphics[width=0.49\columnwidth]{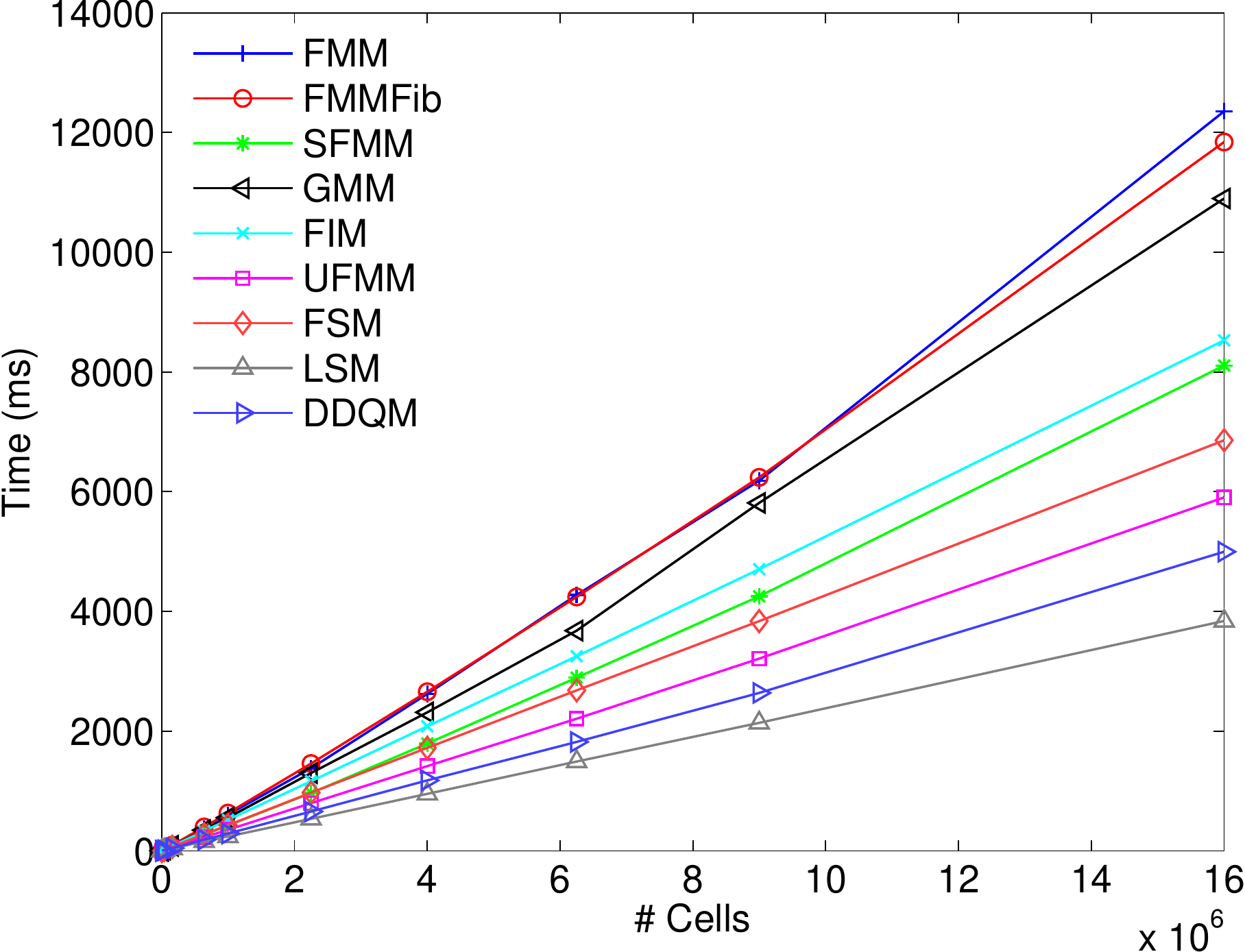}}
	\subfloat[Time ratios against FMM.]{\includegraphics[width=0.47\columnwidth]{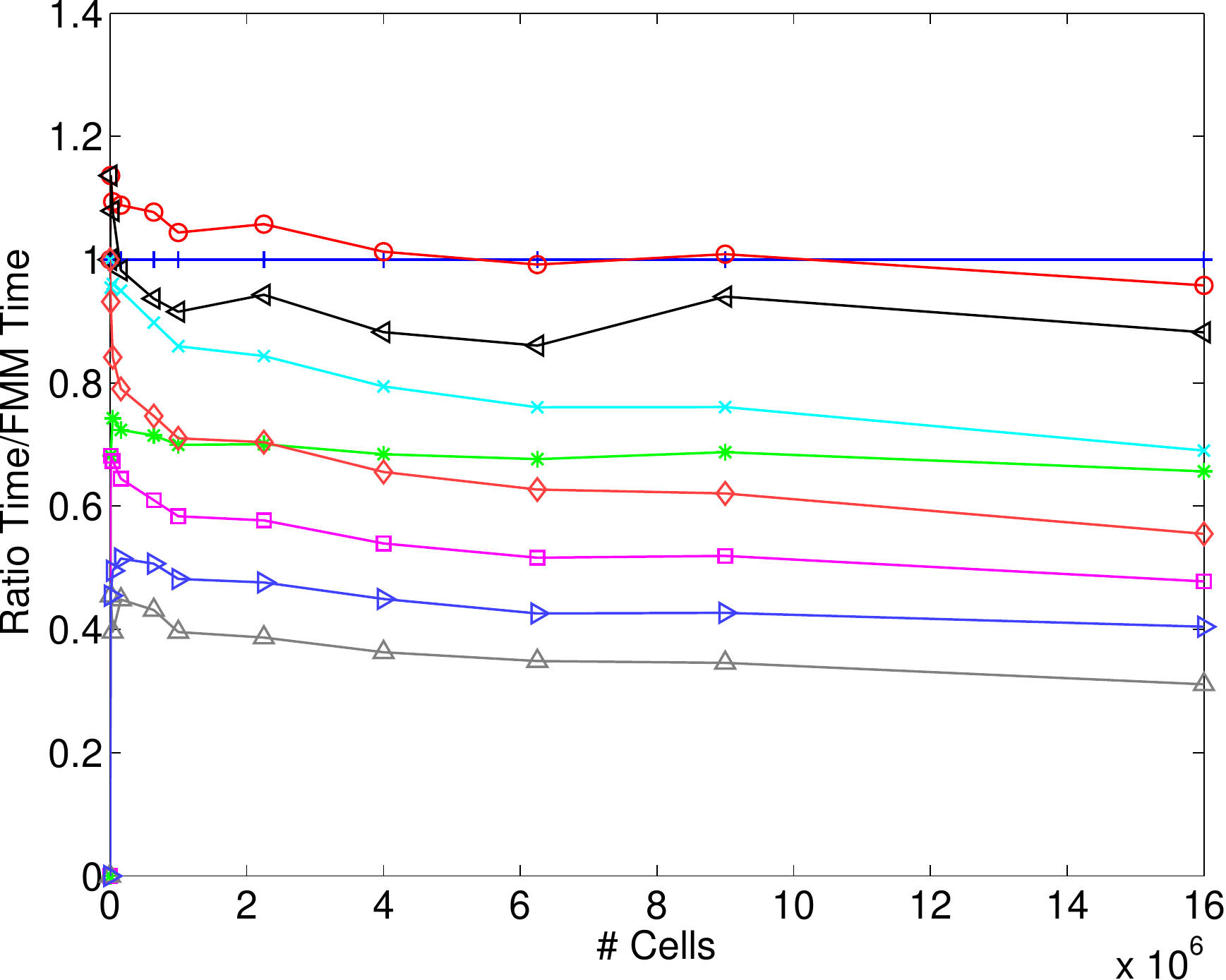}}
	\caption{Computation times and ratios for the empty map environment in 2D.}
\label{fig:res:empty2d}
\end{figure}

\begin{figure}[ht]
	\centering
	\subfloat[Computation times.]{\includegraphics[width=0.49\columnwidth]{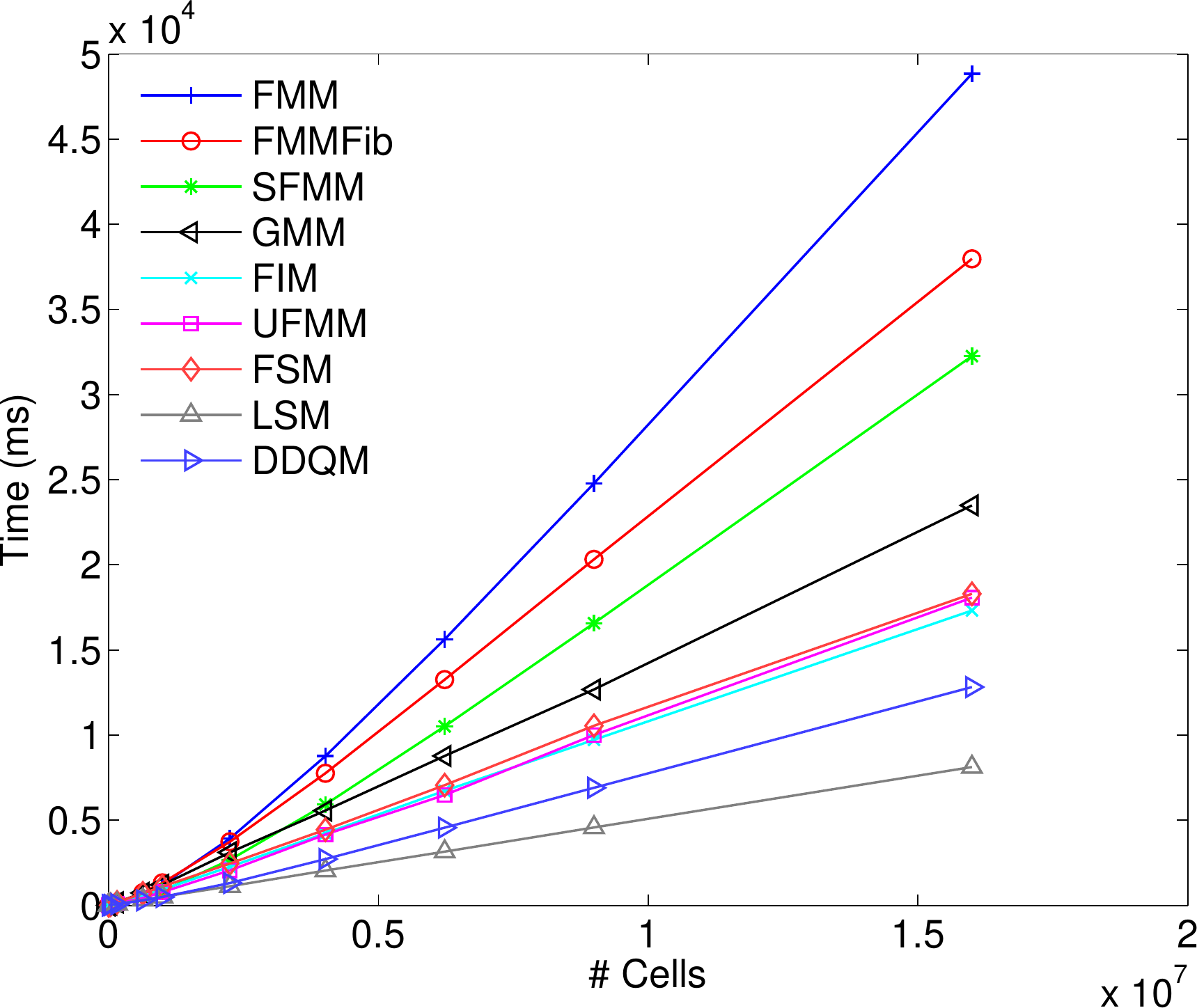}}
	\subfloat[Time ratios against FMM.]{\includegraphics[width=0.49\columnwidth]{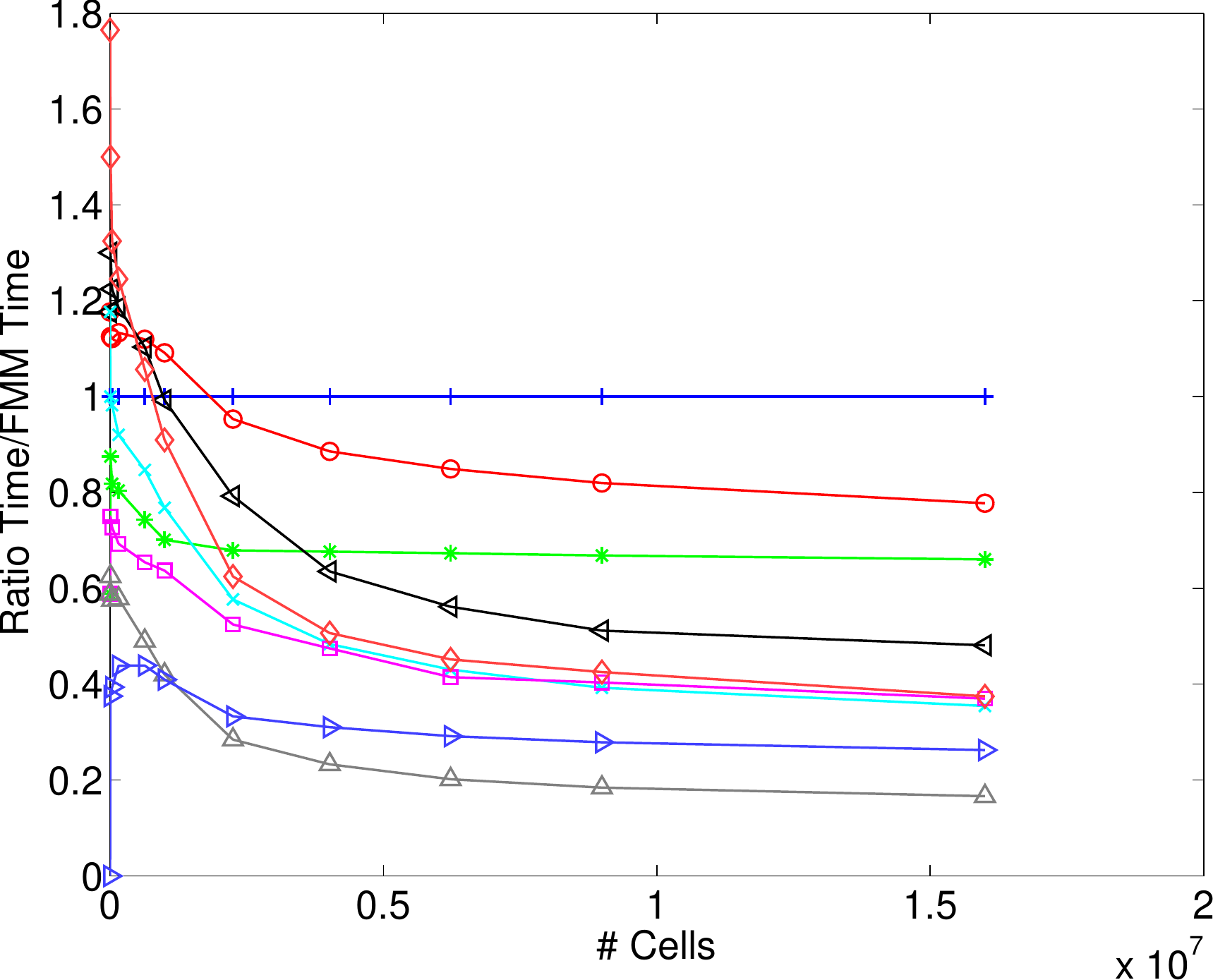}}
	\caption{Computation times and ratios for the empty map experiment in 3D.}
\label{fig:res:empty3d}
\end{figure}

\begin{figure}[ht]
	\centering
	\subfloat[Computation times.]{\includegraphics[width=0.48\columnwidth]{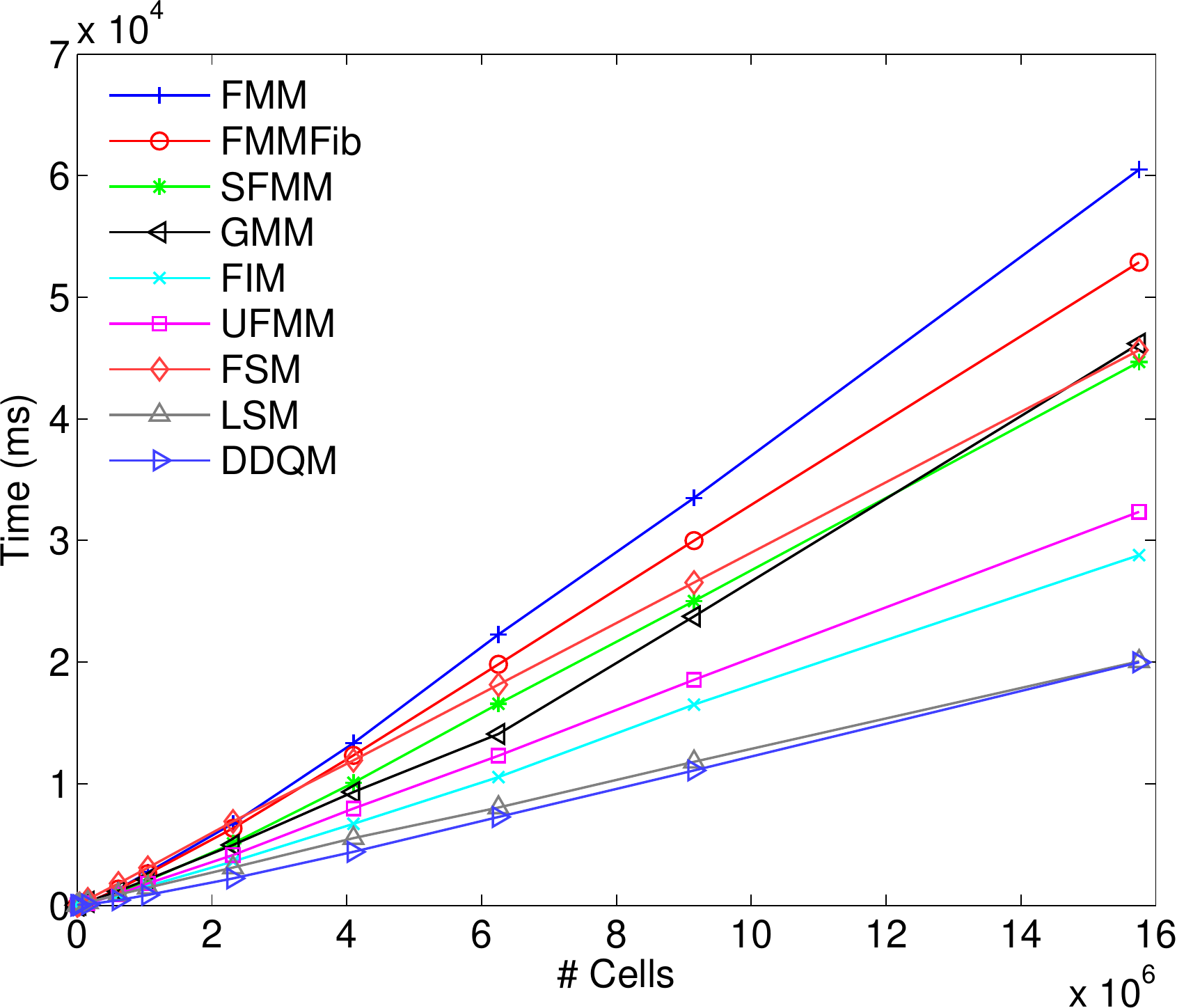}}
	\subfloat[Time ratios against FMM.]{\includegraphics[width=0.49\columnwidth]{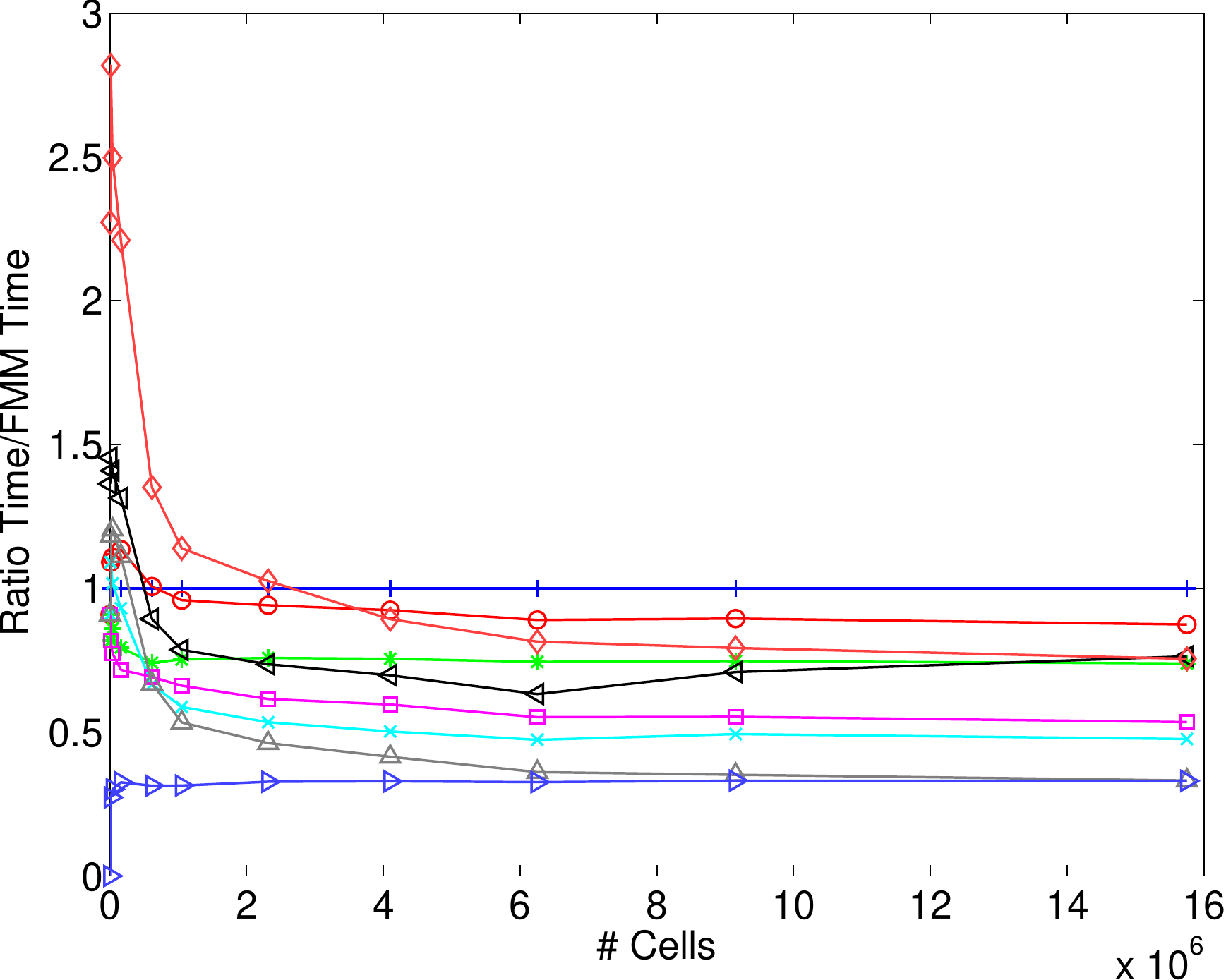}}
	\caption{Computation times and ratios for the empty map experiment in 4D.}
\label{fig:res:empty4d}
\end{figure}

The results for the empty map experiment are shown in \cref{fig:res:empty2d} for 2D, \cref{fig:res:empty3d} for 3D, and \cref{fig:res:empty4d} for 4D. In all cases 2 plots are included: raw computation times for each algorithm, and time ratios computed as:

\begin{equation}
ratio = \frac{\text{Alg.~Time}}{\text{FMM~Time}}
\end{equation}

In all cases, both LSM and DDQM are the fastest algorithms. DDQM tends to perform better in smaller grids, specially en 4D. This was an expected result since this is the case in which they perform less sweeps (queue recomputations in DDQM). Besides, FSM is slower than UFMM in all cases and than FIM in 3D and 4D. As velocities are constant, UFMM provides the same solution as other methods,

GMM slightly improves FMM and FMMFib in 2D. However in higher dimensions this improvement becomes larger, but it always performs worse than UFMM, FIM, DDQM and LSM. In a previous comparison between GMM and FMM \cite{Jeong08}, GMM was about 50\% faster in all cases than FMM. In this results GMM is at most 40\% better. We attribute this difference to implementation, as the heaps for FMM and FMMFib are highly optimized. FIM is always faster as it only needs one iteration trough the narrow band, while GMM always performs two.

SFMM results are of special interest since it is a minor modification of FMM but, however, highly outperforms FMMin all cases, and even FIM in 2D  and small 3D grids. In 4D FIM becomes faster. Besides UFMM outperforms FSM.

As expected, FMMFib is worse than FMM for almost all sizes. However, when dimensions increase FMMFib quickly outperforms FMM as the number of elements in the narrow band increases exponentially with the number of dimensions, therefore the better amortized times of Fibonacci heap become useful.

\subsubsection{Alternating barriers}
\begin{figure}[ht]
	\centering
	\subfloat[1 barrier.]{\includegraphics[width=0.3\columnwidth]{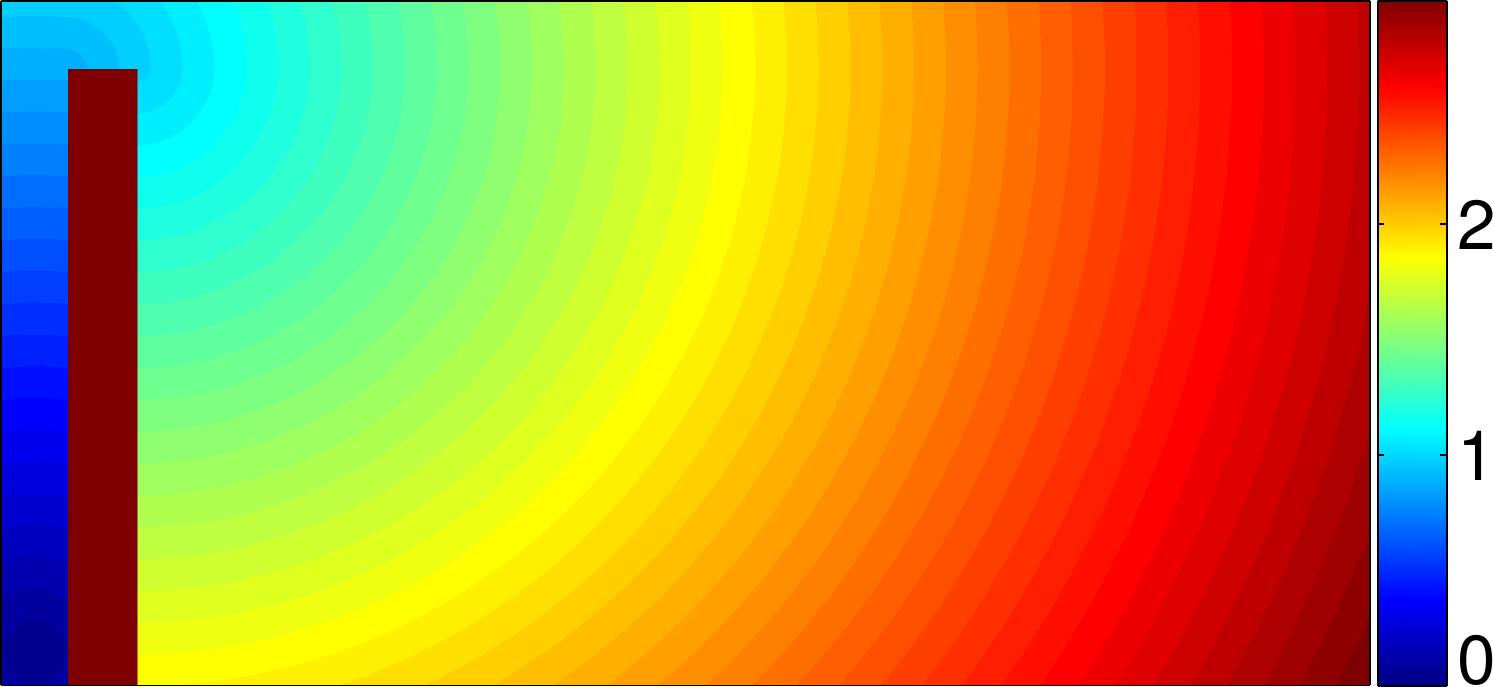}\hspace{0.2cm}}
	\subfloat[5 barriers.]{\includegraphics[width=0.3\columnwidth]{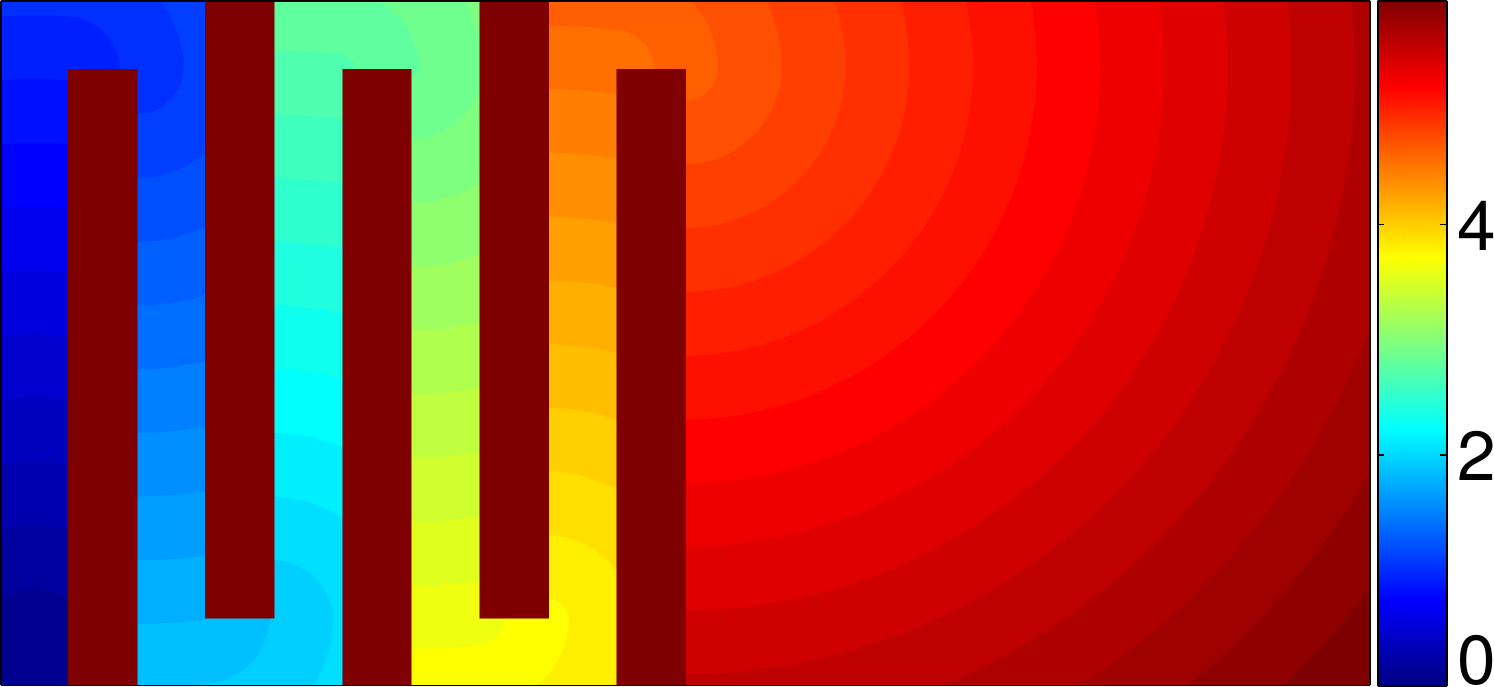}\hspace{0.2cm}}
	\subfloat[9 barriers.]{\includegraphics[width=0.3\columnwidth]{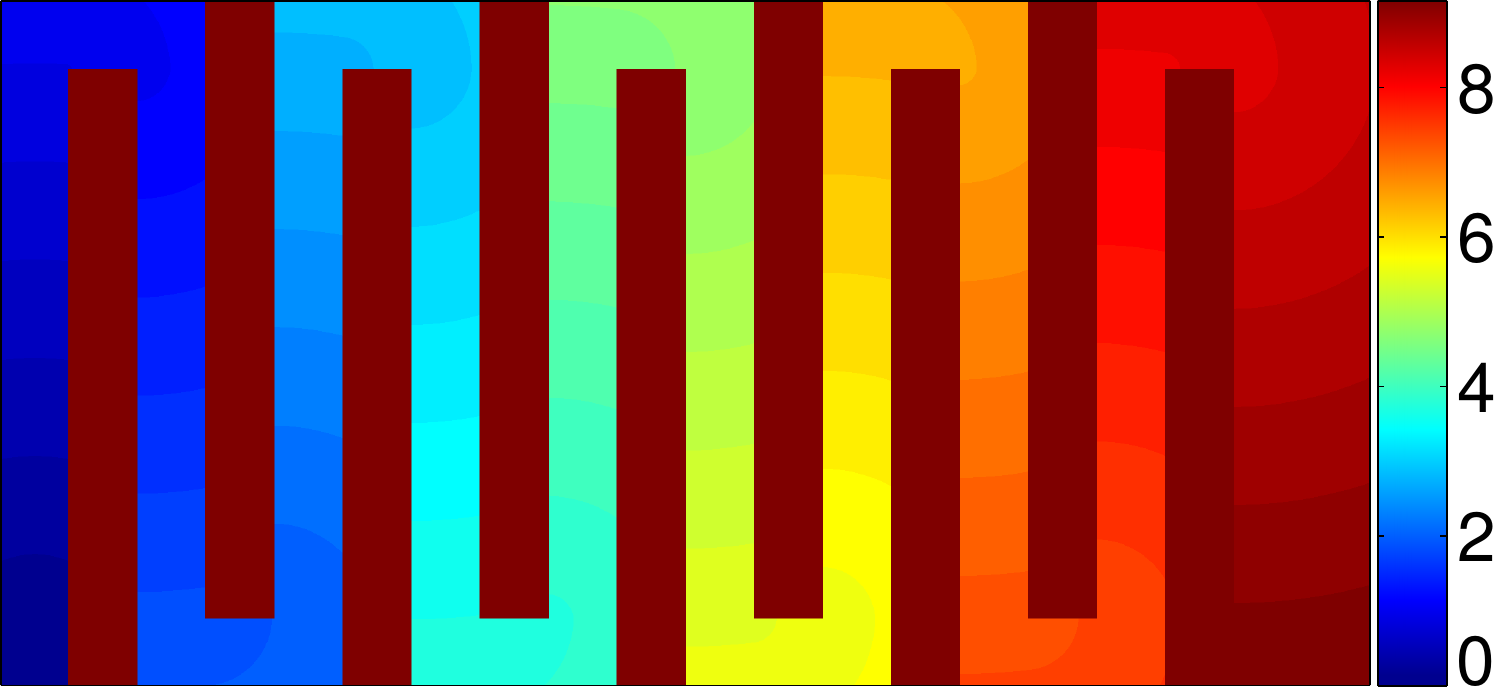}}
	\caption{Example of the resulting times-of-arrival maps applying FMM to some of the alternating barriers environment in 2D.}
\label{fig:res:barriers}
\end{figure}

\begin{figure}[ht]
	\centering
	\subfloat[Computation times.]{\includegraphics[width=0.49\columnwidth]{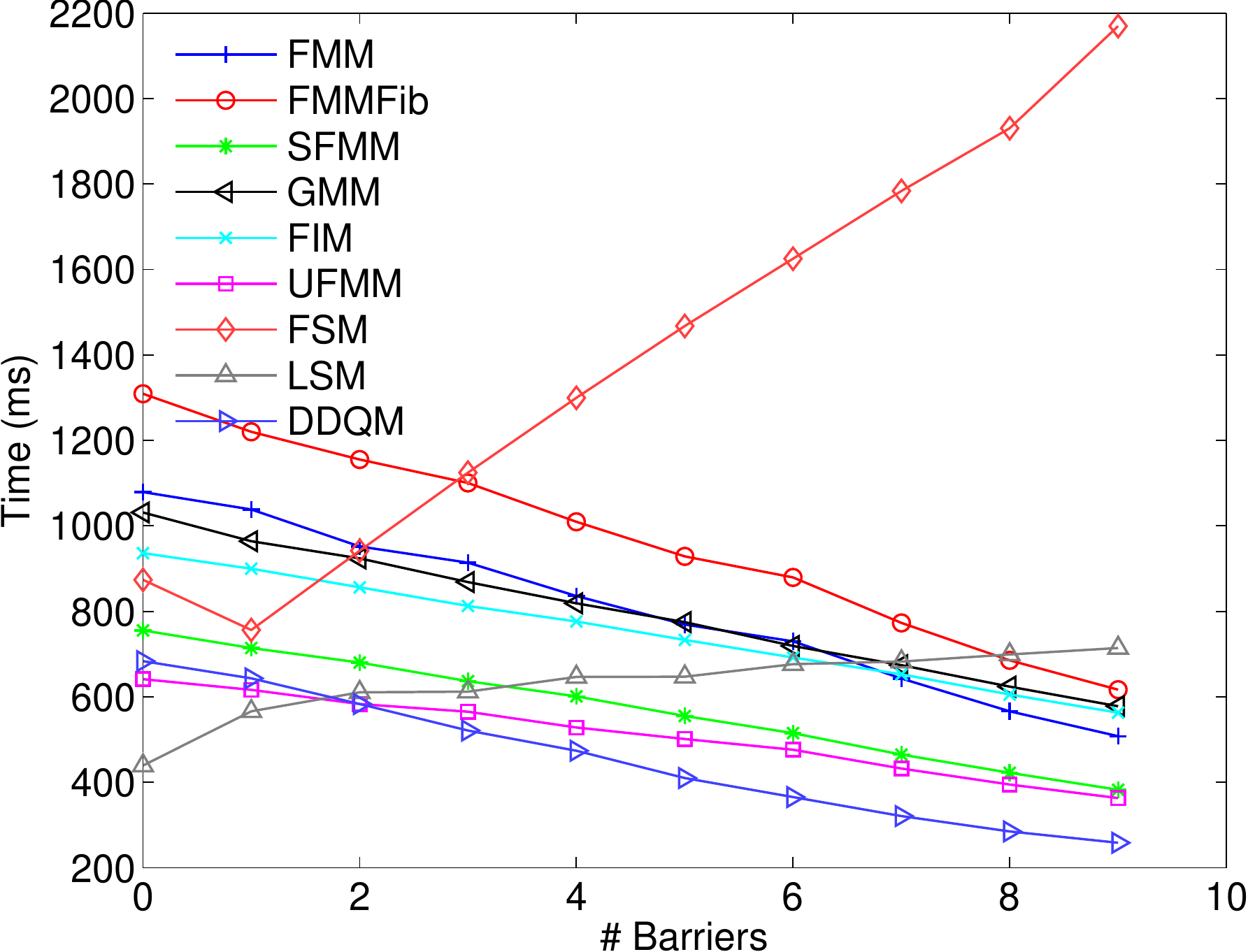}}
	\subfloat[Time ratios against FMM.]{\includegraphics[width=0.48\columnwidth]{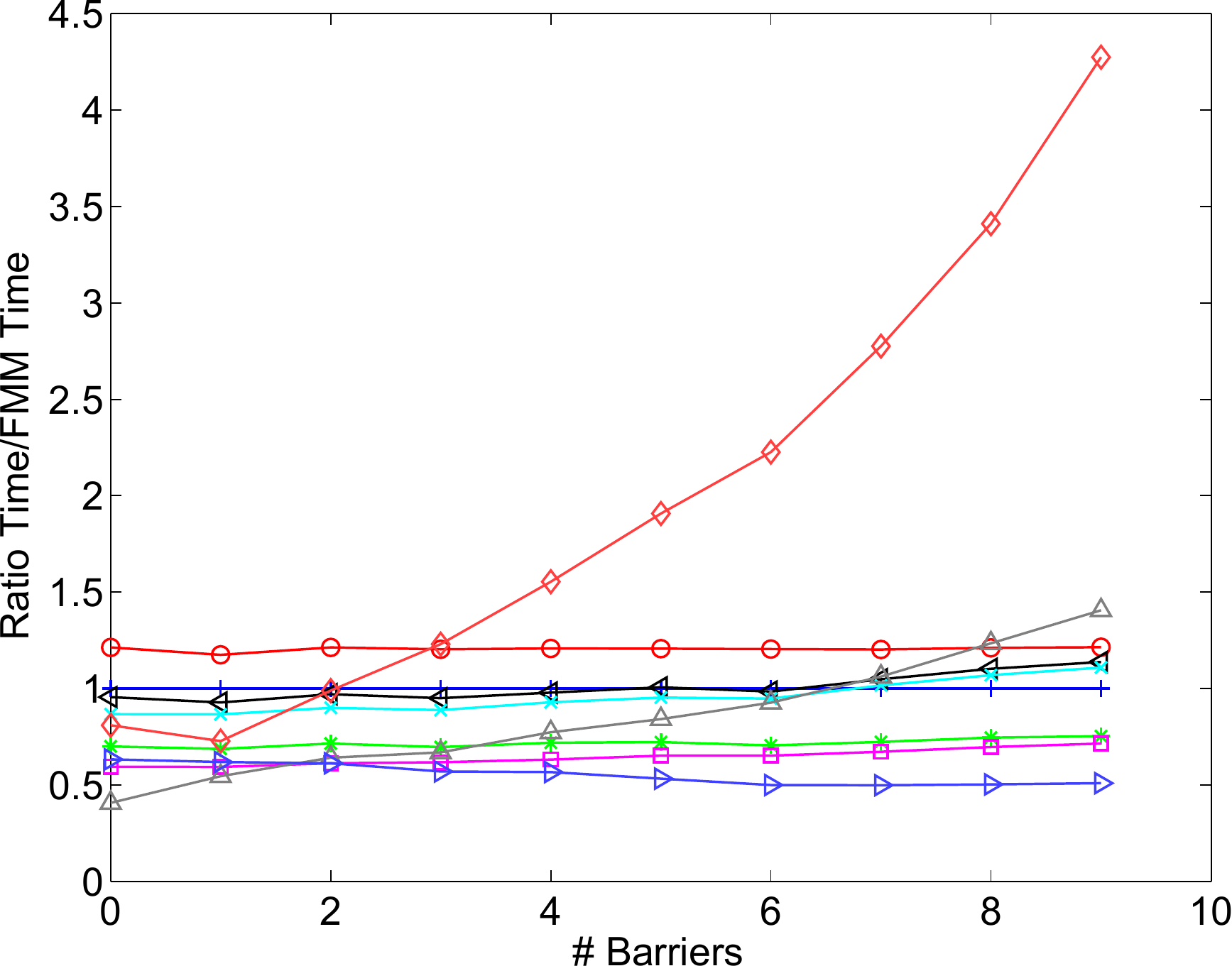}}
	\caption{Computation times and ratios for the alternating barriers experiment in 2D.}
\label{fig:res:barriers2d}
\end{figure}

\begin{figure}[ht]
	\centering
	\subfloat[Computation times.]{\includegraphics[width=0.49\columnwidth]{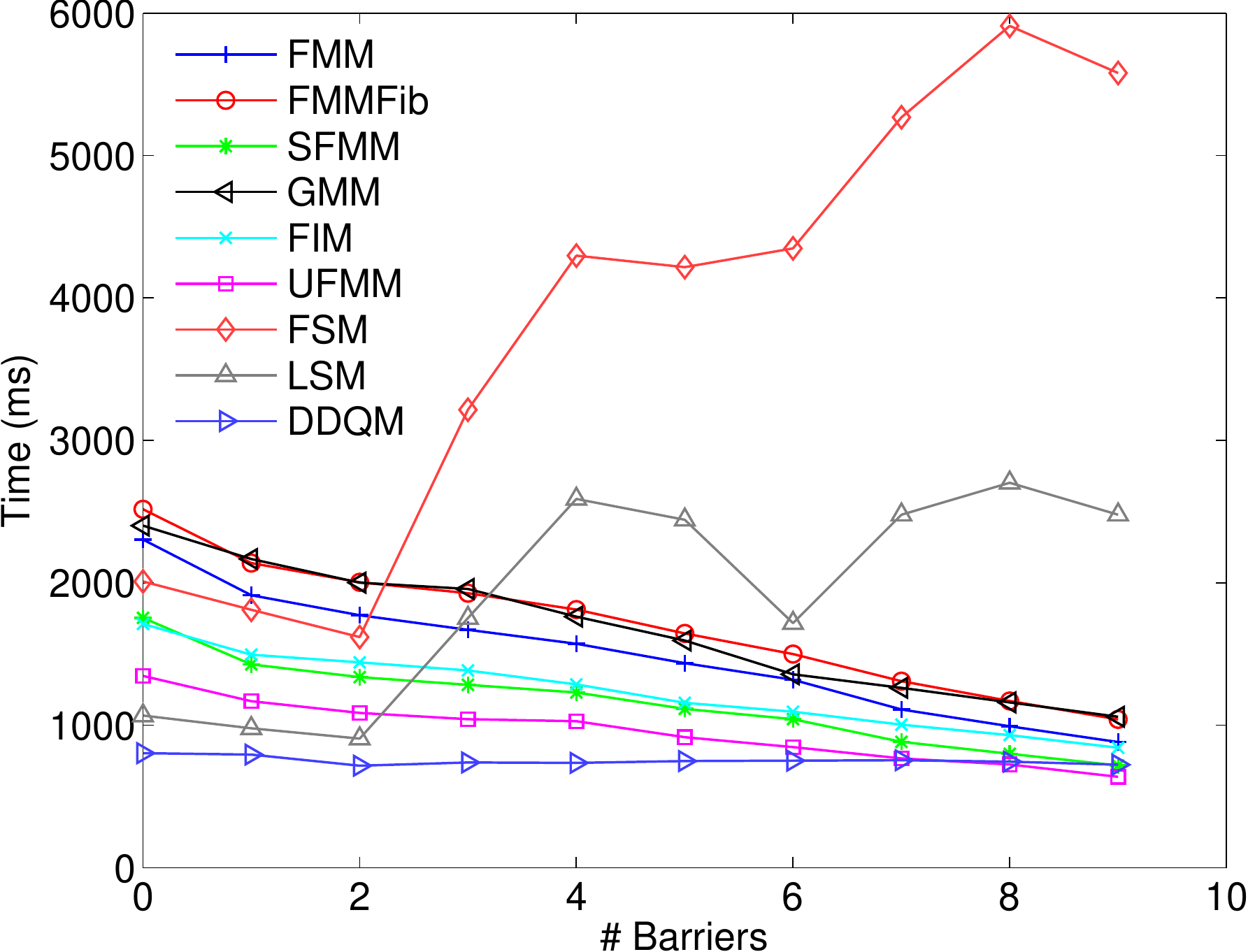}}
	\subfloat[Time ratios against FMM.]{\includegraphics[width=0.47\columnwidth]{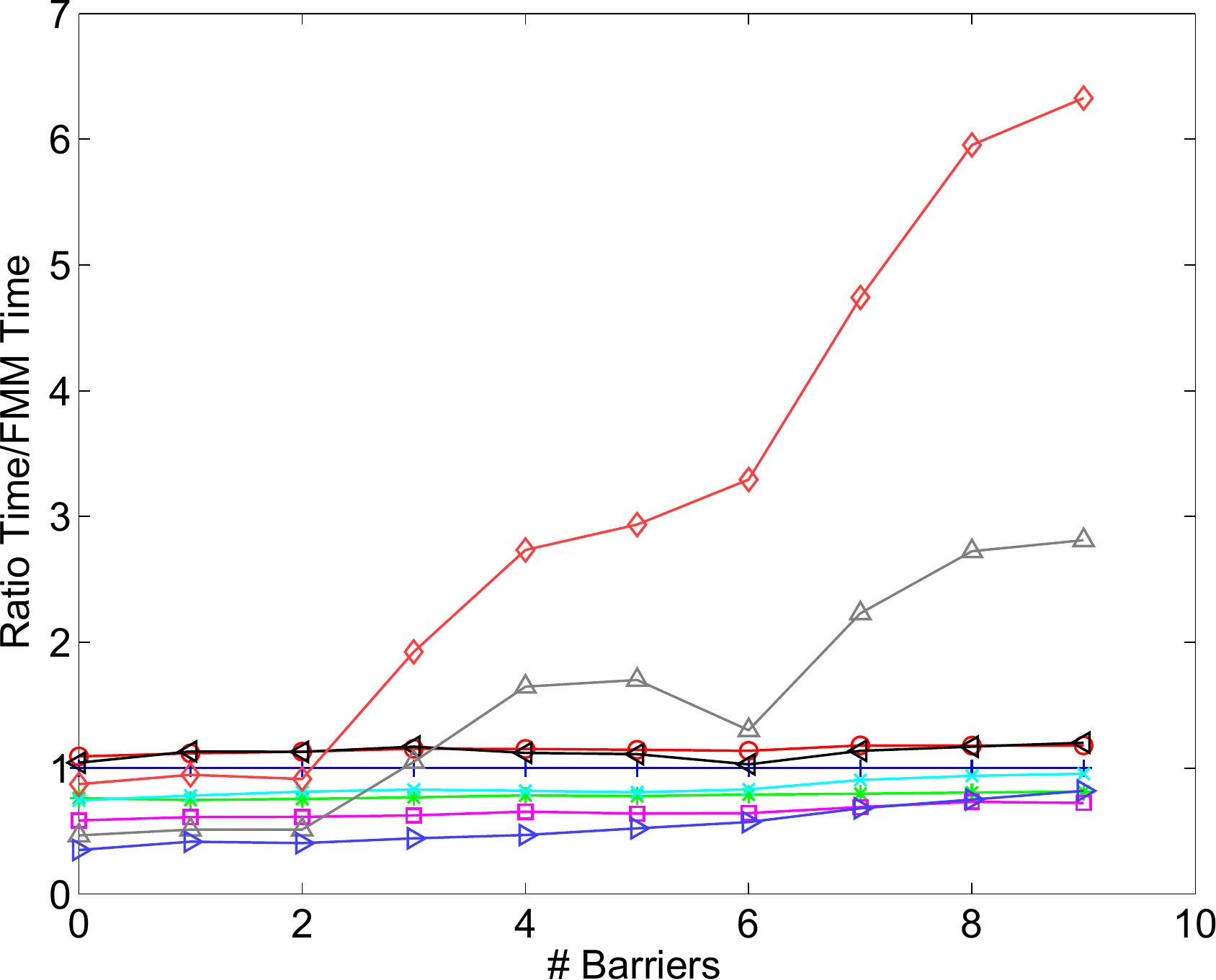}}
	\caption{Computation times and ratios for the alternating barriers experiment in 3D.}
\label{fig:res:barriers3d}
\end{figure}

An example of FMM results in some of the alternating barriers environment is shown in \cref{fig:res:barriers}, while performance results (times and ratios) for 2D and 3D are shown in \cref{fig:res:barriers2d} and \cref{fig:res:barriers3d} correspondingly.

In this case, results are very similar to the empty map experiment. However, as the number of barrier increases and the environment becomes more complex, FSM and LSM decrease their performance. However, LSM is still faster than some algorithms in most cases.

DDQM also suffers from map complexity, however it affects much less and only in 3D. Changes in characteristic directions require more sweeps to compute the final times-of-arrival map. Therefore, FSM and LSM become worse as the number of barriers increases, despite the fact that the more barriers the less cells are to evaluate. This is the reason why all other algorithms tend to lower times as the number of barriers increases.

UFMM provides again the same solution as other methods. It is the second fastest algorithm behind DDQM, closely followed by SFMM and FIM.

\subsubsection{Random velocities}
\begin{figure}[ht]
	\centering
	\subfloat[Max. velocity=30]{\includegraphics[width=0.3\columnwidth]{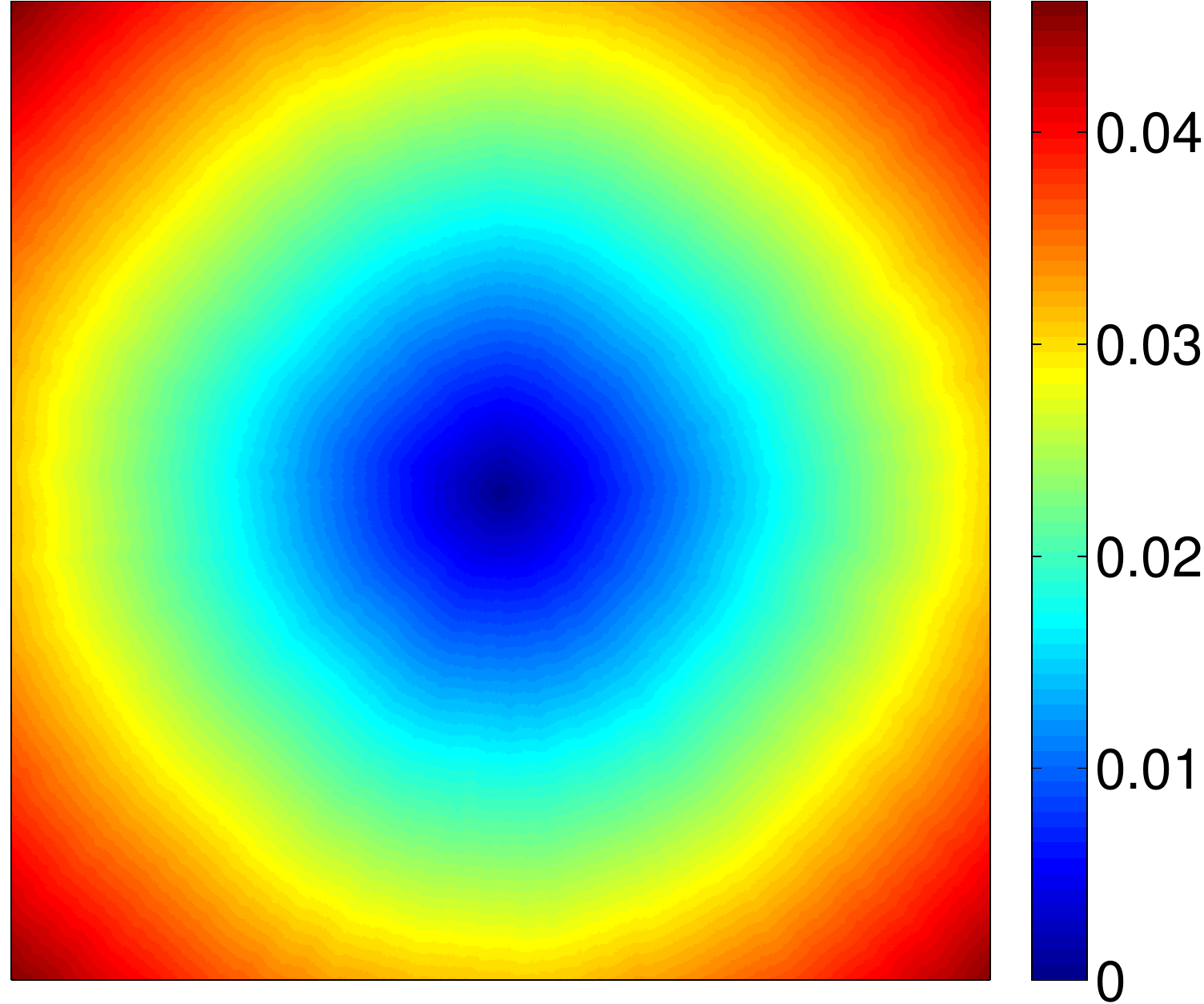}\hspace{0.2cm}}
	\subfloat[Max. velocity=60]{\includegraphics[width=0.3\columnwidth]{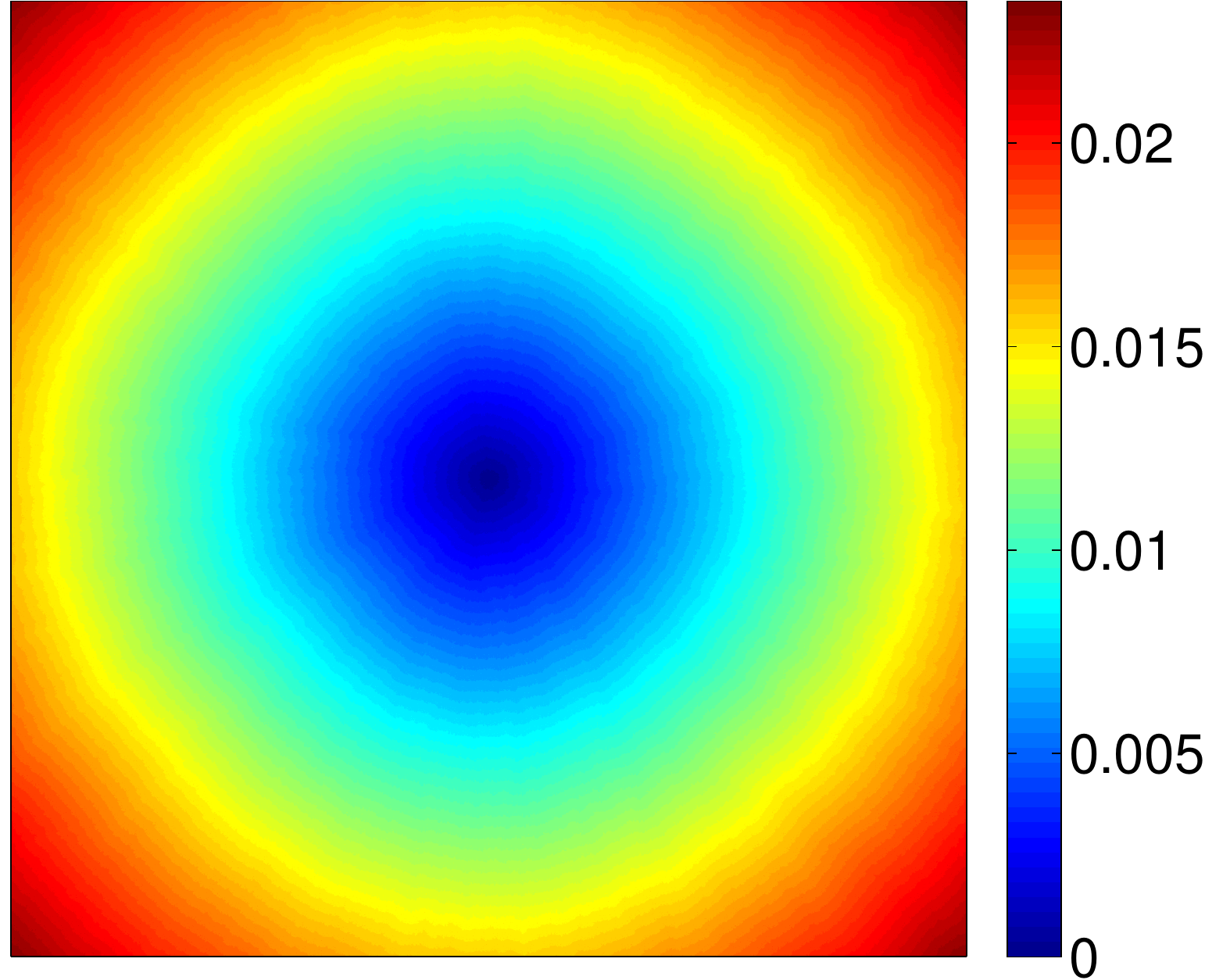}\hspace{0.2cm}}
	\subfloat[Max. velocity=100]{\includegraphics[width=0.3\columnwidth]{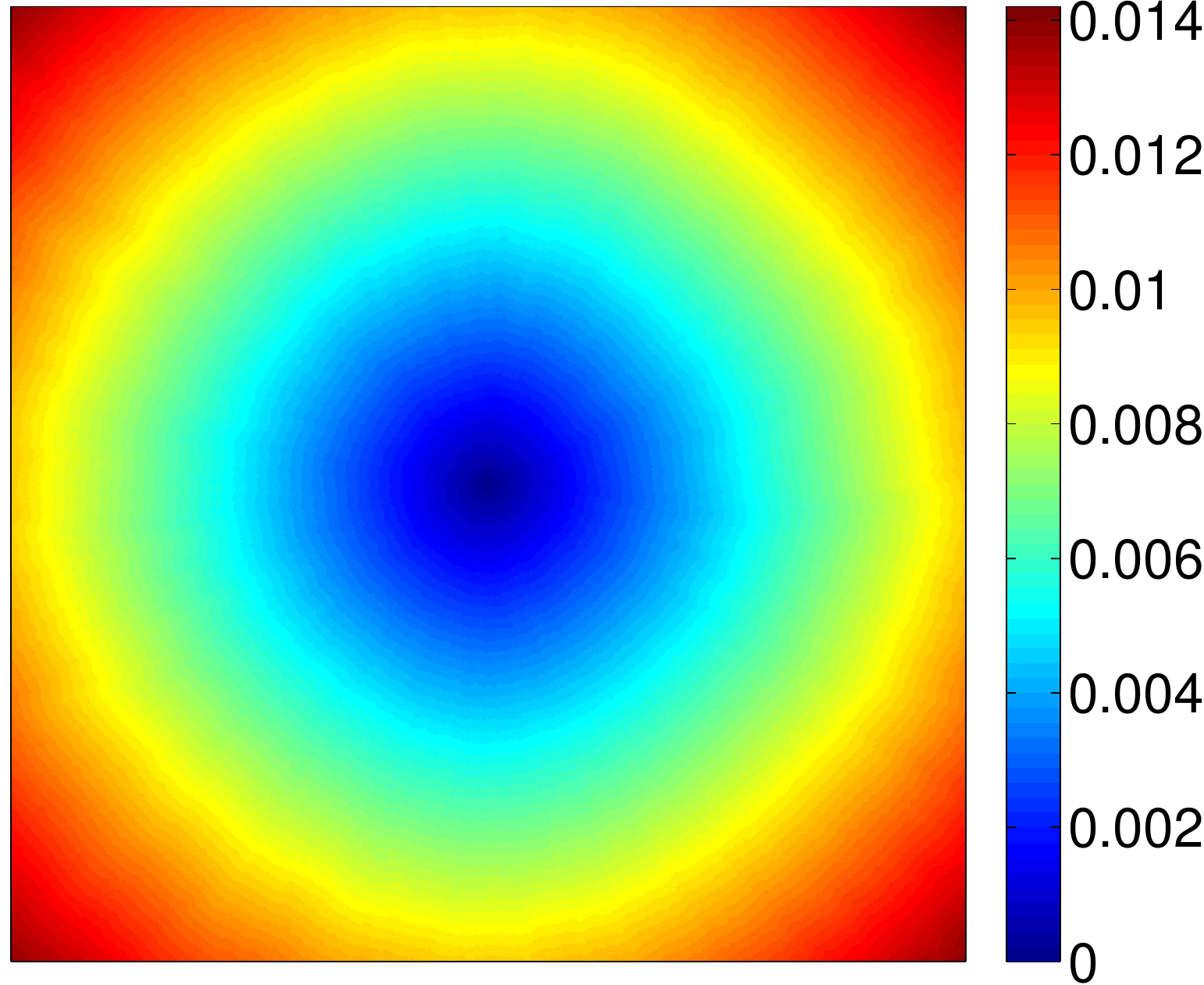}}
	\caption{Example of the resulting times-of-arrival maps applying FMM to the random velocities environment in 2D.}
\label{fig:res:random}
\end{figure}

The output of the FMM for the random velocities map is apparently close to the one of the empty map (\cref{fig:res:random}), but with the wavefronts slightly distorted because of the velocity changes.

\begin{figure}[ht]
	\centering
	\subfloat[Computation times.]{\includegraphics[width=0.47\columnwidth]{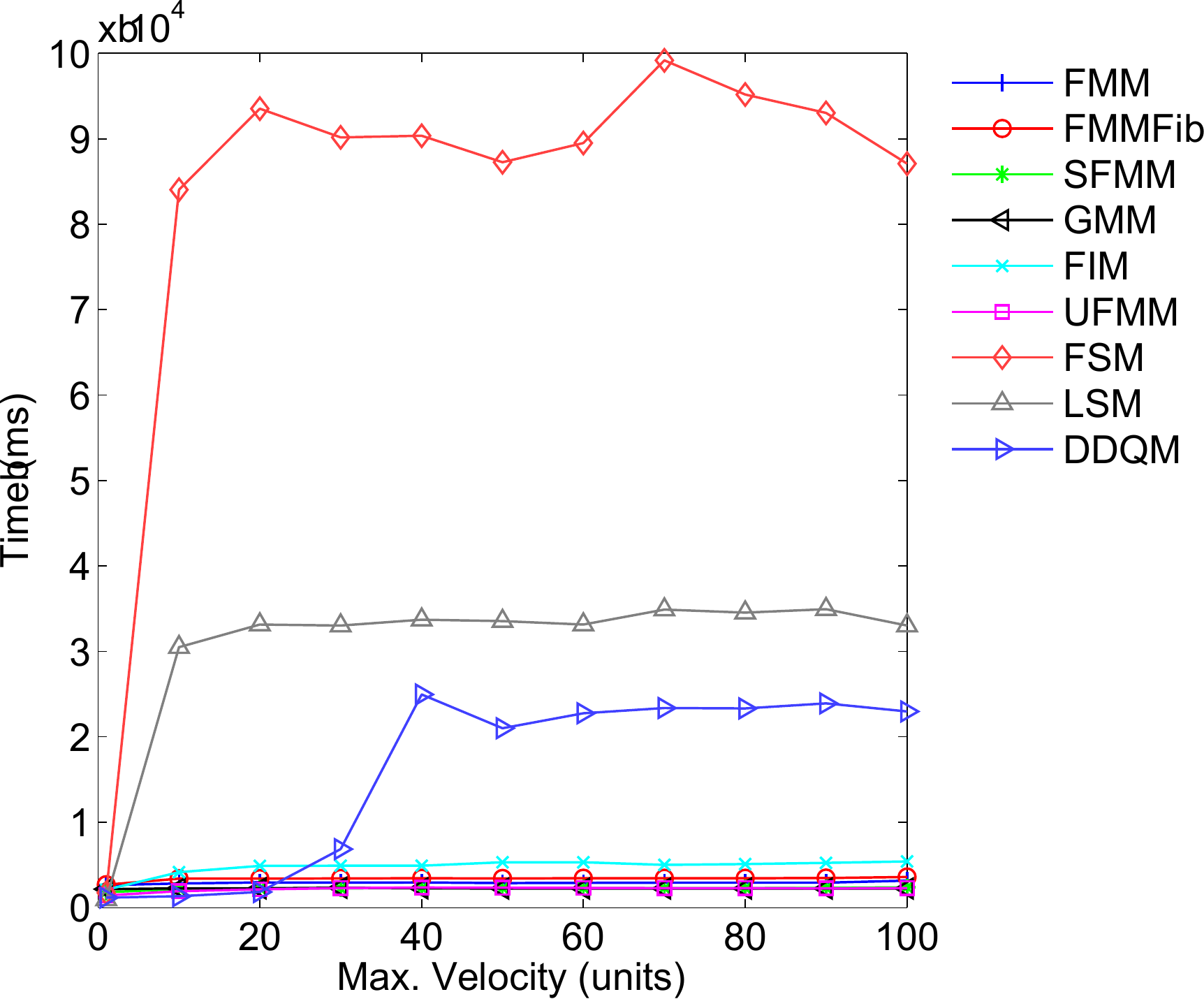}}
	\subfloat[Zoomed view.]{\includegraphics[width=0.49\columnwidth]{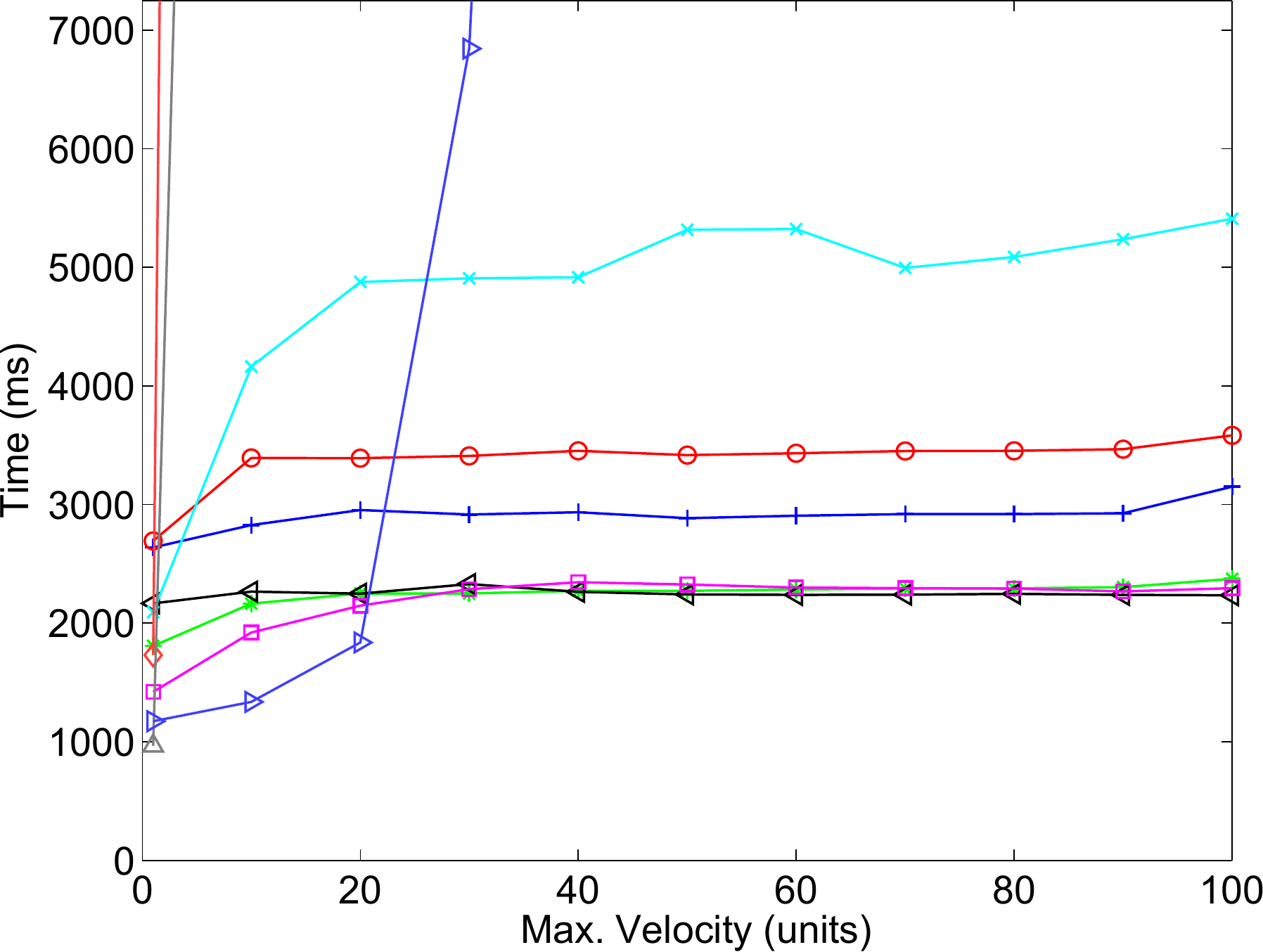}}
	\caption{Computation times and ratios for the random velocities experiment in 2D.}
\label{fig:res:random2d}
\end{figure}

\begin{figure}[ht]
	\centering
	\subfloat[Computation times.]{\includegraphics[width=0.46\columnwidth]{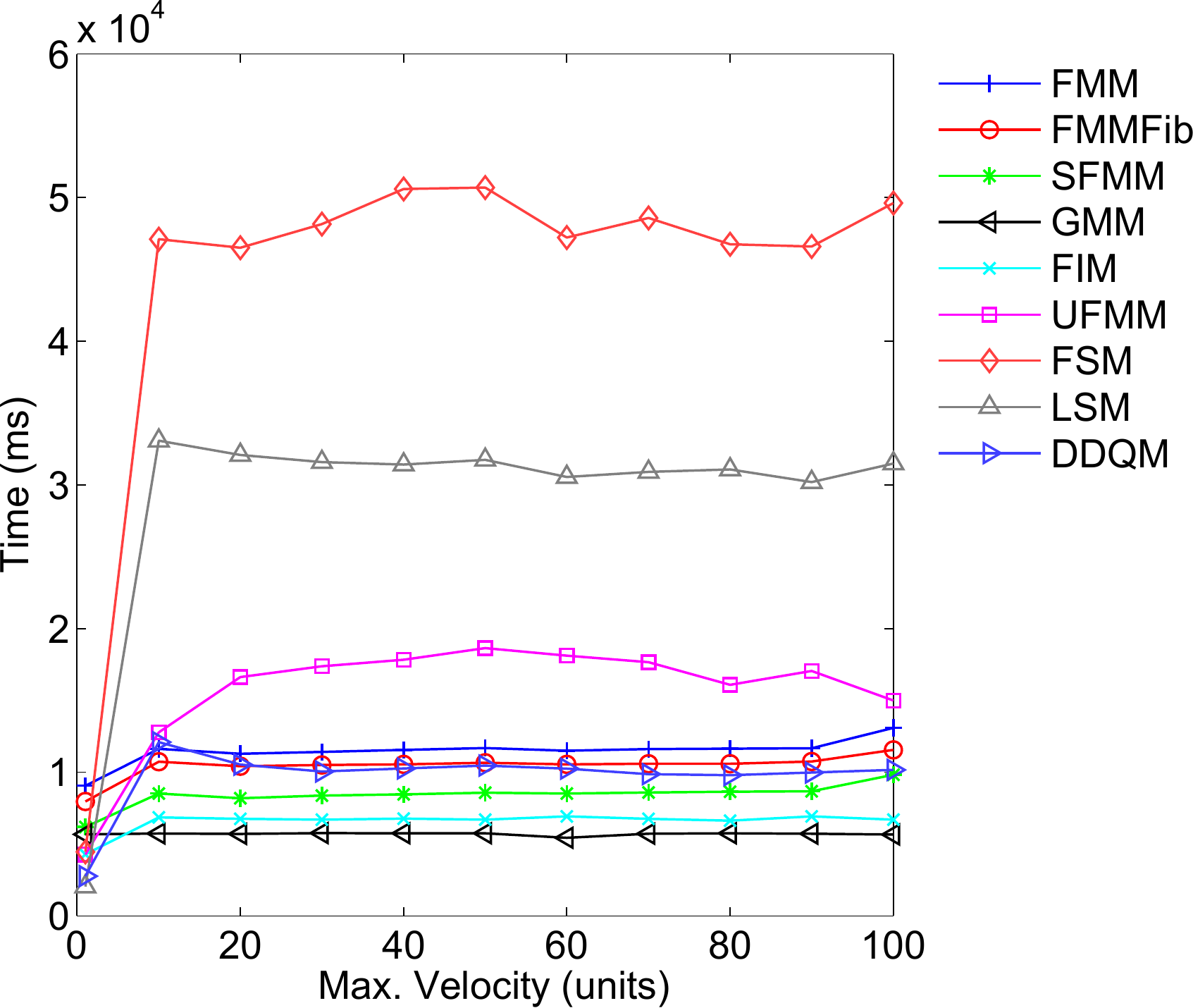}}
	\subfloat[Zoomed view.]{\includegraphics[width=0.49\columnwidth]{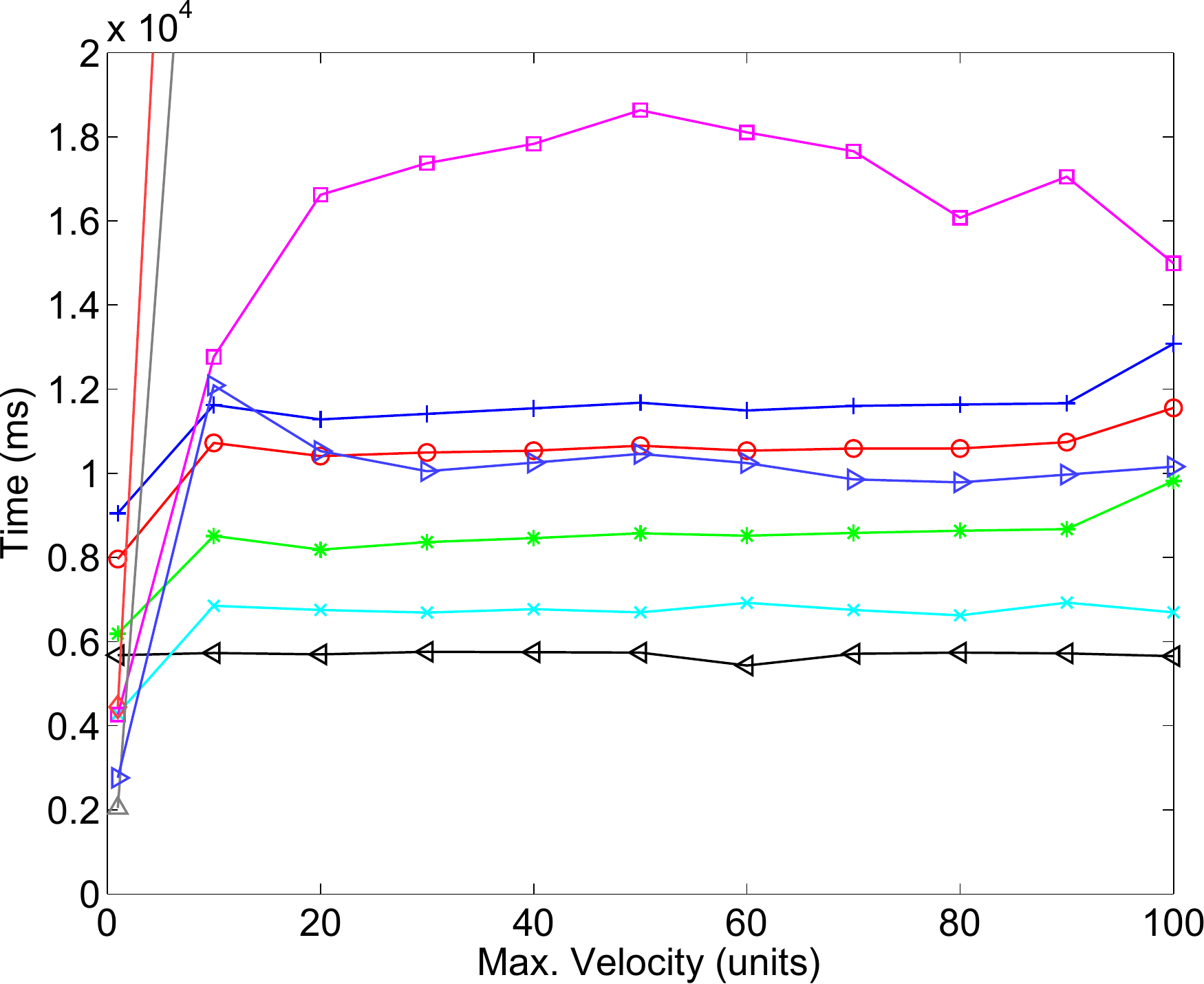}}
	\caption{Computation times and ratios for the random velocities experiment in 3D.}
\label{fig:res:random3d}
\end{figure}

\begin{figure}[ht]
	\centering
	\subfloat[Computation times.]{\includegraphics[width=0.49\columnwidth]{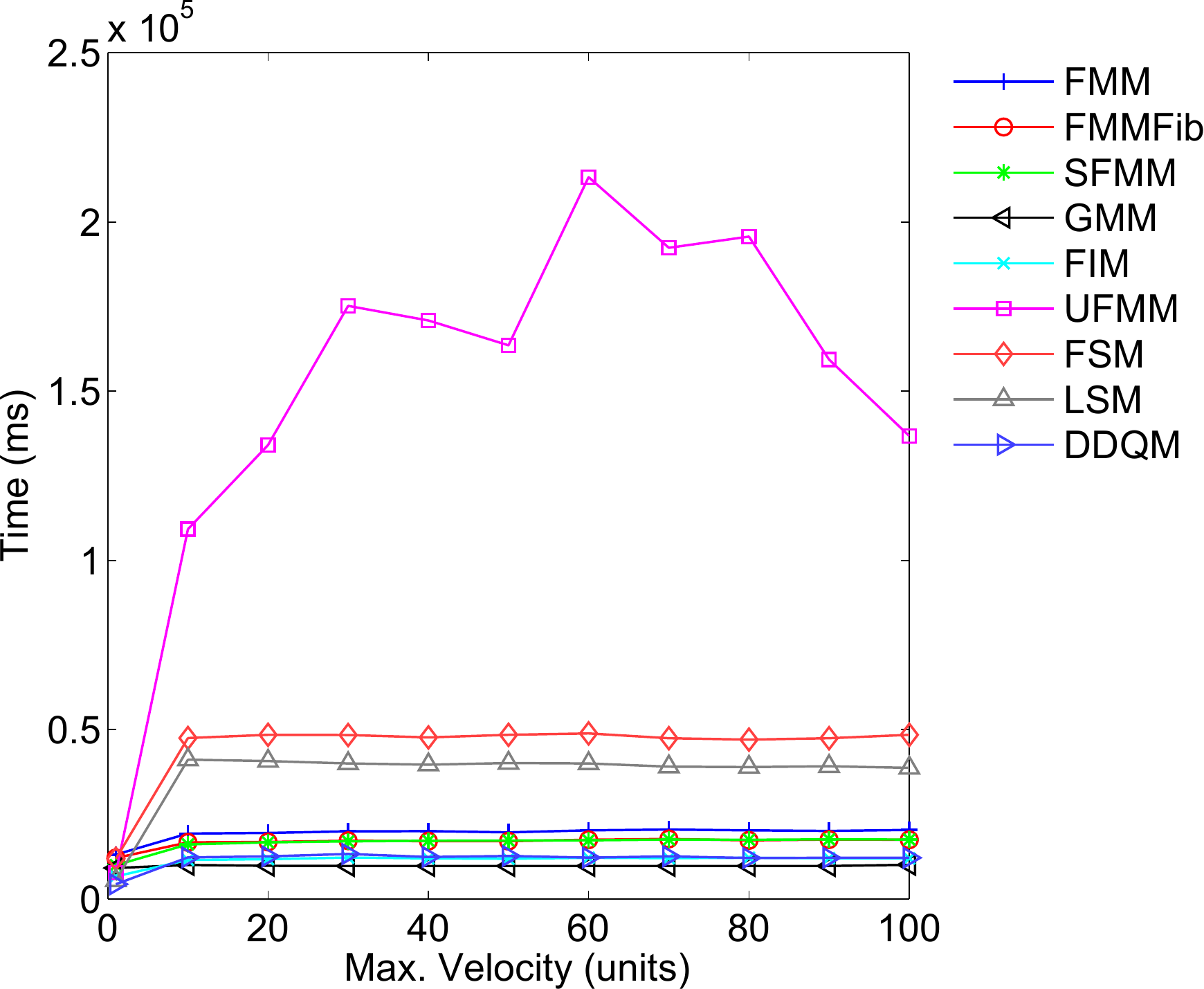}}
	\subfloat[Zoomed view.]{\includegraphics[width=0.49\columnwidth]{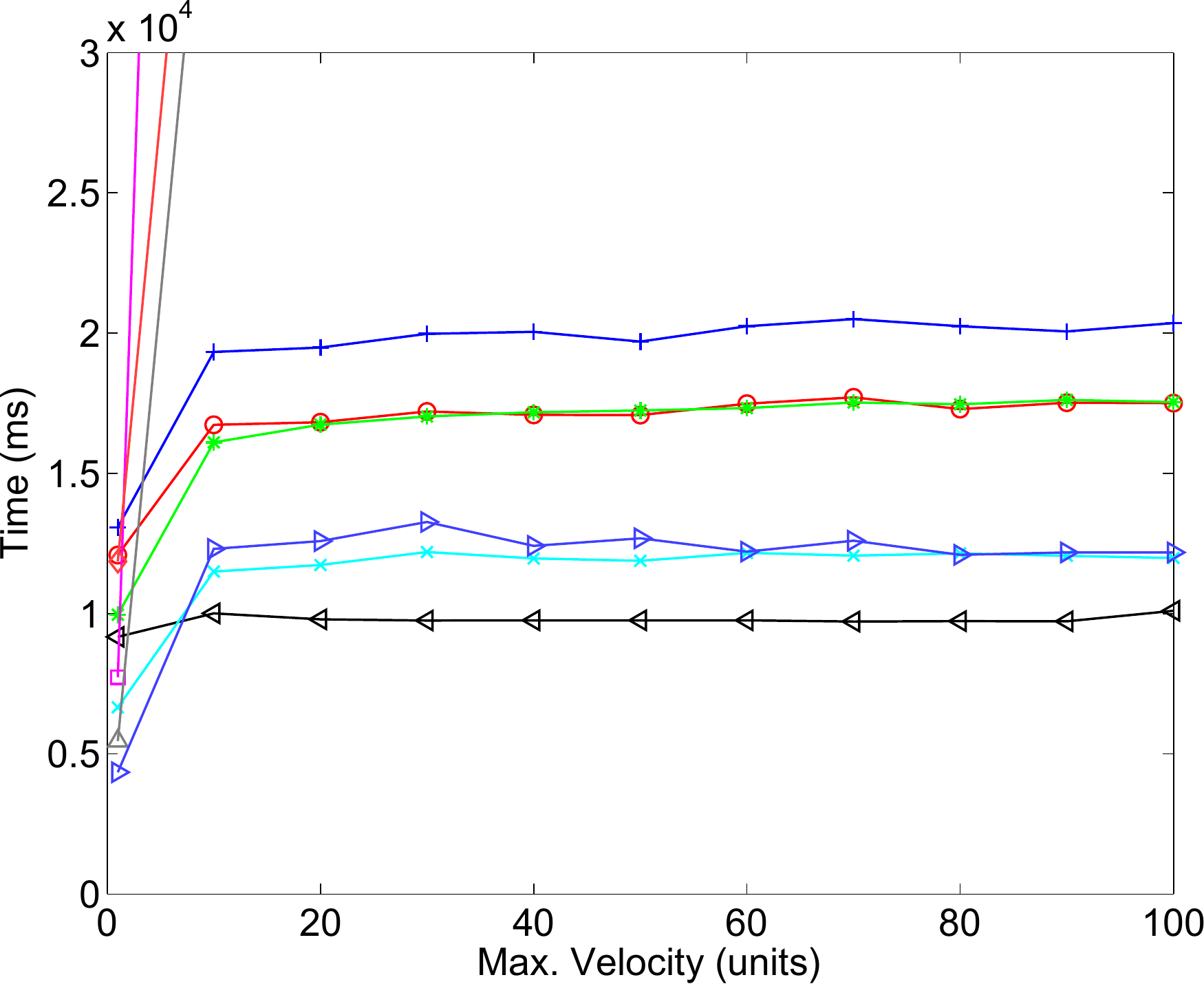}}
	\caption{Computation times and ratios for the random velocities experiment in 4D.}
\label{fig:res:random4d}
\end{figure}

Although the results on the times-of-arrival map are slightly different from the results obtained in the empty map, the performance of the algorithms are highly modified. 2D, 3D and 4D results are shown respectively in \cref{fig:res:random2d}, \cref{fig:res:random3d}, and \cref{fig:res:random4d}. In this case, raw computation times are shown together with a zoomed view of the fastest algorithms to make the analysis easier. Note that all methods become slower with non-constant velocities. The reason is that the narrow band tends to have more elements in this cases.

Some algorithms become unstable with non-uniform velocities: FSM, LSM and DDQM. DDQM is able to maintain the fastest time for slight velocity changes. But when these are sharper the double-queue threshold becomes unstable. However, for a high number of dimensions this effect vanishes (its computation time is barely modified with the number of dimensions) and DDQM becomes the fastest algorithms together with FIM and GMM.

SFMM provides one of the best performances across all dimensions. And FIM becomes relatively faster in 3D and 4D. However, in most of the cases GMM is the fastest algorithm. FIM requires now multiple iterations to converge to the solution while GMM guarantees convergence with only 2 iterations.

Finally, UFMM requires special attention as it does not return the same solution than the other Fast Methods. Its performance is highly affected by the number of dimensions. The main reason is the election of the parameters: they were experimentally chosen to optimize 2D performance. However, these parameters are no longer useful in other number of dimensions. We consider UFMM parameter tuning to be a complex task.

\cref{tab:random:ufmm_err} summarizes the largest errors for this experiment. As the number of dimension is increased, the error decreases exponentially while the computation time increases exponentially. Therefore, by properly tuning parameters for 3D and 4D better times could be achieved while keeping a negligible error in most practical cases.

\begin{table}[ht]
\caption{Largest $L_1$ and $L_\infty$ errors for UFMM in the random velocities experiment.}
\begin{tabular}{r|ccc}
\multicolumn{1}{l|}{} & 2D                & 3D                      & 4D  \\ \hline
$L_1$             	  & $10^{-3}$         & $1.3\cdot 10^{-10}$			&$6.9\cdot 10^{-12}$ \\
$L_\infty$            & $4.8\cdot10^{-3}$	& $10^{-6}$	              & $10^{-7}$ \\ \hline
\end{tabular}
\label{tab:random:ufmm_err}
\end{table}

\subsubsection{Checkerboard}
\begin{figure}[ht]
	\centering
	\subfloat[Max. velocity=30]{\includegraphics[width=0.3\columnwidth]{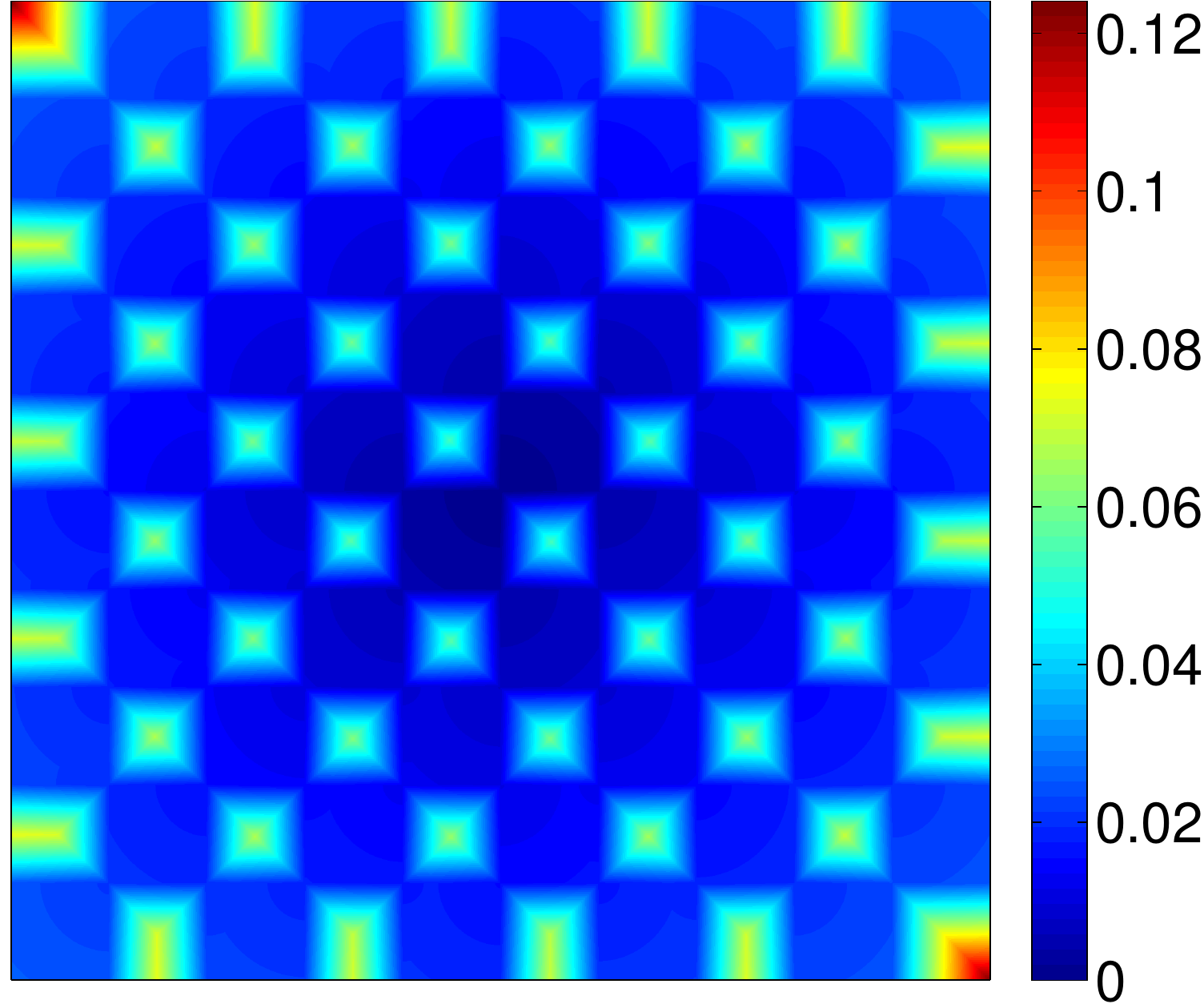}\hspace{0.2cm}}
	\subfloat[Max. velocity=60]{\includegraphics[width=0.3\columnwidth]{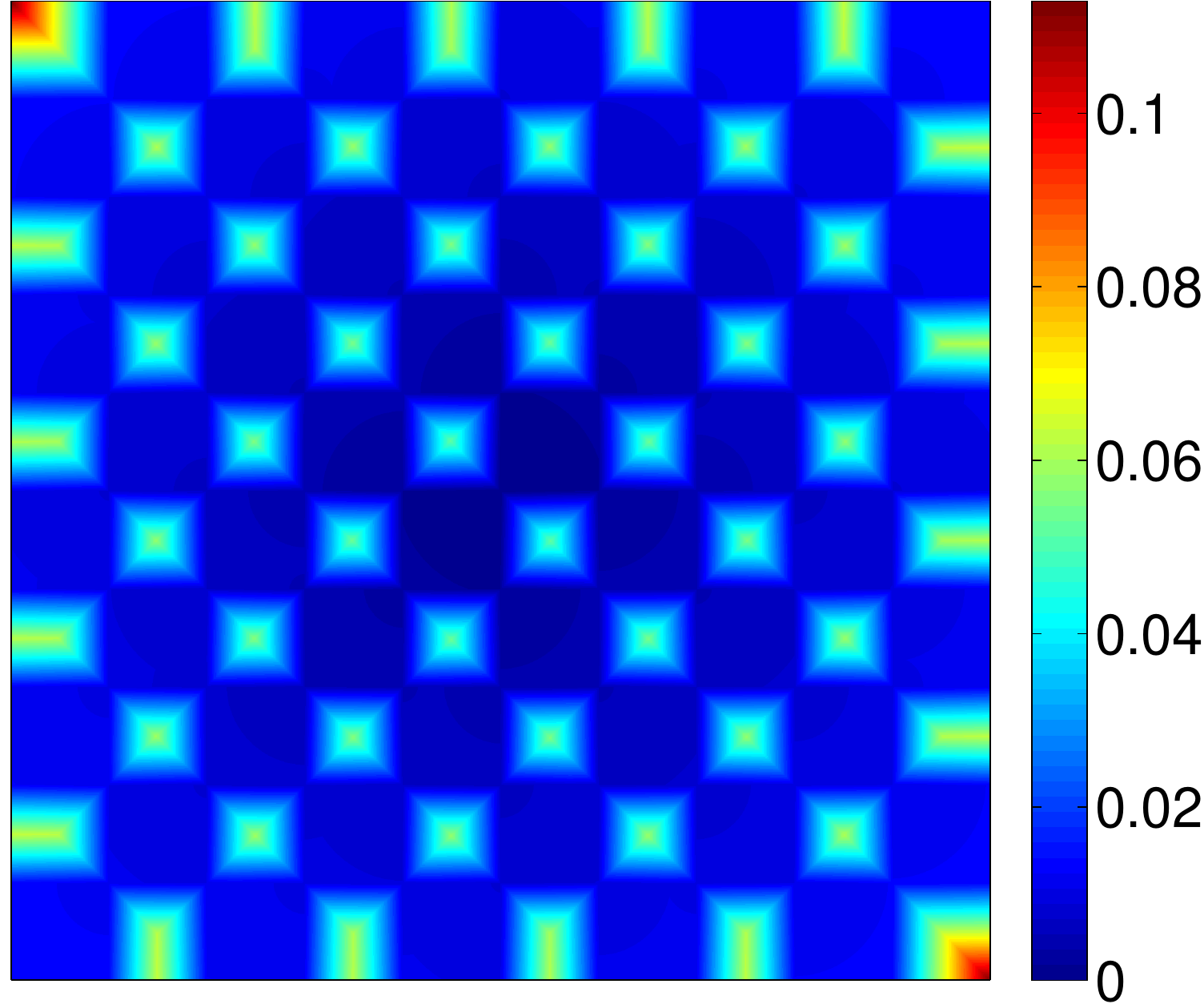}\hspace{0.2cm}}
	\subfloat[Max. velocity=100]{\includegraphics[width=0.3\columnwidth]{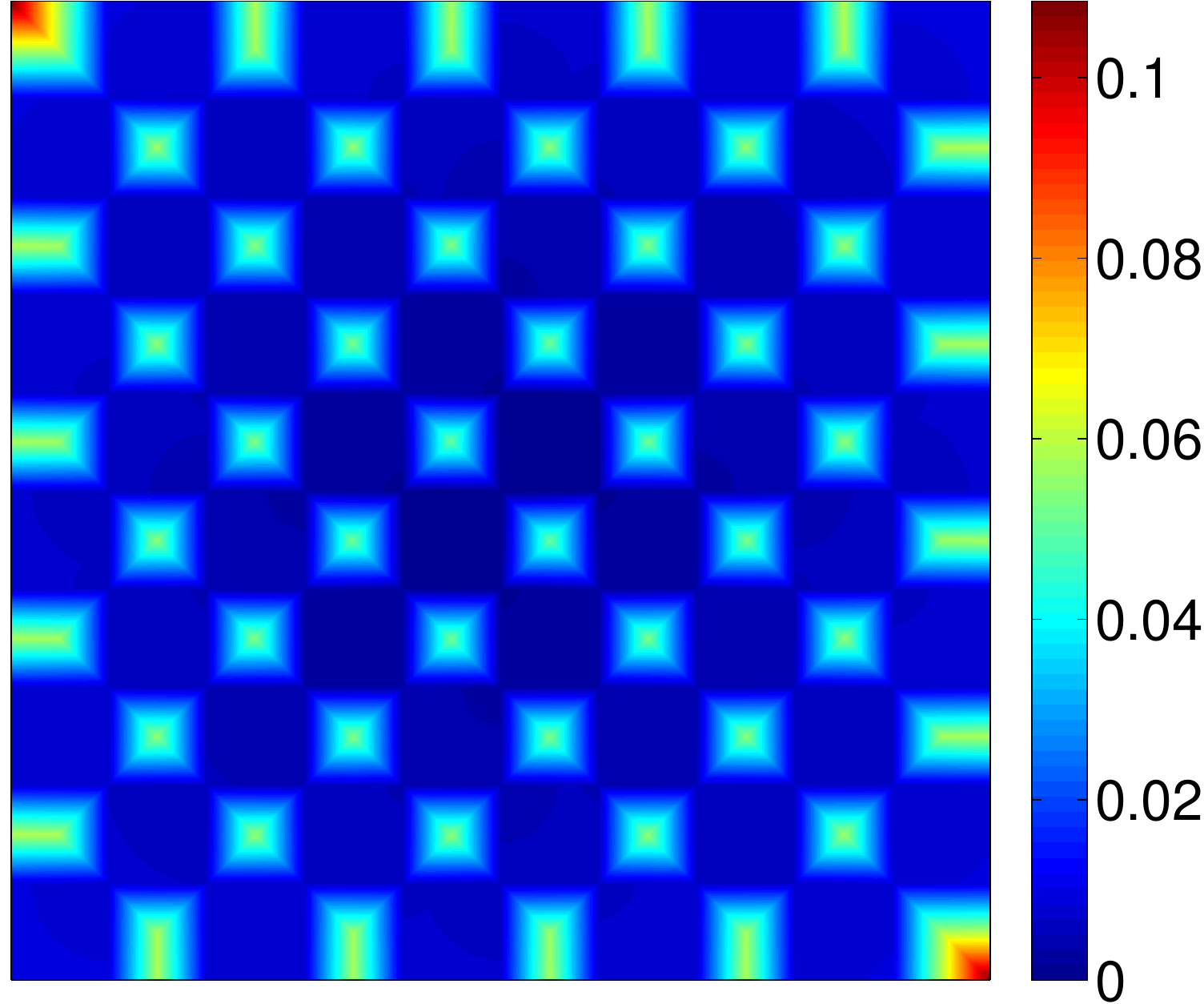}}
	\caption{Example of the resulting times-of-arrival maps applying FMM to the checkerboard environment in 2D.}
\label{fig:res:checker}
\end{figure}

\begin{figure}[ht]
	\centering
		\subfloat[Computation times.]{\includegraphics[width=0.47\columnwidth]{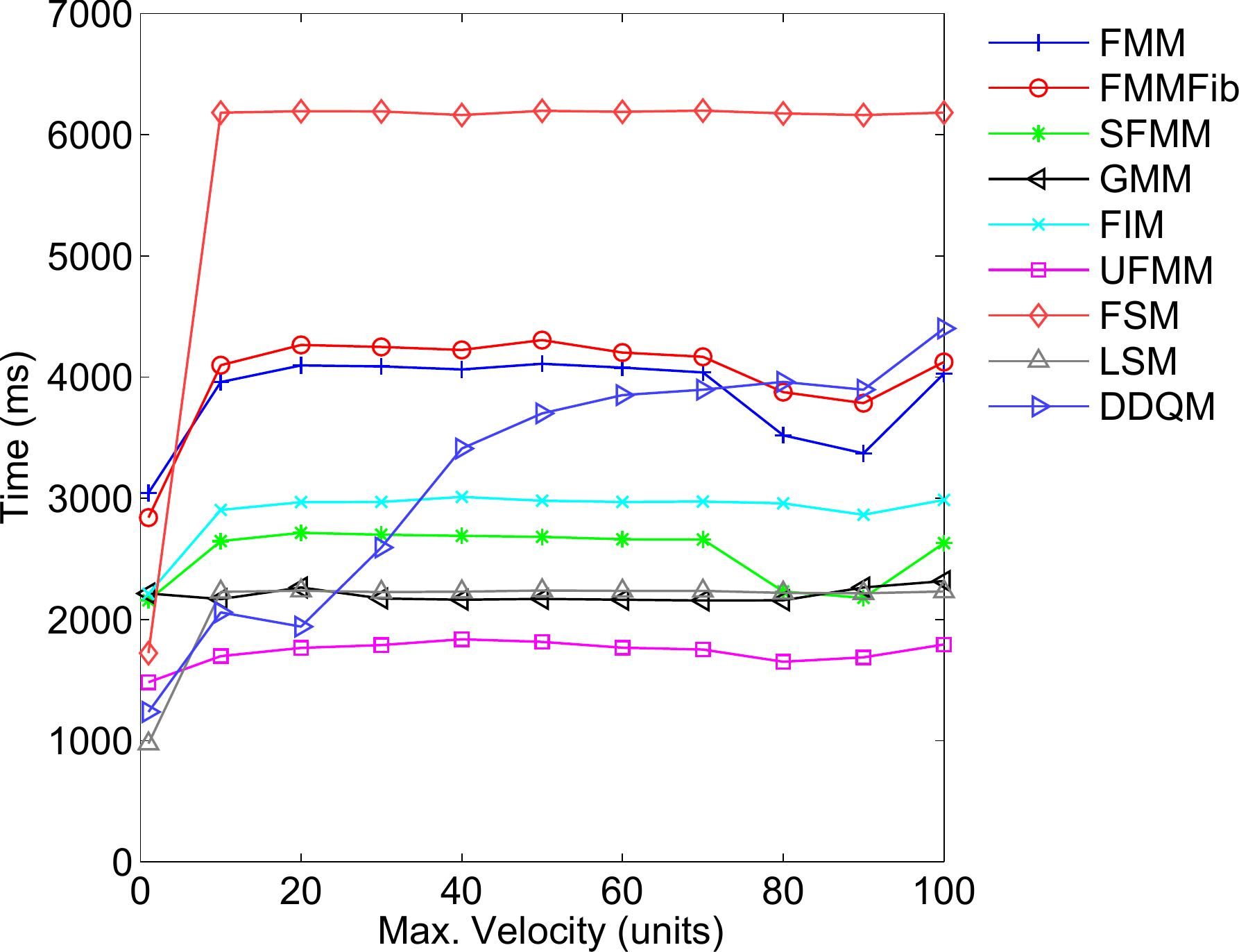}}
		\subfloat[Time ratios against FMM.]{\includegraphics[width=0.47\columnwidth]{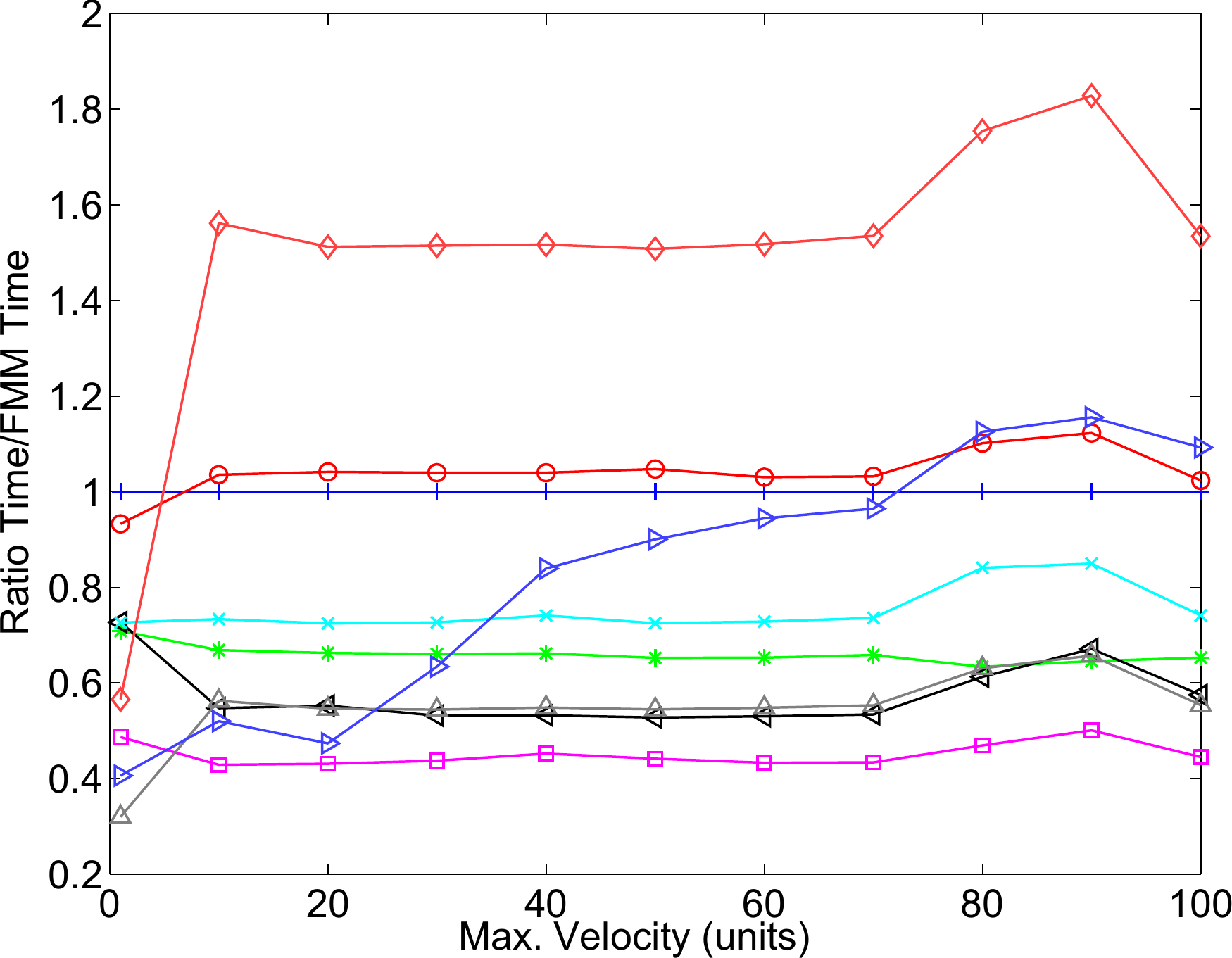}}
	\caption{Computation times for the checkerboard experiment in 2D.}
\label{fig:res:checker2d}
\end{figure}

\begin{figure}[ht]
	\centering
		\subfloat[Computation times.]{\includegraphics[width=0.47\columnwidth]{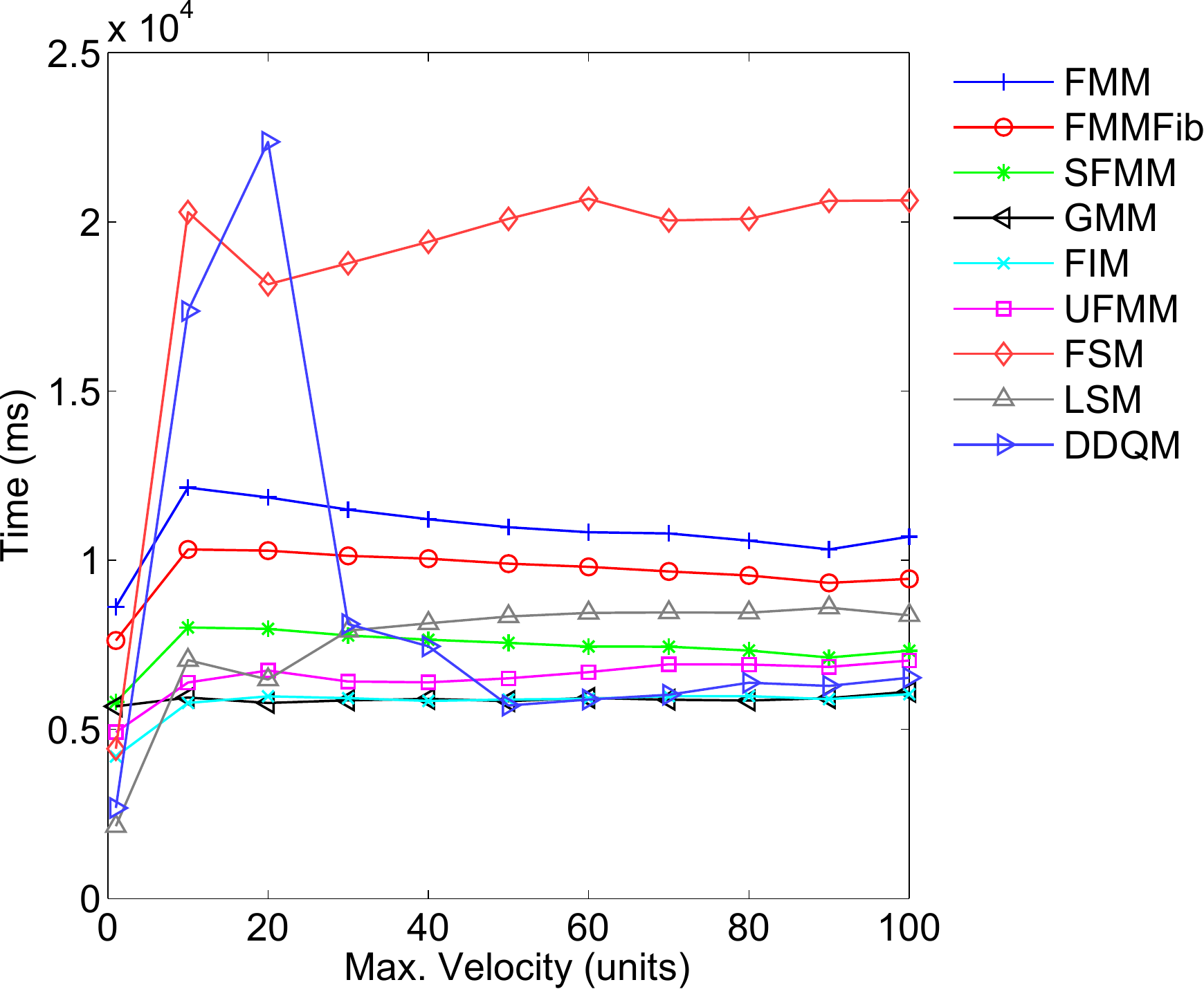}}
		\subfloat[Time ratios against FMM.]{\includegraphics[width=0.47\columnwidth]{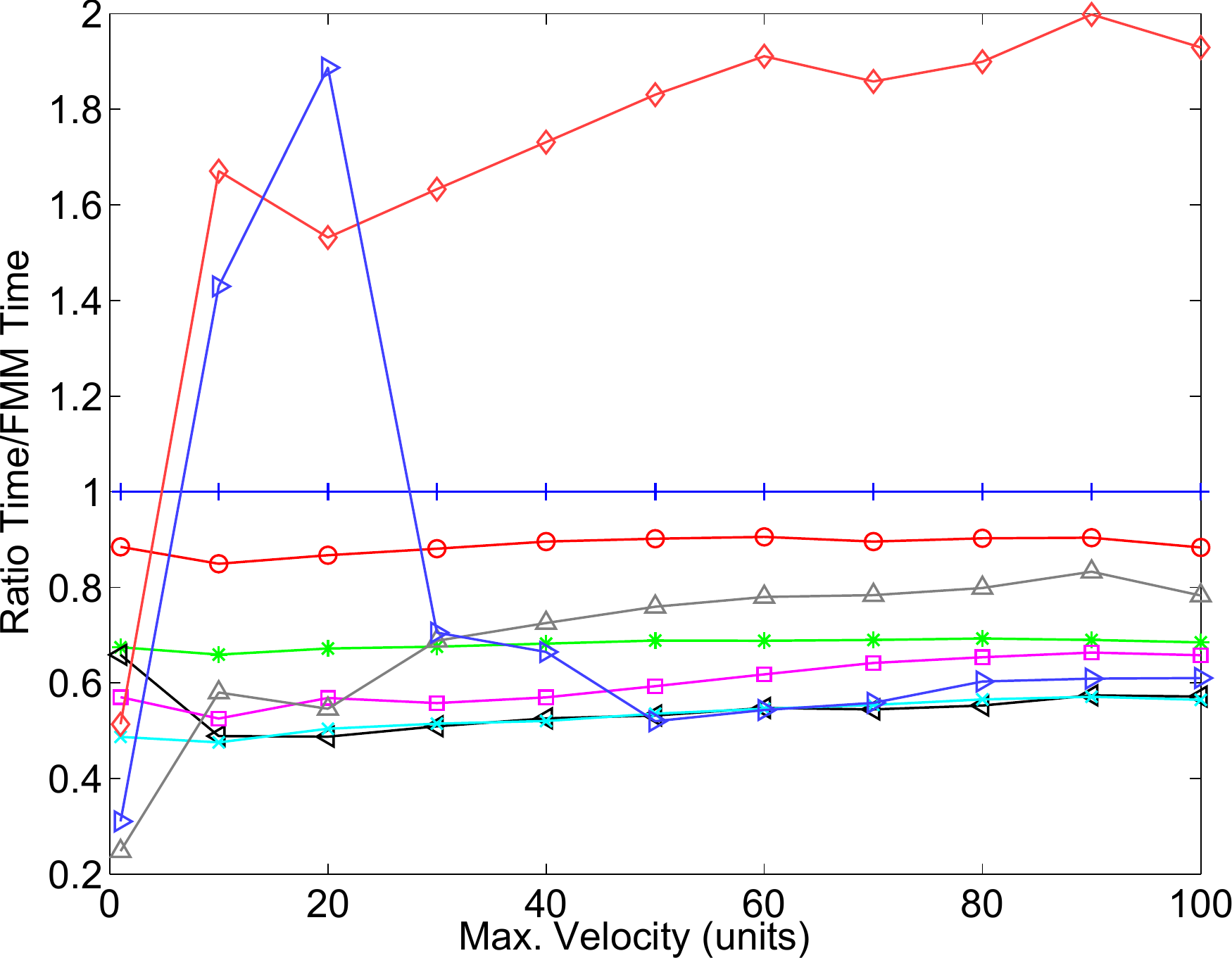}}
	\caption{Computation times for the checkerboard experiment in 2D.}
\label{fig:res:checker3d}
\end{figure}

\begin{figure}[ht]
	\centering
		\subfloat[Computation times.]{\includegraphics[width=0.47\columnwidth]{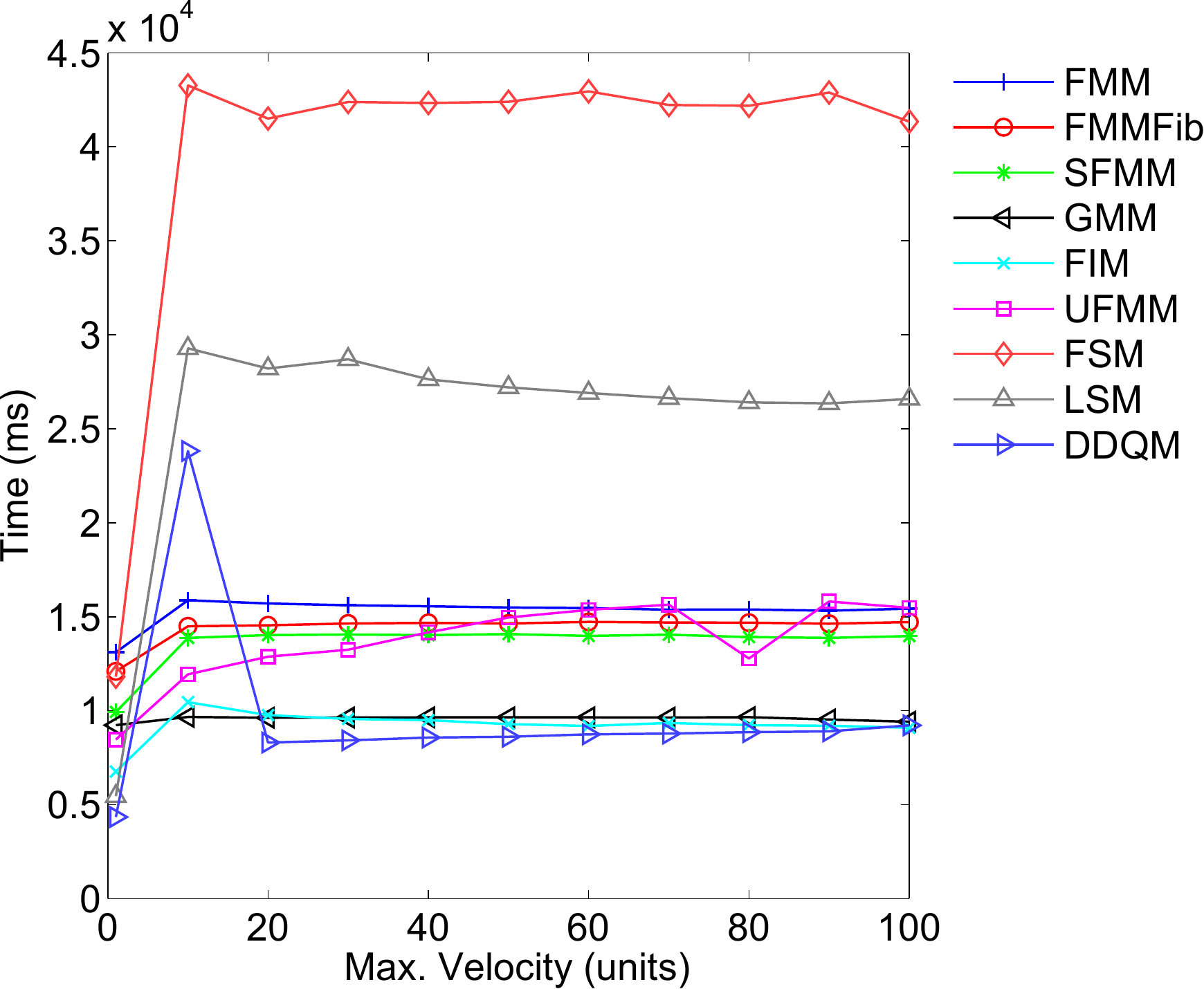}}
		\subfloat[Time ratios against FMM.]{\includegraphics[width=0.47\columnwidth]{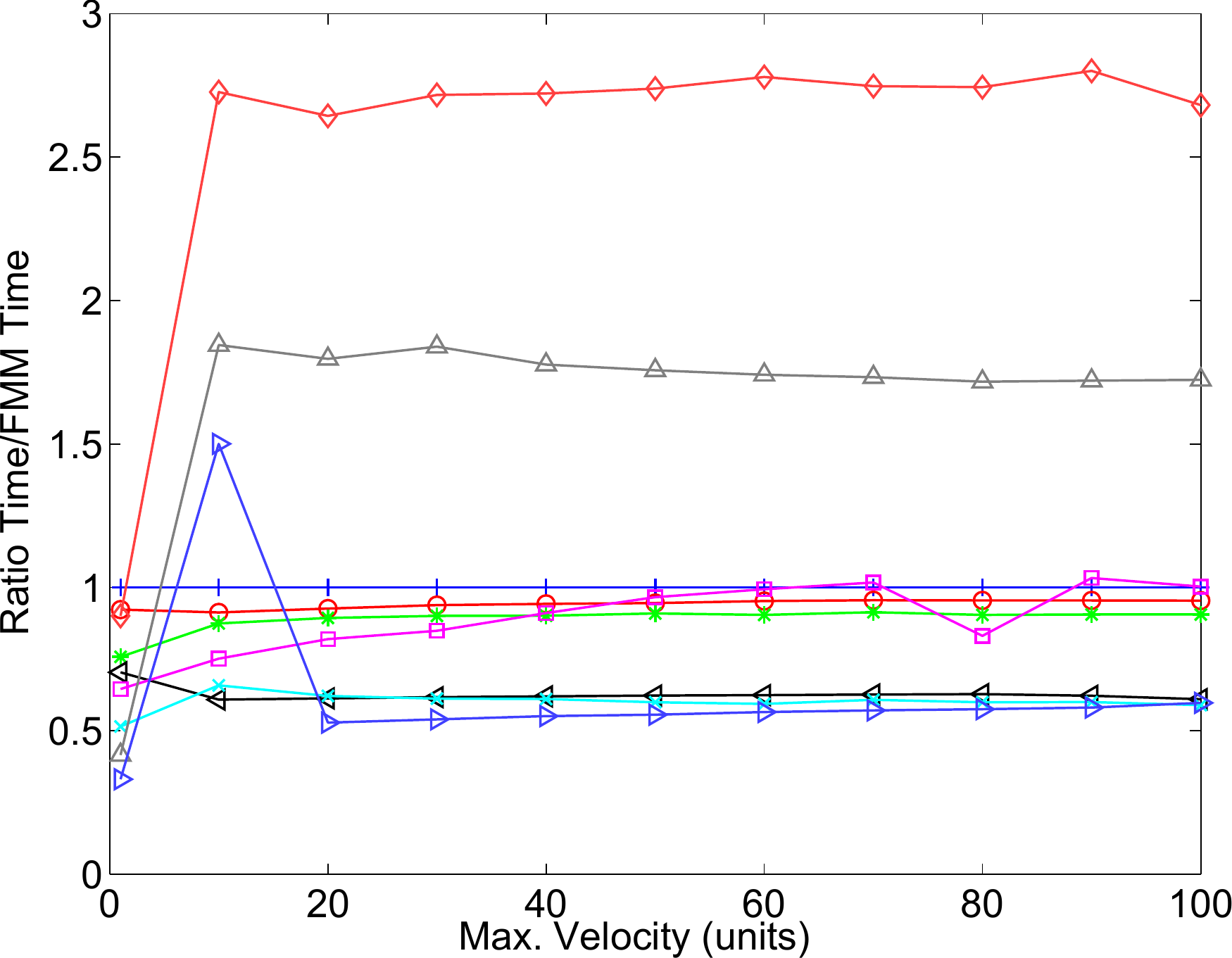}}
	\caption{Computation times for the checkerboard experiment in 4D.}
\label{fig:res:checker4d}
\end{figure}

Finally, different times-of-arrival map returned by FMM applied to the checkerboard map are shown in \cref{fig:res:checker}. Numerical results of the computation times are included in \cref{fig:res:checker2d} for 2D, \cref{fig:res:checker3d} for 3D, and \cref{fig:res:checker4d} for 4D.

The results are relatively close to the random velocities experiment. However, the differences for FSM, LSM and DDQM are much smaller. In fact, DDQM presents a poor performance in 2D, but in 3D and 4D it becomes the fastest algorithm for higher velocity modifications.

GMM and FIM have very similar results for 3D and 4D. Since the environment is structured, FIM does not require too many iterations to compute the final value.

UFMM errors are shown in \cref{tab:checker:ufmm_err}. In this case, UFMM becomes worse with the number of dimensions but it is among the fastest algorithms in all cases. 

The differences between this experiments and random velocities is that the environment presents a well-defined structure, and locally acts as a constant velocity environment, as in empty map experiment.

\begin{table}[ht]
\caption{Largest $L_1$ and $L_\infty$ errors for UFMM in the checkerboard experiment.}
\begin{tabular}{r|ccc}
\multicolumn{1}{l|}{} & 2D                & 3D                & 4D  \\ \hline
$L_1$             	  & $1.7\cdot10^{-7}$ &$1.2\cdot 10^{-9}$ &$1.9\cdot 10^{-10}$ \\
$L_\infty$            & $2.5\cdot10^{-6}$	& $5\cdot10^{-7}$	  & $10^{-6}$ \\ \hline
\end{tabular}
\label{tab:checker:ufmm_err}
\end{table}

\section{Discussion}
\label{sec:disc}
With the four experiments designed, the main characteristics of the Fast Methods have been shown. Any other environment can be thought as a combination of free space with obstacles and high-frequency or low-frequency velocity changes of different magnitude. Note also that the grid sizes covered by the experiments vary from extremely small to extremely big grids. In practice, it is hard to find applications requiring more than 16 million cells.

FIM results can be speeded up for non-constant velocity problems if larger errors are allowed. UFMM can be probably improved as well. However, our experience is that the configuration of its parameters is complex and requires a deep knowledge of the environment to be applied on.

Several conclusions can be extracted from the conducted experiments: 1) There is no practical reasons to use FMM or FMMFib as SFMM is faster in all cases with the same behaviour as its counterparts. 2) If a sweep-based method is required, LSM should be always chosen, as it greatly outperforms FSM. 3) In problems with constant velocity DDQM should be chosen, as it has shown the best performance for the empty map and the alternating barriers environments. 4) For variable velocities, but simple scenarios, GMM is the algorithm to choose. 5) UFMM is hard to tune and the result has errors. Also, it has been outperformed in most of the cases by DDQM in constant velocity scenarios, or by SFMM or FIM in experiments with variable velocities. 6) There is not a clear winner for complex scenarios with variable velocity. UFMM can perform well in all cases if tuned properly. Otherwise, SFMM is a safe choice, specially in cases where there is not too much information about the environment.

If a goal point is selected, cost-to-go heuristics could be applied \cite{Valero13}, and thus enormously affect the results. Heuristics for FMM, FMMFib and SFMM are straightforward. They would improve the results in most of the cases. They could also be applied to UFMM. However, it is not clear if they can be applied to other Fast Methods. It is also of interest the solution of anisotropic problems \cite{Clement07}, solved only by FMM-based methods.

\section{Conclusions}
\label{sec:conc}
Along this paper we have introduced the main Fast Methods in a common mathematical framework, focusing on a practical point of view.

With this work we aimed at closing the discussion about which method should be used in which case, as this is the first exhaustive comparison of the main Fast Methods (up to our knowledge).

The code is publicly available as well as the automatic benchmark programs. This code has been deeply tested and it can serve as a base for future algorithm design, as it provides all the tools required to easily implement and compare novel Fast Methods.

The future work focuses in 3 different aspects: develop the analogous work for parallel Fast Methods \cite{Detrixhe13}, study the application of these methods to anisotropic problems and also to the new Fast Marching-based solutions focused on path planning applications \cite{Janson14,Yershov15}. Finally, the combination of UFMM and SFMM seems straightforward and it would presumably outperform both algorithms.

\begin{acknowledgements}
Authors want to gratefully acknowledge the contribution of Pablo Gely Mu\~noz to the GMM,  FIM and UFMM implementations, and to Adam Chacon for the interesting discussions and suggestions towards improving the work.
\end{acknowledgements}

\bibliographystyle{./styles/spmpsci}      
\bibliography{./styles/fmmbib}   

\begin{thebibliography}{10}
\providecommand{\url}[1]{{#1}}
\providecommand{\urlprefix}{URL }
\expandafter\ifx\csname urlstyle\endcsname\relax
  \providecommand{\doi}[1]{DOI~\discretionary{}{}{}#1}\else
  \providecommand{\doi}{DOI~\discretionary{}{}{}\begingroup
  \urlstyle{rm}\Url}\fi

\bibitem{Ahmed11}
Ahmed, S., Bak, S., McLaughlin, J., Renzi, D.: {A third order accurate Fast
  Marching Method for the Eikonal equation in two dimensions}.
\newblock SIAM Journal on Scientific Computing \textbf{33}(5), 2402--2420
  (2011)

\bibitem{Alzaben15}
Al~Zaben, N., Madusanka, N., Al~Shdefat, A., Choi, H.: Comparison of active
  contour and fast marching methods of hippocampus segmentation.
\newblock In: 6th International Conference on Information and Communication
  Systems, pp. 106--110 (2015)

\bibitem{Bak10}
Bak, S., McLaughlin, J., Renzi, D.: Some improvements for the fast sweeping
  method.
\newblock SIAM Journal on Scientific Computing \textbf{32}(5), 2853--2874
  (2010)

\bibitem{Basu14}
Basu, S., Racoceanu, D.: Reconstructing neuronal morphology from microscopy
  stacks using fast marching.
\newblock In: IEEE International Conference on Image Processing, pp. 3597--3601
  (2014)

\bibitem{Bellman}
Bellman, R.: Dynamic Programming.
\newblock Princeton University Press, Princeton, NJ (1957)

\bibitem{Capozzoli13}
Capozzoli, A., Curcio, C., Liseno, A., Savarese, S.: A comparison of fast
  marching, fast sweeping and fast iterative methods for the solution of the
  eikonal equation.
\newblock In: 21st Telecommunications Forum, pp. 685--688 (2013)

\bibitem{Chacon14}
Chacon, A.: Eikonal euqations: New two-scale algorithms and error analysis.
\newblock Ph.D. thesis, Cornell Univeristy (2014)

\bibitem{Clement07}
Clement, P., Pailhas, Y., Patron, P., Petillot, Y., Evans, J., Lane, D.: Path
  planning for autonomous underwater vehicles.
\newblock IEEE Transactions on Robotics \textbf{23}(2), 331--341 (2007)

\bibitem{Detrixhe13}
Detrixhe, M., Gibou, F., Min, C.: A parallel fast sweeping method for the
  eikonal equation.
\newblock Journal of Computational Physics \textbf{237}, 46--55 (2013)

\bibitem{Dijkstra}
Dijkstra, E.: {A note on two problems in connection with graphs}.
\newblock Numerische Mathematik \textbf{1}, 269--271 (1959)

\bibitem{Do14}
Do, Q., Mita, S., Yoneda, K.: Narrow passage path planning using fast marching
  method and support vector machine.
\newblock In: IEEE Intelligent Vehicles Symposium Proceedings, pp. 630--635
  (2014)

\bibitem{Forcadel08}
Forcadel, N., Le~Guyader, C., Gout, C.: Generalized fast marching method:
  applications to image segmentation.
\newblock Numerical Algorithms \textbf{48}(1--3), 189--211 (2008)

\bibitem{Fredman87}
Fredman, M., Tarjan, R.: {Fibonacci heaps and their uses in improved network
  optimization algorithms}.
\newblock Journal of the Association for Computing Machinery \textbf{34}(3),
  596--615 (1987)

\bibitem{Fu11}
Fu, Z., Jeong, W., Pan, Y., Whitaker, R.: A fast iterative method for solving
  the eikonal equation on triangulated surfaces.
\newblock SIAM Journal on Scientific Computing \textbf{33}(5), 2468--2488
  (2011)

\bibitem{Gomez12Thesis}
G{\'o}mez, J.V.: {Advanced Applications of the Fast Marching Square Planning
  Method}.
\newblock Master's thesis, Carlos III University (2012)

\bibitem{Gremaud06}
Gremaud, P., Kuster, C.: Computational study of fast methods for the eikonal
  equation.
\newblock SIAM Journal on Scientific Computing \textbf{27}(6), 1803--1816
  (2006)

\bibitem{Porter74}
Jain, A.K., Hong, L., Pankanti, S.: Random insertion into a priority queue
  structure.
\newblock Tech. rep., Stanford Univeristy Reports (1974)

\bibitem{Janson14}
Janson, L., Schmerling, E., Clark, A., Pavone, M.: Fast marching trees: a fast
  marching sampling-based method for optimal motion planning in many
  dimensions.
\newblock International Journal of Robotics Research \textbf{Submitted} (2014)

\bibitem{Jeong08}
Jeong, W., Whitaker, R.: A fast iterative method for eikonal equations.
\newblock SIAM Journal on Scientific Computing \textbf{30}(5), 2512--2534
  (2008)

\bibitem{Jones06}
Jones, M., Baerentzen, J., Sramek, M.: {3D distance fields: a survey of
  techniques and applications}.
\newblock IEEE Transactions on Visualization and Computer Graphics
  \textbf{12}(7), 581--599 (2006)

\bibitem{Kim01}
Kim, S.: An o(n) level set method for eikonal equations.
\newblock SIAM Journal on Scientific Computing \textbf{22}(6), 2178--2193
  (2001)

\bibitem{Sethian98}
Kimmel, R., Sethian, J.A.: {Computing geodesic paths on manifolds}.
\newblock Proceesings of the National Academy of Sciences \textbf{95}(15),
  8431--8435 (1998)

\bibitem{Liu15}
Liu, Y., Bucknall, R.: Path planning algorithm for unmanned surface vehicle
  formations in a practical maritime environment.
\newblock Ocean Engineering \textbf{97}, 126--144 (2015)

\bibitem{Luo14}
Luo, S.: {High-order factorizations and high-order schemes for point-source
  Eikonal equations}.
\newblock SIAM Journal on Numerical Analysis \textbf{52}(1), 23--44 (2014)

\bibitem{Osher88}
Osher, S., Sethian, J.A.: {Fronts propagating with curvature dependent speed:
  Algorithms based on Hamilton-Jacobi formulations}.
\newblock Journal of Computational Physics \textbf{79}(1), 12--49 (1988)

\bibitem{Xinxin14}
Qu, X., Liu, S., Wang, F.: A new ray tracing technique for crosshole radar
  traveltime tomography based on multistencils fast marching method and the
  steepest descend method.
\newblock In: 15th International Conference on Ground Penetrating Radar, pp.
  503--508 (2014)

\bibitem{Rasch08}
Rasch, C., Satzger, T.: Remarks on the o(n) implementation of the fast marching
  method (2008)

\bibitem{Rouy92}
Rouy, E., Tourin, A.: A viscosity solutions approach to shape-from shading.
\newblock SIAM Journal on Numerical Analysis \textbf{29}, 867--884 (1992)

\bibitem{Sethian96}
Sethian, J.A.: {A fast marching level set method for monotonically advancing
  fronts}.
\newblock Proceesings of the National Academy of Sciences \textbf{93}(4),
  1591--1595 (1969)

\bibitem{Sethian99b}
Sethian, J.A.: {Fast Marching Methods}.
\newblock SIAM Review \textbf{41}(2), 199--235 (1999)

\bibitem{Sethian99}
Sethian, J.A.: Level Set Methods and Fast Marching Methods.
\newblock Cambridge University Press (1999)

\bibitem{Sethian00}
Sethian, J.A., Vladimirsky, A.: { Fast Methods for the Eikonal and related
  Hamilton-Jacobi equations on unstructured meshes}.
\newblock Proceesings of the National Academy of Sciences \textbf{97}(11),
  5699--5703 (2000)

\bibitem{Sethian03}
Sethian, J.A., Vladimirsky, A.: {Ordered upwind methods for static
  Hamilton-Jacobi equations: theory and algorithms}.
\newblock SIAM Journal on Numerical Analysis \textbf{41}(1), 325--363 (2003)

\bibitem{Tsai03}
Tsai, Y., Cheng, L., Osher, S., Zhao, H.: Fast sweeping algorithms for a class
  of hamilton-jacobi equations.
\newblock SIAM Journal of Numerical Analysis \textbf{41}(2), 659--672 (2003)

\bibitem{Tsitsiklis95}
Tsitsiklis, J.N.: {Fast Marching Methods}.
\newblock IEEE Transactions on Automatic Control \textbf{40}, 1528--1538 (1995)

\bibitem{Valero13}
Valero-G{\'o}mez, A., G{\'o}mez, J.V., Garrido, S., Moreno, L.: {The Path to
  Efficiency: Fast Marching Method for Safer, More Efficient Mobile Robot
  Trajectories}.
\newblock IEEE Robotics and Automation Magazine \textbf{20}(4) (2013)

\bibitem{Yatziv05}
Yatziv, L., Bartesaghi, A., Sapiro, G.: O(n) implementation of the fast
  marching algorithm.
\newblock Journal of Computational Physics \textbf{212}, 393--399 (2005)

\bibitem{Yershov15}
Yershov, D., Otte, M., Frazzolli, E.: Planning for optimal feedback control in
  the volume of free space.
\newblock ArXiv  (2015)

\bibitem{Zhang11}
Zhang, X., Bording, R.: Fast marching method seismic traveltimes with
  reconfigurable field programmable gate arrays.
\newblock Canadian Journal of Exploration Geophysics \textbf{36}(1), 60--68
  (2011)

\bibitem{Zhao05}
Zhao, H.: Fast sweeping method for eikonal equations.
\newblock Mathematics of Computation \textbf{74}, 603--627 (2005)

\end{thebibliography}

\end{document}